\documentclass[a4paper]{article}

\usepackage[T1]{fontenc}
\usepackage{a4}
\usepackage[english]{babel}
\usepackage[dvipsnames,usenames]{color}

\usepackage{bbm}

\usepackage{verbatim,graphics,indentfirst} 
 \usepackage{url}
\usepackage{stmaryrd}
\usepackage[nosumlimits]{amsmath}
 \usepackage{amssymb}
\usepackage{amsthm}

\usepackage{latexsym}
\usepackage{enumerate}
\usepackage{graphicx}
\usepackage[T1]{fontenc}
\usepackage{lmodern}

\usepackage[english]{babel}
\usepackage{caption}
\usepackage{afterpage}

\usepackage{hyperref}
\usepackage{mathrsfs}
\usepackage{setspace}
\usepackage[all]{xy}

\usepackage[table]{xcolor}

\usepackage{tikz}
\usetikzlibrary{shapes.arrows,chains}
\usetikzlibrary {positioning}
\usetikzlibrary{arrows}
\usepackage{tkz-graph}
\usepackage{float}

\theoremstyle{definition}
\newtheorem{definition}{\vspace{1mm}Definition}%[section]

\theoremstyle{plain}% default
\newtheorem{lemma}[definition]{\vspace{1mm}Lemma}
\newtheorem{theorem}[definition]{\vspace{1mm}Theorem}
\newtheorem{corollary}[definition]{\vspace{1mm}Corollary}
\newtheorem{proposition}[definition]{\vspace{1mm}Proposition}

\newtheorem{example}[definition]{\vspace{1mm}Example}

\newcommand{\cat}[1]{\mathbf{#1}}

\newcommand{\pl}   {{\e\!\!\!\!\ou}}
\newcommand{\conn}{{\copyright}}

\newcommand{\val}{v}
\newcommand{\cB}{\mathcal{B}}

\newcommand{\TWO}{\mathbbm{2}}

\newcommand{\nats}                  {{\mathbb N}}

\newcommand{\cL}{\mathcal{L}}

\newcommand{\cM}{\mathcal{M}}
\newcommand{\Mt}{\mathbb{M}}

\newcommand{\KSMt}{\mathbb{KS}}

\newcommand{\der}    						  {\vartriangleright}

\newcommand{\tuple}[1]                         {{\langle #1\rangle}}

\newcommand{\calB} {{\mathcal B}}

\DeclareMathOperator*{\sub}{\mathsf{sub}}
\DeclareMathOperator*{\var}{\mathsf{var}}
\DeclareMathOperator*{\atm}{\mathsf{atm}}
\DeclareMathOperator*{\sing}{\mathsf{sing}}

\DeclareMathOperator*{\head}{\mathsf{hd}}
\DeclareMathOperator*{\ctx}{\mathsf{ctx}}
\DeclareMathOperator*{\skel}{\mathsf{skel}}
\DeclareMathOperator*{\unskel}{\mathsf{unskel}}
\DeclareMathOperator*{\mon}{\mathsf{Mon}}

\DeclareMathOperator*{\posit}{\mathsf{min}}

\DeclareMathOperator*{\mult}{\mathsf{mult}}
\DeclareMathOperator*{\Atm}{\mathsf{Atm}}

\newcommand{\ou}          {\vee}
\newcommand{\e}          {\wedge}

\DeclareMathOperator*{\fin}{\mathsf{fin}}

\newcommand{\Ax}             {\mathsf{Ax}}
\DeclareMathOperator*{\inst}{\mathsf{inst}}

\newcommand{\Val}{\textrm{Val}}
\newcommand{\BVal}{\textrm{BVal}}
\newcommand{\Lind}{\textrm{Lind}}

\newcommand{\botop}{\ensuremath{\bot\mkern-14mu\top}}

\newcommand{\Int}             {\mathsf{Int}}

\newcommand{\MPt} {\mathbb{MP}}

\date{}

\title{Modular many-valued semantics for combined logics\thanks{Research funded by FCT/MCTES through national funds and when applicable co-funded by EU under the project UIDB/50008/2020. 
Accepted article for The Journal of Symbolic Logic on 2023-04-20.  \url{https://doi.org/10.1017/jsl.2023.22}}} 

\author{Carlos Caleiro and S\'ergio Marcelino\\
{\tt \{ccal,smarcel\}@math.tecnico.ulisboa.pt} \\
{SQIG - Instituto de Telecomunica\c c\~oes}\\
{Dep. Matem\'atica - Instituto Superior T\'ecnico}\\
{Universidade de Lisboa, Portugal}}

\begin{document}

\maketitle

%\tableofcontents

\begin{abstract}
We obtain, for the first time, a modular many-valued semantics for combined logics, which is built directly from many-valued semantics for the logics being combined, by means of suitable universal operations over partial non-deterministic logical matrices. Our constructions preserve finite-valuedness in the context of multiple-conclusion logics whereas, unsurprisingly, it may be lost in the context of single-conclusion logics.

Besides illustrating our constructions over a wide range of examples, we also develop concrete applications of our semantic characterizations, namely regarding the semantics of strengthening a given many-valued logic with additional axioms, the study of conditions under which a given logic may be seen as a combination of simpler syntactically defined fragments whose calculi can be obtained independently and put together to form a calculus for the whole logic, and also general conditions for decidability to be preserved by the combination mechanism.
\end{abstract}

\section{Introduction}

Modularly putting together logics, or their fragments, in particular by joining corresponding calculi, while keeping control of the metatheoretic properties induced and, in particular, of the resulting underlying semantics, is the core idea of the general mechanism for combining logics known as fibring~\cite{gabbay1996,GAbFib99,acs:css:ccal:98a,ccal:car:jfr:css:04,wcarnielli:mconiglio:gab:mpg:css:06,ccal:acs:09}. 
Given its fundamental character and abstract formulation, this mechanism is a key ingredient of the general theory of universal logic~\cite{bez:UL:94,Beziau01082011}. Further, due to its compositional nature, a deep understanding of combined logics is crucial for the construction and analysis of complex logics, a subject of ever growing importance in application fields (see~\cite{frocos}, for instance). Given logics $\cL_1$ and $\cL_2$, fibring combines them into the smallest logic $\cL_1\bullet\cL_2$ on the combined language which extends both $\cL_1$ and $\cL_2$. This simple idea, however, is far from well understood, to date, despite the long track of work on the subject.  An interesting running example of the difficulties at hand consists of the combination of the conjunction-only and disjunction-only fragments of classical logic, which does not coincide with its conjunction-disjunction fragment (see, for instance,~\cite{journals/igpl/BeziauC11,softcomp,humberstone}), and which we will also address. 

To date, we have many interesting general results regarding conservativity, decidability, finite model properties, or interpolation, as well as soundness and completeness preservation with respect to different forms of symbolic calculi~\cite{gabbayoverview,ccal:jabr:05,acs:css:zan:01,css:jfr:car:01,Cerro96combiningclassical,Beckert:1998,ccal:car:jfr:css:04,wcarnielli:mconiglio:gab:mpg:css:06,Rasga01062002,acs:css:ccal:98a,zan:acs:css:99,mconiglio:acs:css:10,schechter,igpl,charfinval}
 for combined logics, but there is no generally usable tool support for the obtained logics, due mainly to the absence of a satisfactory semantic counterpart of fibring that naturally relates models of the component logics with models of the combined logic. With the honorable exception of fusion of modal logics~\cite{Thomason1980,Wolter98,gkwz03,fine,kracht}, a very particular case of fibring, the available general semantics for combined logics, so far, are either not constructible from the semantics of the component logics~\cite{zan:acs:css:99,acs:css:zan:01,Rasga01062002} (due to the use of fullness assumptions), or use highly uncommon semantic structures~\cite{journals/logcom/SernadasSRC09a,schechter}. 

For these reasons, general fibred semantics is still an open problem: how to combine, in the general case, the semantics $\cM_1$ (adequate for logic $\cL_1$) and $\cM_2$ (adequate for logic $\cL_2$) into a semantics 
$\cM_1 \star\cM_2$ built directly from $\cM_1$ and $\cM_2$,  providing an adequate semantics for
 $\cL_1\bullet\cL_2$? We have known for some time that this question is far from straightforward. Indeed, when taking logical matrices as models, as is most common, we know that combining two logics, each given by a single finite matrix, can result in a logic that cannot even be given by a finite set of finite matrices, nor by a single matrix, even if infinite~\cite{charfinval,finval}. This fact led us to considering non-deterministic logical matrices (Nmatrices), as introduced by Avron and coauthors~\cite{Avron01062005,avr:zam:surveyNDS,avronBook}, and in~\cite{wollic17} we understood how this expressive gain could solve the problem in a neat way, but just for disjoint combinations, that is, when the logics being combined do not share any connectives.\smallskip

In this paper, finally, we define a simple and usable general semantics for combined propositional-based logics. We do so by further enriching our semantic structures with partiality, and adopting partial non-deterministic logical matrices (PNmatrices), as introduced in~\cite{Baaz2013}. As we shall see, the added expressivity brought by partiality is crucial not just in keeping our semantic structures as compact as possible, but mostly in dealing with shared language. 
Our work has an additional fundamental ingredient. We consider a very enlightening step forward from traditional Tarski-style $\textsf{Set}\times\textsf{Formula}$ single-conclusion consequence relations (see~\cite{Wojcicki88}) to Scott-like $\textsf{Set}\times\textsf{Set}$ multiple-conclusion consequence relations as introduce in~\cite{Sco:CaA:74}. This abstraction really sheds new light into the overall problem, as shall be made clear.

We will thus show that semantics for combined logics can always be obtained in the form of a PNmatrix obtained directly from given PNmatrices characterizing the original logics. Further, the resulting semantics will be finite as long as the given PNmatrices are finite, at least in the setting of multiple-conclusion logics. The connections between single-conclusion and multiple-conclusion consequence relations (see~\cite{SS}) are fundamental, at this point, in understanding how much more demanding it is to obtain a similar result in the single-conclusion setting, where indeed finiteness may be lost.\smallskip

Our approach is anchored on the observation that semantics, be they given by means of matrices, Nmatrices, or PNmatrices, are simply clever ways of defining suitable collections of bivaluations (see~\cite{Bez:SB}). Since a collection of bivaluations characterizes in a unique way a multiple-conclusion logic (see~\cite{SS}), it is relatively simple to express the combination of multiple-conclusion logics in terms of bivaluations. Then, we just need to come up with a matching construction over PNmatrices, generalizing the finiteness-preserving strict-product operation proposed in~\cite{wollic17}. The single-conclusion case is harder, simply because the characterization of a Tarski-style consequence in terms of bivaluations is no longer one-to-one. The situation can be restored, however, when the 
collection of bivaluations is meet-closed, once we consider a simple (though not finiteness-preserving) corresponding operation over PNmatrices. Furthermore, we show that the constructions above enjoy universal properties, consistent with the definition of $\cL_1\bullet\cL_2$ as the least logic that extends $\cL_1$ and $\cL_2$, as advocated in~\cite{acs:css:ccal:98a,ccal:car:jfr:css:04}.
Besides providing a range of meaningful concrete examples, we also explore three concrete applications of the semantic characterizations obtained, which can be seen as relevant contributions in their own right: a semantic characterization of logics obtained by imposing new axioms to a given many-valued logic (with some perhaps surprising consequences, such as a denumerable semantics for intuitionistic propositional logic); an analysis of when and how to split a logic into syntactical fragments whose combination is the original logic, as a method for obtaining a calculus for a given many-valued logic by putting together axiomatizations for simpler syntactic fragments; and some general, preliminary, results on the preservation of decidability when combining logics.\smallskip

The rest of the paper is organized as follows. In Section~\ref{sec2} we recall the necessary notions regarding logics, their syntax, semantics and calculi, and define the relevant notions regarding combined syntax and combined logics, along with examples illustrating the particularites both of single-conclusion versus multiple-conclusion consequence, and of partiality and non-determinism. The section is a bit long, but the reader can browse through faster, and come back for details when necessary. Section~\ref{sect:semcomblog} is devoted to the main contributions of the paper: the semantic characterization of combined logics in terms of bivaluations, and PNmatrices, first in the multiple-conclusion setting and then in the more complex single-conclusion scenario, by means of simple strict-product and $\omega$-power operations on PNmatrices. Section~\ref{secuniv}, then, shows that the characterizations previously obtained enjoy natural universal properties. Section~\ref{secapps} showcases the three mentioned applications of our characterization, including some relevant important contributions by themselves. The paper concludes, in Section~\ref{sec:conc}, with a summary of the results obtained, their implications, and an outlook of further research. 

\section{Logics and their combination}\label{sec2}

We start by introducing our main objects of interest, fixing notation, and setting up the technical framework necessary for studying the semantics of combined logics, as well as of their associated calculi and semantics. 

\subsection{Syntax}

The \emph{syntax} of a \emph{(propositional) logic} is defined, as usual, by means of a \emph{signature}, an indexed family $\Sigma=\{\Sigma^{(n)}\}_{n\in \nats_0}$ of countable sets where each $\Sigma^{(n)}$ contains all allowed $n$-place connectives, and a  denumerable set $P$ of \emph{variables} (which we consider fixed once and for all). As standard, $L_{\Sigma}(P)$ denotes the set of all \emph{formulas}\footnote{In our setting, the set $L_{\Sigma}(P)$ of all formulas is always denumerable. This is not just a (relatively common) choice. This cardinality constraint plays an essential role in some technical results, where it is crucial to avoid the pitfalls of ill-defined \emph{natural extensions} as observed in~\cite{czelak,cintula}.} constructed from the variables in $P$ using the connectives in $\Sigma$. We will use $p,q,r,\dots$  to denote variables, $A,B,C,\dots$  to denote formulas, and $\Gamma,\Delta,\Theta,\dots$  to denote sets of formulas, in all cases, possibly, with annotations.

We use $\var(A)$, $\sub(A)$, $\head(A)$ to denote, respectively, the set of   \emph{variables occurring in $A$}, the set of  \emph{subformulas of $A$}, and the \emph{head constructor of $A$}, given a formula $A\in L_{\Sigma}(P)$. These notations have simple recursive definitions: $\var(p)=\sub(p)=\{p\}$, and $\head(p)=p$ for $p\in P$;  and if $c\in\Sigma^{(n)}$ and $A_1,\dots,A_n\in L_\Sigma(P)$, $\var(c(A_1,\dots,A_n))=\bigcup_{i=1}^n\var(A_i)$, $\sub(c(A_1,\dots,A_n))=\{c(A_1,\dots,A_n)\}\cup\bigcup_{i=1}^n\sub(A_i)$, and 
$\head(c(A_1,\dots,A_n))=c$. We also consider the extensions of the $\var$ and $\sub$ notations to sets of formulas given, for $\Gamma\subseteq L_\Sigma(P)$, by $\var(\Gamma)=\bigcup_{A\in\Gamma}\var(A)$, and \emph{mutatis mutandis} for $\sub$.\smallskip

As we will consider combined logics, with mixed syntax, we need to consider different signatures, as well as relations and operations between signatures. 
Signatures being  families of sets, the usual set-theoretic notions can be smoothly extended to signatures. We will sometimes abuse notation, and confuse  a signature $\Sigma$ with the set $(\biguplus_{n\in\nats_0}\Sigma^{(n)})$ of all its connectives, and write $c\in\Sigma$ when $c$ is some $n$-place connective $c\in\Sigma^{(n)}$. For this reason, the \emph{empty signature}, with no connectives at all, will be simply denoted by $\emptyset$. 

Let $\Sigma,\Sigma_0$ be two signatures. We say that $\Sigma_0$ is a \emph{subsignature} of $\Sigma$, and write $\Sigma_0\subseteq\Sigma$, whenever $\Sigma_0^{(n)}\subseteq{\Sigma}^{(n)}$ for every $n\in\nats_0$. Expectedly, given signatures $\Sigma_1,\Sigma_2$, we can also define their \emph{shared subsignature} 
$\Sigma_1\cap\Sigma_2=\{\Sigma_1^{(n)}\cap{\Sigma_2}^{(n)}\}_{n\in\nats_0}$, their \emph{combined signature} $\Sigma_1\cup\Sigma_2=\{\Sigma_1^{(n)}\cup{\Sigma_2}^{(n)}\}_{n\in\nats_0}$, and their \emph{difference signature} $\Sigma_1\setminus\Sigma_2=\{{\Sigma_1}^{(n)}\setminus\Sigma_2^{(n)}\}_{n\in\nats_0}$. Clearly, $\Sigma_1\cap\Sigma_2$ is the largest subsignature of both $\Sigma_1$ and $\Sigma_2$, and contains the connectives shared by both. When there are no shared connectives we have that $\Sigma_1\cap\Sigma_2=\emptyset$. Analogously, $\Sigma_1\cup\Sigma_2$ is the smallest signature that has both $\Sigma_1$ and $\Sigma_2$ as subsignatures, and features all the connectives from both $\Sigma_1$ and $\Sigma_2$ in a combined signature. Furthermore, $\Sigma_1\setminus\Sigma_2$ is the largest subsignature of $\Sigma_1$ which does not share any connectives with $\Sigma_2$.\smallskip

A \emph{substitution} is a function $\sigma:P\to L_\Sigma(P)$, that of course extends freely to a function $\sigma:L_{\Sigma_0}(P)\to L_\Sigma(P)$ for every $\Sigma_0\subseteq\Sigma$. As usual, we use $A^\sigma$ to the denote the formula that results from $A\in L_{\Sigma_0}(P)$ by uniformly replacing each variable $p\in\var(A)$ by $\sigma(p)$, and $\Gamma^\sigma=\{A^\sigma:A\in\Gamma\}$ for each $\Gamma\subseteq L_\Sigma(P)$.\smallskip

Note that if $\Sigma_0\subseteq\Sigma$ then $L_{\Sigma_0}(P)\subseteq L_{\Sigma}(P)$. Still, $L_{\Sigma_0}(P)$ and $L_{\Sigma}(P)$ are both denumerable. In fact, the pair can be endowed with a very useful bijection capturing the view of an arbitrary $L_{\Sigma}(P)$ formula from the point of view of $\Sigma_0$, the \emph{skeleton} function $\skel_{\Sigma_0}:L_{\Sigma}(P)\to L_{\Sigma_0}(P)$ (or simply $\skel_0$, or even $\skel$), whose underlying idea we borrow from~\cite{charfinval}. Note that, given $A\in L_{\Sigma}(P)$, $\head(A)$ may be in $\Sigma\setminus\Sigma_0$, in which case we dub $A$ a \emph{$\Sigma_0$-monolith} or simply a \emph{monolith}. The idea is simply to replace monoliths by dedicated variables, just renaming the original variables. Let $\mon(\Sigma_0,\Sigma)$ be the set of all monoliths. It is easy to see that $\mon(\Sigma_0,\Sigma)$ is always countable, though it can be finite when $\Sigma\setminus\Sigma_0$ contains nothing but a finite set of $0$-place connectives.
In any case, $\mon(\Sigma_0,\Sigma)\cup P$ is always denumerable, because $P$ is, and thus we can fix a bijection $\eta:\mon(\Sigma_0,\Sigma)\cup P\to P$. The 
$\skel$ bijection is now definable from $\eta$, recursively, by letting $\skel(p)=\eta(p)$ for $p\in P$, and for $c\in \Sigma^{(n)}$ and $A_1,\dots,A_n\in L_\Sigma(P)$, 
$\skel(c(A_1,\dots,A_n))=c(\skel(A_1),\dots,\skel(A_n))$ if $c\in\Sigma_0$, and $\skel(c(A_1,\dots,A_n))=\eta(c(A_1,\dots,A_n))$ if $c\in\Sigma\setminus\Sigma_0$. 

The $\skel$ bijection thus defined can be inverted by means of the substitution $\unskel_{\Sigma_0}:P\to L_{\Sigma}(P)$ (or simply $\unskel_0$, or even $\unskel$) defined by $\unskel(p)=\eta^{-1}(p)$. Note that $\skel(A)^{\unskel}=A$ for every $A\in L_\Sigma(P)$.

Note also that the restriction of $\skel$ to $P$, $\skel:P \to L_{\Sigma_0}(P)$ (with a slight abuse of notation, we will use the same name) is a substitution, and $\skel(A)=A^{\skel}$ for every $A\in L_{\Sigma_0}(P)$.

Henceforth, we will often apply $\skel$ and $\unskel$ not just to formulas, but also directly to sets of formulas, or even binary relations involving formulas and sets of formulas. For that purpose, we let  
$\skel(X)=\{\skel(x):x\in X\}$ if $X$ is a set, and $\skel(Y,Z)=(\skel(Y),\skel(Z))$ if $(Y,Z)$ is a pair, and \emph{mutatis mutandis} for $\unskel$.

\subsection{Consequence relations}

With respect to the very notion of logic, we will not only consider the traditional Tarski-style set-formula  notion of consequence (\emph{single-conclusion}), but also the more general Scott-style set-set notion of consequence (\emph{multiple-conclusion}), which plays an essential role in our results.\smallskip

A \emph{multiple-conclusion logic} is a pair $\tuple{\Sigma,\der}$ where $\Sigma$ is a signature, and $\der$ is a \emph{multiple-conclusion consequence relation} over $L_\Sigma(P)$, that is, $\der\ \subseteq \wp({L_{\Sigma}(P)})\times \wp({L_{\Sigma}(P)})$ is a relation satisfying the properties below, for every $\Gamma,\Delta,\Gamma',\Delta',\Omega\subseteq L_{\Sigma}(P)$ and $\sigma:P\to L_\Sigma(P)$: 
 \begin{itemize}
 \item[(O)] %[overlap] 
 
if $\Gamma\cap \Delta\neq \emptyset$ then $\Gamma\der \Delta$, % (overlap),
 
  \item[(D)] %[dilution] 
    
if $\Gamma\der \Delta$ then $\Gamma\cup\Gamma'\der \Delta\cup\Delta'$, % (dilution),

  \item[(C)] %[cut for sets] 
  
  if $\Gamma\cup\underline{\Omega}\der \overline{\Omega}\cup\Delta$ for each partition\footnote{Here and elsewhere, $\tuple{\underline{\Omega},\overline{\Omega}}$ is partition of $\Omega$ if $\underline{\Omega}\cup\overline{\Omega}=\Omega$ and $\underline{\Omega}\cap\overline{\Omega}=\emptyset$.} $\tuple{\underline{\Omega},\overline{\Omega}}$ of $\Omega$, then $\Gamma\der \Delta$, % (cut for sets).

 \item[(S)] if $\Gamma\der \Delta$ then $\Gamma^\sigma\der \Delta^\sigma$.  % (sub. invariance).

\end{itemize}

Often, the relation also satisfies the following property, for every $\Gamma,\Delta\subseteq L_\Sigma(P)$:

\begin{itemize}
 \item[(F)] %[si] 
if $\Gamma\der \Delta$ then there exist finite sets $\Gamma_{\textrm{fin}}\subseteq \Gamma$ and
$\Delta_{\textrm{fin}}\subseteq \Delta$ such that $\Gamma_{\textrm{fin}}\der \Delta_{\textrm{fin}}$.
\end{itemize}

Property (C) is best known as \emph{cut for sets}.
The other properties are usually known as \emph{overlap} (O), \emph{dilution} (D), \emph{substitution invariance} (S),
and  \emph{compactness} (F)
(see~\cite{Sco:CaA:74,SS,Wojcicki88}). As is well known, a multiple-conclusion logic $\tuple{\Sigma,\der}$ has a \emph{compact version} 
$\tuple{\Sigma,\der_{\textrm{fin}}}$ defined to be the largest compact multiple-conclusion logic such that $\der_{\textrm{fin}}{\subseteq}\der$.

A pair of sets of formulas $\tuple{\Gamma,\Delta}$ is said to be a \emph{theory-pair of  $\tuple{\Sigma,\der}$} (see~\cite{blasio}) when, for every $A\in L_\Sigma(P)$, if 
${\Gamma}\der{\{A\}\cup\Delta}$ then $A\in\Gamma$, and also, if 
${\Gamma\cup\{A\}}\der{\Delta}$ then $A\in\Delta$. It is clear that, given sets $\Gamma,\Delta$, the pair of sets
$\tuple{\Gamma,\Delta}^\der=\tuple{\{A:{\Gamma}\der{\{A\}\cup\Delta}\},\{A:{\Gamma}\cup{\{A\}\der\Delta}\}}$ is the least theory-pair containing $\tuple{\Gamma,\Delta}$.

A theory-pair $\tuple{\Gamma,\Delta}$ is \emph{consistent} if ${\Gamma}\not\der\Delta$ (otherwise, by dilution, we necessarily have $\Gamma=\Delta=L_\Sigma(P)$). A consistent theory-pair is \emph{maximal} if there is no consistent theory-pair that properly contains it, that is, if $\tuple{\Gamma',\Delta'}$ is a consistent theory with $\Gamma\subseteq\Gamma'$ and $\Delta\subseteq\Delta'$ then $\Gamma=\Gamma'$ and $\Delta=\Delta'$. Equivalently, using cut for the set of all formulas,  a consistent theory-pair $\tuple{\Gamma,\Delta}$ is maximal precisely if $\tuple{\Gamma,\Delta}$ is a partition of $L_\Sigma(P)$. This implies, obviously, that a consistent theory-pair can always be extend to a maximal one.

\smallskip

We say that a pair $\tuple{\Sigma,\vdash}$ is a \emph{single-conclusion logic} if $\Sigma$ is a signature and $\vdash$ is a \emph{single-conclusion consequence relation} over $L_\Sigma(P)$, that is, $\vdash\ \subseteq \wp({L_{\Sigma}(P)})\times {L_{\Sigma}(P)}$ is a relation, and $\posit(\tuple{\Sigma,\vdash})$ is a multiple-conclusion logic, where $\posit(\tuple{\Sigma,\vdash})=\tuple{\Sigma,\der_\vdash}$ is defined, for $\Gamma,\Delta\subseteq L_\Sigma(P)$, by $\Gamma\der_\vdash\Delta$ if and only if there exists $A\in \Delta$ such that $\Gamma\vdash A$. It is well known (see~\cite{SS}) that this constitutes an alternative definition of the usual notion of Tarski-style consequence relation, inheriting the usual properties of \emph{reflexivity}, \emph{monotonicity}, \emph{transitivity}, and \emph{structurality} from properties (ODCS), as well as \emph{compactness} from (F), when it holds. Again,  a single-conclusion logic $\tuple{\Sigma,\vdash}$ has a \emph{compact version} 
$\tuple{\Sigma,\vdash_{\textrm{fin}}}$ defined to be the largest compact single-conclusion logic such that $\vdash_{\textrm{fin}}{\subseteq}\vdash$.

More standardly, now, a set of formulas ${\Gamma}$ is said to be a \emph{theory of  $\tuple{\Sigma,\vdash}$}
 when, for every $A\in L_\Sigma(P)$, if  ${\Gamma}\vdash{A}$ then $A\in\Gamma$. 
Given a set $\Gamma$, we know that $\Gamma^\vdash=\{A:{\Gamma}\vdash A\}$ is the least theory that contains $\Gamma$.

 As usual, a theory ${\Gamma}$ is \emph{consistent} if $\Gamma\neq L_\Sigma(P)$. A consistent theory $\Gamma$ is \emph{maximal relatively to $A\notin\Gamma$} if every theory that properly contains $\Gamma$ must also contain $A$, that is, if $\Gamma\cup\{B\}\vdash A$ for every $B\notin\Gamma$. A consistent theory $\Gamma$ is \emph{relatively maximal} if it is maximal relatively to some formula $A\notin\Gamma$. 
According to the usual \emph{Lindenbaum Lemma} (see, for instance,~\cite{Wojcicki88}), if a single-conclusion logic is compact then a consistent theory always has a relatively maximal extension.

 \smallskip

Every multiple-conclusion logic $\tuple{\Sigma,\der}$ has of course a \emph{single-conclusion companion} defined simply by $\sing(\tuple{\Sigma,\der})=\tuple{\Sigma,\der_{\sing}}$ where, for $\Gamma\cup\{A\}\subseteq L_\Sigma(P)$, $\Gamma\der_{\sing} A$ if and only if $\Gamma\der\{A\}$. 

When  $\tuple{\Sigma,\vdash}$ is a single-conclusion logic, it is clear that the single-conclusion companion of 
$\posit(\tuple{\Sigma,\vdash})$ is precisely $\tuple{\Sigma,\vdash}$. There may however be many multiple-conclusion logics whose single-conclusion companion coincides with 
$\tuple{\Sigma,\vdash}$, which we dub as \emph{multiple-conclusion counterparts of $\tuple{\Sigma,\vdash}$}. Indeed, 
among the many possible multiple-conclusion counterparts of $\tuple{\Sigma,\vdash}$, we have that $\posit(\tuple{\Sigma,\vdash})$ is precisely the \emph{minimal}~\cite{SS}. 
Note that whenever $\sing(\tuple{\Sigma,\der})=\tuple{\Sigma,\vdash}$, it is easy to see that $\Gamma$ is a theory of 
$\tuple{\Sigma,\vdash}$ if and only if there exists $\Delta$ such that the pair $\tuple{\Gamma,\Delta}$ is a theory-pair of $\tuple{\Sigma,\der}$.\smallskip

For both types of logics (we will use $\propto$ as a placeholder for either a multiple-conclusion $\der$ or a single-conclusion $\vdash$), it is well-known that the logics with a common signature form a complete lattice (under the inclusion ordering on the consequence relations), as in both cases it is relatively straightforward to check that intersections of consequence relations are consequence relations (see~\cite{Wojcicki88,SS}). These facts make it relatively easy to enrich the signature of a logic. If $\Sigma_0\subseteq\Sigma$ and $\tuple{\Sigma_0,\propto_0}$ is a logic then the \emph{extension of $\tuple{\Sigma_0,\propto_0}$ to $\Sigma$}, denoted by
$\tuple{\Sigma,\propto_0^\Sigma}$, is the least logic (of the same type) with signature $\Sigma$ such that 
${\propto_0}\subseteq{\propto_0^\Sigma}$. 
It is relatively simple to see that $\propto_0^\Sigma=\bigcup_{\sigma:P\to L_\Sigma(P)}\sigma(\propto_0)=\unskel(\propto_0)$. Just note that for each substitution $\sigma:P\to L_\Sigma(P)$ we have that $({\skel}\circ{\sigma}):P\to L_{\Sigma_0}(P)$ is also a substitution, and also $\unskel\circ({\skel}\circ{\sigma})=\sigma$. The fact that $\skel,\unskel$ are bijections make it straightforward to show that $\propto_0^\Sigma$ is indeed a consequence relation of the correct type.

Let $\Gamma,\Delta,\{A\}\subseteq L_\Sigma(P)$. Concretely, in the multiple-conclusion case, we have that ${\Gamma}\der_0^{\Sigma}{\Delta}$ if and only if ${\skel(\Gamma)}\der_0{\skel(\Delta)}$ if and only if there exist $\Gamma_0,\Delta_0\subseteq L_{\Sigma_0}(P)$ and $\sigma:P\to L_\Sigma(P)$ such that $\Gamma_0^\sigma\subseteq\Gamma$, $\Delta_0^\sigma\subseteq\Delta$, and $\Gamma_0\der_0\Delta_0$. 

Analogously, in the single-conclusion case, we have that ${\Gamma}\vdash_0^{\Sigma}{A}$ if and only if ${\skel(\Gamma)}\vdash_0{\skel(A)}$ if and only if there exist $\Gamma_0,\{A_0\}\subseteq L_{\Sigma_0}(P)$ and $\sigma:P\to L_\Sigma(P)$ such that $\Gamma_0^\sigma\subseteq\Gamma$, $A_0^\sigma=A$, and $\Gamma_0\vdash_0\ A_0$.  \smallskip

It is also quite natural to formulate the \emph{combination of logics} (also known as \emph{fibring}) as follows.

\begin{definition}
Let $\textrm{type}\in\{\textrm{single},\textrm{multiple}\}$.

The \emph{combination} of $\textrm{type}$-conclusion logics $\tuple{\Sigma_1,\propto_1}$, $\tuple{\Sigma_2,\propto_2}$, which we denote by 
$\tuple{\Sigma_1,\propto_1}\bullet\tuple{\Sigma_2,\propto_2}$, is the least $\textrm{type}$-conclusion logic $\tuple{\Sigma_1\cup\Sigma_2,\propto_{12}}$ such that $\propto_1,\propto_2{\subseteq}\propto_{12}$. The combination is said to be \emph{disjoint} if $\Sigma_1\cap\Sigma_2=\emptyset$.
\end{definition}

Note that it follows easily that the combination of compact logics is necessarily compact. Indeed, if $\tuple{\Sigma_1,\propto_1}$, $\tuple{\Sigma_2,\propto_2}$ are compact then the least logic $\tuple{\Sigma_1\cup\Sigma_2,\propto_{12}}$ such that $\propto_1,\propto_2{\subseteq}\propto_{12}$ is also the 
least logic such that $\propto_{1,\textrm{fin}},\propto_{2,\textrm{fin}}{\subseteq}\propto_{12}$. Since it is clear that $\propto_{1,\textrm{fin}},\propto_{2,\textrm{fin}}{\subseteq}\propto_{12,\textrm{fin}}$, it follows that $\propto_{12}{=}\propto_{12,\textrm{fin}}$.

\subsection{Calculi}

Logics are often defined by syntactic means, using symbolic calculi. Again, we will consider multiple as well as single-conclusion calculi.\smallskip

A \emph{multiple-conclusion calculus} is pair $\tuple{\Sigma,R}$ where $\Sigma$ is a signature, and $R\subseteq \wp(L_\Sigma(P))\times \wp(L_\Sigma(P))$ is a set of \emph{(schematic) (multiple-conclusion) inference rules}, each rule $\tuple{\Gamma,\Delta}\in R$ being usually represented as $\frac{\;\Gamma\;}{\Delta}$ where $\Gamma$ is the set of \emph{premises} and $\Delta$ the set of \emph{conclusions} of the rule, also represented as $\frac{\;A_1\;\dots\;A_n\;}{\;B_1\;\dots\;B_m\;}$ when $\Gamma=\{A_1,\dots,A_n\}$ and $\Delta=\{B_1,\dots,B_m\}$. We can associate a multiple-conclusion logic $\tuple{\Sigma,\der_R}$  to a given calculus $\tuple{\Sigma,R}$ by means of a suitable tree-shaped notion of derivation (see~\cite{SS,SYNTH,wollic19}). For the purpose of this paper, however, it is sufficient to characterize $\tuple{\Sigma,\der_R}$ as the least multiple-conclusion logic such that ${R}\subseteq{\der_R}$, being compact when the rules in $R$ all have finite sets of premises and conclusions.

A \emph{single-conclusion calculus} is pair $\tuple{\Sigma,R}$ where $\Sigma$ is a signature, and $R\subseteq \wp(L_\Sigma(P))\times L_\Sigma(P)$ is a set of \emph{(schematic) (single-conclusion) inference rules}, each rule $\tuple{\Gamma,A}\in R$ being usually represented as $\frac{\;\Gamma\;}{A}$ where $\Gamma$ is the set of \emph{premises} and $A$ the  \emph{conclusion} of the rule. It is clear that single-conclusion calculi can be rephrased as particular cases of multiple-conclusion calculi, in particular those whose rules all have singleton sets of conclusions, and that the corresponding notion of derivation will now be linear-shaped, and coincide with the usual notion of proof in Hilbert-style calculi, giving rise to an associated single conclusion logic $\tuple{\Sigma,\vdash_R}$. Again, in any case, $\tuple{\Sigma,\vdash_R}$ is the least single-conclusion logic such that ${R} \subseteq {\vdash_R}$, and again it is finitary if all the rules in $R$ have finitely many premises. In that case, if we use $R$ for both a single conclusion-calculus and its singleton-conclusion multiple-conclusion rephrasal, it is straightforward to see that 
$\tuple{\Sigma,\der_R}=\posit(\tuple{\Sigma,\vdash_R})$.

When $\Sigma_0\subseteq\Sigma$ and $\tuple{\Sigma_0,R_0}$ is a calculus, it is immediate that ${\propto_{\tuple{\Sigma_0,R_0}}^{\Sigma}}={\propto_{\tuple{\Sigma,R_0}}}$, that is, the extension of $\tuple{\Sigma_0,\propto_{R_0}}$ to $\Sigma$ is precisely $\tuple{\Sigma,\propto_{R_0}}$.\smallskip 

Combining calculi at this level is quite simple, as has been known for a long time in the single-conclusion case~\cite{me}. The logic associated to joining the rules of two calculi is precisely the combination of the logics defined by each calculi, as we show next.

\begin{proposition}\label{fibringcalculi}
Let $\textrm{type}\in\{\textrm{single},\textrm{multiple}\}$.

If $\tuple{\Sigma_1,R_1},\tuple{\Sigma_2,R_2}$ are $\textrm{type}$-conclusion calculi, $\tuple{\Sigma_1,\propto_{R_1}},\tuple{\Sigma_2,\propto_{R_2}}$ their associated $\textrm{type}$-conclusion logics, then $\tuple{\Sigma_1,\propto_{R_1}}\bullet\tuple{\Sigma_2,\propto_{R_2}}=\tuple{\Sigma_1\cup\Sigma_2,\propto_{R_1\cup R_2}}$.
\end{proposition}
\proof{The reasoning is straightforward. Let $\tuple{\Sigma_1\cup\Sigma_2,\propto}=\tuple{\Sigma_1,\propto_{R_1}}\bullet\tuple{\Sigma_2,\propto_{R_2}}$, which means that $\propto$ is the least consequence with ${\propto_{R_1},\propto_{R_2}}\subseteq{\propto}$.

As ${R_1}\subseteq{\propto_{R_1}}$ and ${R_2}\subseteq{\propto_{R_2}}$ it follows that $R_1\cup R_2\subseteq{\propto}$, and since $\propto_{R_1\cup R_2}$ is the least 
consequence with this property we can conclude that ${\propto_{R_1\cup R_2}}\subseteq{\propto}$.

Conversely, as $R_1,R_2\subseteq R_1\cup R_2$ it  follows that ${\propto_{R_1},\propto_{R_2}}\subseteq{\propto_{R_1\cup R_2}}$, and since $\propto$ is the least consequence with this property we can conclude that ${\propto}\subseteq{\propto_{R_1\cup R_2}}$.\qed}

For simplicity, whenever $\propto{=}\propto_R$ (in either the single or the multiple-conclusion scenarios), we will say that $R$ constitutes an \emph{axiomatization} of the logic $\tuple{\Sigma,\propto}$.

\subsection{Semantics}

Another common way of charactering logics is by semantic means. For their universality we shall consider \emph{logical matrices}, but will consider an extension of the usual notion which also incorporates two less common ingredients: \emph{non-determinism} and \emph{partiality}, following~\cite{Avron01062005,Baaz2013}. These two ingredients play essential roles in our forthcoming results.\smallskip

A \emph{partial non-deterministic matrix (PNmatrix)} is a pair $\tuple{\Sigma,\Mt}$ where $\Sigma$ is a signature, and $\Mt=\tuple{V,D,\cdot_\Mt}$ is such that
$V$ is a
set (of \emph{truth-values}), $D\subseteq V$ is the set of \emph{designated} values, and 
for each $n\in\nats_0$ and $c\in \Sigma^{(n)}$, $\cdot_\Mt$ provides the \emph{truth-table} $c_\Mt:V^n\to  \wp(V)$ of~$c$ in $\Mt$.  When appropriate we shall refer to $\Mt$ as a $\Sigma$-PNmatrix.
When $c_\Mt(x_1,\dots,x_n)\neq\emptyset$ for all $x_1,\dots,x_n\in V$ we say that the truth-table of $c$ in 
$\Mt$ is \emph{total}. When $c_\Mt(x_1,\dots,x_n)$ has at most one element for all $x_1,\dots,x_n\in V$ we say that the truth-table of $c$ in 
$\Mt$ is \emph{deterministic}. 
When all the truth-tables are total, we say that $\tuple{\Sigma,\Mt}$  is also total, or equivalently that it is a \emph{non-deterministic matrix (Nmatrix)}.
Analogously, when all the truth-tables are deterministic, we say that $\tuple{\Sigma,\Mt}$  is also deterministic, or equivalently that it is a \emph{partial matrix (Pmatrix)}.
Finally, if $\tuple{\Sigma,\Mt}$ is total and deterministic we say that it is simply a \emph{matrix} (the usual notion of a logical matrix, up to an isomorphism).

{ 
Given $X\subseteq V$, $\Mt_X=\tuple{X,\cdot_{\Mt_X},D\cap X}$ is the substructure of $\Mt=\tuple{V,D,\cdot_\Mt}$ obtained by restricting it to values in $X$, that is, 
 $\conn_{\Mt_X}(x_1,\ldots,x_k)=\conn_\Mt(x_1,\ldots,x_k)\cap X$
 for $\conn\in \Sigma^{(k)}$ and $x_1,\ldots,x_k\in X$.
 We say that $X\neq\emptyset$ is a \emph{total component} of $\Mt$ whenever $\Mt_X$ is an Nmatrix. 
We denote by $\mathsf{Tot}_\Mt$ the set of total components of $\Mt$.
A truth-value $x\in V$ is said to be \emph{spurious} in $\Mt$ if there is no total component $X\in\mathsf{Tot}_\Mt$ such that $x\in X$.

}

A \emph{$\Mt$-valuation} is a function $v:L_\Sigma(P)\to V$ such that $v(c(A_1,\dots,A_n))\in c_\Mt(v(A_1),\dots,v(A_n))$ for every $n\in\nats_0$, every $n$-place connective $c\in\Sigma$, and every $A_1,\dots,A_n\in L_\Sigma(P)$. We denote the set of all $\Mt$-valuations by $\Val(\Mt)$. { Note that if $x$ is spurious in $\Mt$ and $v\in\Val(\Mt)$ then $x\notin v(L_\Sigma(P))\in\mathsf{Tot}_\Mt$. Consequently, we have that $\Val(\Mt)=\bigcup_{X\in\mathsf{Tot}_\Mt}\Val(\Mt_X)$.}

As is well known, if $\tuple{\Sigma,\Mt}$ is a matrix then every function $f:Q\to V$ with $Q\subseteq P$ can be extended to a $\Mt$-valuation (in an essentially unique way for all formulas $A$ with $\var(A)\subseteq Q$). When $\Mt$ is a Nmatrix, a function $f$ as above can possibly be extended in many different ways. Still, we know from~\cite{Avron01062005} that a function $f:\Gamma\to V$ with $\Gamma\subseteq L_\Sigma(P)$ can be extended to a $\Mt$-valuation provided that $\sub(\Gamma)\subseteq\Gamma$ and 
that $f(c(A_1,\dots,A_n))\in c_\Mt(f(A_1),\dots,f(A_n))$ whenever $c(A_1,\dots,A_n)\in\Gamma$. We dub such a function a \emph{prevaluation} of the PNmatrix $\Mt$.
In case $\Mt$ is not total, in general, one does not even have such a guarantee~\cite{Baaz2013}, unless $f(\Gamma)\subseteq X$ for some { $X\in\mathsf{Tot}_{\Mt}$. 
}
\smallskip

Given a signature $\Sigma$, a \emph{bivaluation} is a function $b:L_\Sigma(P)\to\{0,1\}$. The set of all such bivaluations is denoted by $\BVal(\Sigma)$.
A set of bivaluations $\calB\subseteq\BVal(\Sigma)$ is known to characterize a multiple-conclusion relation ${\der_{\calB}}\subseteq{\wp(L_\Sigma(P))\times\wp(L_\Sigma(P))}$ defined by 
$\Gamma\der_{\calB}\Delta$ when, for every $b\in \calB$, either $0\in b(\Gamma)$ or $1\in b(\Delta)$. Of course, $\calB$ also characterizes the more usual single conclusion relation
${\vdash_{\calB}}\subseteq{\wp(L_\Sigma(P))\times L_\Sigma(P)}$ such that $\Gamma\vdash_{\calB} A$ when $\Gamma\der_{\calB} \{A\}$, i.e., $\tuple{\Sigma,\vdash_\calB}=\sing(\tuple{\Sigma,\der_\calB})$. In both cases, 
$\tuple{\Sigma,\propto_{\calB}}$ is a logic whenever $\calB$ is closed under substitutions, that is, if $b\in \calB$ and $\sigma: P\to L_\Sigma(P)$ then $(b\circ \sigma)\in \calB$. By definition, if ${\calB'}\subseteq {\calB}$ then ${\propto_{\calB}}\subseteq{\propto_{\calB'}}$.

A PNmatrix $\tuple{\Sigma,\Mt}$ with $\Mt=\tuple{V,D,\cdot_\Mt}$ defines a set of induced bivaluations, closed under substitutions, $\BVal(\Mt)=\{t\circ v: v\in \Val(\Mt)\}$, where $t:V\to\{0,1\}$ such that $t(x)=1$ if $x\in D$, and $t(x)=0$ if $x\notin D$, simply captures the distinction between designated and undesignated values.
We will write $\propto_\Mt$ instead of $\propto_{\BVal(\Mt)}$, and say that $\tuple{\Sigma,\propto_\Mt}$ is the $\textrm{type}$-conclusion logic characterized by $\Mt$, for $\textrm{type}\in\{\textrm{single},\textrm{multiple}\}$, which is always compact (finitary) when $V$ is finite. Clearly, $\sing(\tuple{\Sigma,\der_\Mt})=\tuple{\Sigma,\vdash_\Mt}$.

Note that if $\tuple{\Gamma,\Delta}$ is a maximal consistent theory-pair of $\tuple{\Sigma,\der_\Mt}$ then there must exist $b\in\BVal(\Mt)$ such that $b^{-1}(1)=\Gamma$, and consequently also $b^{-1}(0)=\Delta$. In the single-conclusion case, it is well-known that if $\Gamma$ is a relatively maximal theory of $\tuple{\Sigma,\vdash_\Mt}$ then there must also exist $b\in\BVal(\Mt)$ such that $b^{-1}(1)=\Gamma$ (see~\cite{SS}).

\smallskip

When $\Sigma_0\subseteq\Sigma$ and $\calB_0\subseteq\BVal(\Sigma_0)$, it is straightforward to check that 
$\calB_0^\Sigma=\{b\circ\skel:b\in\calB_0\}\subseteq\BVal(\Sigma)$ characterizes the extension to $\Sigma$ of the logic characterized by $\calB_0$.
It is worth noting that $\calB_0^\Sigma$ is still closed under substitutions because $\skel{\circ}\sigma{\circ}\unskel:P\to L_{\Sigma_0}(P)$ is a substitution, and $\skel(A^\sigma)=\skel(A)^{\skel\circ\sigma\circ\unskel}$; consequently, we have that $b\circ\skel\circ\sigma=(b\circ\skel\circ\sigma\circ \unskel)\circ\skel$ and $b\circ\skel\circ\sigma\circ \unskel$ is in $\calB_0$ if $b$ is.
When $\tuple{\Sigma_0,\Mt_0}$ is a PNmatrix it is very easy, making essential use of non-determinism, to define a PNmatrix $\tuple{\Sigma,\Mt_0^\Sigma}$ that characterizes the extension to $\Sigma$ of the logic characterized by $\Mt_0$. If $\Mt_0=\tuple{V_0,D_0,\cdot_{\Mt_0}}$, one  defines $\Mt_0^\Sigma=\tuple{V_0,D_0,\cdot_{\Mt_0}^\Sigma}$ such that $c_{\Mt_0}^\Sigma=c_{\Mt_0}$ if $c\in\Sigma_0$, and  $c_{\Mt_0}^\Sigma(x_1,\dots,x_n)=V_0$ for every $x_1,\dots,x_n\in V_0$ if $c\in\Sigma\setminus\Sigma_0$ is a $n$-place connective. It is straightforward to check that $\Val({\Mt_0^\Sigma})=\{v\circ\skel:v\in\Val({\Mt_0})\}$.\smallskip

Our work in the next Section of the paper is to prove results that enable us to describe the semantics of the combination of logics characterized by (P)(N)matrices, analogous to Proposition~\ref{fibringcalculi}. Before tackling these problems, it is certainly useful to ground all the relevant notions to concrete illustrative examples.

\subsection{Illustrations}\label{subsect:illustrations}

In this subsection we present a few examples illustrating the relevant notions, as introduced earlier, including the new semantic phenomena brought by partiality and non-determinism, the divide between the single and multiple-conclusion settings, and also the idiosyncrasies of combined logics.

\renewcommand{\arraystretch}{1.2}

\begin{example}\label{classical1}(Classical logic, fragments and combinations).\\
%We start by considering 
Let  $\Sigma_\mathsf{cls}$ be the usual connectives $\top,\neg,\wedge,\vee,\to$, that is, $\Sigma_\mathsf{cls}^{(0)}=\{\top\}$, $\Sigma_\mathsf{cls}^{(1)}=\{\neg\}$, 
$\Sigma_\mathsf{cls}^{(2)}=\{\wedge,\vee,\to\}$, and $\Sigma_\mathsf{cls}^{(n)}=\emptyset$ for $n>2$.

For simplicity, given some connectives $\Theta$ of the signature, we will write $\Sigma_\mathsf{cls}^\Theta$ to denote the subsignature of $\Sigma_\mathsf{cls}$ containing only the connectives in $\Theta$.\smallskip

Take, for instance, $\Sigma_0=\Sigma_\mathsf{cls}^{\neg,\to}\subseteq \Sigma_\mathsf{cls}$ and 
$\Sigma'_0=\Sigma_\mathsf{cls}^{\to}\subseteq \Sigma_\mathsf{cls}$. The formula $A=\neg((\neg p)\to(p\wedge(q\to q)))\in L_{\Sigma_\mathsf{cls}}(P)$ is such that its skeleton from the point of view of $\Sigma_0$ is $\skel_{\Sigma_0}(A)=%\skel_0(A)=
\neg((\neg p')\to r)$, if $\skel_{\Sigma_0}(p)=p'$ and for the monolith we have $\skel_{\Sigma_0}(p\wedge(q\to q))=r$. However, from the point of view of $\Sigma'_0$ we have that $\skel_{\Sigma'_0}(A)$ is itself a variable, as $A$ is a monolith.\smallskip

Consider now the Boolean  $\Sigma_\mathsf{cls}$-matrix
$\mathbbm{2}=\tuple{\{0,1\},\{1\},\cdot_\mathbbm{2}}$
defined by the following tables.  

 \begin{center}
% \begin{tabular}{c}
%$\bot_{\TWO}$ \\%&      \\
%\hline
%$0$
%\end{tabular} 
% \quad
 \begin{tabular}{c}
$\top_{\TWO}$ \\%&      \\
\hline
$1$\\
\quad
\end{tabular} 
 \quad
  \begin{tabular}{c | c c}
    %\hline
    &$\neg_\mathbbm{2}$     \\ 
    \hline
    $0$ & $1$   \\ %\hline
     $1$  &$0$    %\hline
   % $-1$ & $0$   \\
 %   \hline
  \end{tabular}
\quad
   \begin{tabular}{c | c c}
    %\hline
    $\land_\mathbbm{2}$&  $0$ & $1$   \\ 
    \hline
    $0$ & $0$ &  $0$ \\ %\hline
     $1$  &$0$  &  $1$  %\hline
   % $-1$ & $0$   \\
 %   \hline
  \end{tabular}
 \quad
  \begin{tabular}{c | c c}
    %\hline
    $\vee_\mathbbm{2}$&  $0$ & $1$   \\ 
    \hline
    $0$ & $0$ &  $1$ \\ %\hline
     $1$  &$1$  &  $1$  %\hline
   % $-1$ & $0$   \\
 %   \hline
  \end{tabular}
  \quad
  \begin{tabular}{c | c c}
    %\hline
    $\to_\mathbbm{2}$&  $0$ & $1$   \\ 
    \hline
    $0$ & $1$ &  $1$ \\ %\hline
     $1$  &$0$  &  $1$  %\hline
   % $-1$ & $0$   \\
 %   \hline
  \end{tabular}
  \end{center} 
  
It is well-known that  this matrix characterizes classical propositional logic, and we will write $\vdash_\mathsf{cls}$ to denote $\vdash_\mathbbm{2}$.
However, we can also consider a multiple-conclusion version of classical logic, and write $\der_\mathsf{cls}$ to denote $\der_\mathbbm{2}$. Of course, we have that 
$\sing(\tuple{\Sigma_\mathsf{cls},\der_\mathbbm{2}}){=}\tuple{\Sigma_\mathsf{cls},\vdash_\mathbbm{2}}$, but in this case 
$\posit(\tuple{\Sigma_\mathsf{cls},\vdash_\mathbbm{2}}){\subsetneq}\tuple{\Sigma_\mathsf{cls},\der_\mathbbm{2}}$. Note that 
$\emptyset\not\vdash_\mathsf{cls}p$ and $\emptyset\not\vdash_\mathsf{cls}\neg p$, but one has $\emptyset\der_\mathsf{cls}p,\neg p$.\smallskip

The multiple-conclusion version of classical logic $\der_\mathsf{cls}$ is known from~\cite{SS} to be alternatively characterized by the following multiple-conclusion calculus.

 $$\qquad\qquad\qquad\frac{\quad}{ \;\top\;}\qquad\qquad\qquad\frac{p \,,\, \neg p}{ }\qquad\qquad\qquad\frac{ }{\;p \,,\, \neg p\;} \qquad %\qquad\qquad \qquad (R_\neg) 
 $$

$$\qquad\qquad\frac{\; p\e q\;}{p}\qquad\qquad\qquad \frac{\; p\e q\;}{q}\qquad\qquad \qquad\frac{\;p\,,\, q\;}{p\e q}%\qquad\qquad (R_\e)
$$

 $$\qquad\qquad\frac{p}{\; p\vee q\;}\qquad\qquad\qquad \frac{q}{\;p \ou q\;}\qquad\qquad\qquad \frac{p\ou q}{\;p\,,\, q\;}%\qquad\qquad (R_\ou)
 $$
 
 $$\qquad\qquad\frac{}{\;p\,,\,p\to q\;}\qquad\qquad\quad \frac{\;p\,,\,p\to q\;}{q}\qquad\qquad\quad \frac{q}{\;p\to q\;}%\qquad \qquad (R_\to)
 $$\smallskip
 
 It is crucial to observe that the axiomatization above is completely modular with respect to syntax, as each rule involves a single connective. 
 Said another way, the axiomatization is obtained by joining axiomatizations of each of its single-connective fragments, or equivalently 
 we have that $$\tuple{\Sigma_\mathsf{cls},\der_\mathsf{cls}}=\tuple{\Sigma_\mathsf{cls}^\top,\der^\top_\mathsf{cls}}\bullet\tuple{\Sigma^\neg_\mathsf{cls},\der^\neg_\mathsf{cls}}\bullet\tuple{\Sigma^\wedge_\mathsf{cls},\der^\wedge_\mathsf{cls}}\bullet\tuple{\Sigma^\vee_\mathsf{cls},\der^\vee_\mathsf{cls}}\bullet\tuple{\Sigma^\to_\mathsf{cls},\der^\to_\mathsf{cls}}$$
 where $\der^c_\mathsf{cls}{=}\der_\mathsf{cls}{\cap}(\wp(L_{\Sigma_\mathsf{cls}^c}(P))\times\wp(L_{\Sigma_\mathsf{cls}^c}(P)))$ for each $c\in\Sigma_\mathsf{cls}$ is the corresponding single-connective fragment of the logic.
 This fact marks a sharp contrast with respect to the single-conclusion setting, which goes well beyond classical logic, and that will play a key role in our developments. For instance, the single-conclusion version of classical logic is known to be characterized by the following single-conclusion calculus. 
 	 %%%%%
			$$
			\frac{}{\;\top\;}\quad
			\frac{}{\;p\to (q \to p)\;}\quad
			\frac{}{\;(p\to (q \to r))\to ((p\to q)\to(p \to r))\;}\quad
			\frac{\;p \qquad p \to q\;}{q}$$
			$$
			\frac{}{\;(p\e q) \to p\;}\qquad
			\frac{}{\;(p\e q) \to q\;}\qquad
			\frac{}{\;(r\to p) \to ((r\to q) \to (r\to(p\e q)))\;}
			$$
			$$
			\frac{}{\;p \to (p\ou q)\;}\qquad
			\frac{}{\;q \to (p\ou q)\;}\qquad
			\frac{}{\;(p \to r) \to ((q \to r) \to ((p\ou q)\to r)))\;}
			$$
			$$
			\frac{}{\;\neg(p\to p) \to q\;}\qquad
			\frac{}{\;(p\to \neg q)\to (q\to \neg p)\;}\qquad 
			\frac{}{\;\neg \neg p\to p\;}$$
			\smallskip

The fact that the single-conclusion axiomatization is not at all modular with respect to syntax, and that most rules refer to more than one connective, is definitely not a coincidence, as shown in~\cite{softcomp}. Indeed, each of the single-connective fragments $\vdash_\mathsf{cls}^c{=}\vdash_\mathsf{cls}{\cap}(\wp(L_{\Sigma_\mathsf{cls}^c}(P))\times L_{\Sigma_\mathsf{cls}^c}(P))$ of $\vdash_\mathsf{cls}$ for each connective $c\in \Sigma_\mathsf{cls}$ can be axiomatized as follows (see~\cite{rautenberg}).\\

 \begin{tabular}{lcl}
  {$\vdash_\mathsf{cls}^\top$} & \textrm{:} & $\frac{}{\;\;\top\;\;}$\\[2mm]
%
%  {$\vdash_\mathsf{cls}^\bot$} & $\frac{\;\;\bot\;\;}{p}$\\[2mm]
%
  {$\vdash_\mathsf{cls}^\neg$} & \textrm{:} & 
    $\frac{p}{\;\;\neg\neg p\;\;}\quad
    \frac{\;\;\neg \neg p\;\;}{p}\quad 
    \frac{\;\;p\quad  \neg p\;\;}{q}$\\[2mm]
  {$\vdash_\mathsf{cls}^\land$} & \textrm{:} & 
    $\frac{\;p\e q\;}{\;\; p\;\;} \quad
    \frac{\;p\e q\;}{\;\; q\;\;}\quad
    \frac{\;\;p\quad  q\;\;}{p\e q}$\\[2mm]
  {$\vdash_\mathsf{cls}^\lor$} & \textrm{:} & 
    $\frac{p}{\;p\ou q\;}\quad
    \frac{\;p\ou p\;}{p}\quad
    \frac{\;p\ou q\;}{q\ou p}\quad
    \frac{\;p\ou (q \ou r)\;}{(p \ou q)\ou r}$\\[2mm]
  {$\vdash_\mathsf{cls}^{\to}$} & \textrm{:} & 
  $\frac{}{\;p\to(q\to p)}\quad
  \frac{}{\;(p\to(q\to r))\to ((p\to q)\to (p\to r))\;}\quad
  \frac{}{\;(((p\to q)\to p)\to p)\;}\quad 
  \frac{\;p\quad p\to q\;}{q}$\\[2mm]
 \end{tabular}\\
 
 It is clear that joining all these rules will yield a logic much weaker than classical logic (see~\cite{softcomp}), that is, 
 $$\tuple{\Sigma_\mathsf{cls},\vdash_\mathsf{cls}}\supsetneq\tuple{\Sigma_\mathsf{cls}^\top,\vdash^\top_\mathsf{cls}}\bullet\tuple{\Sigma^\neg_\mathsf{cls},\vdash^\neg_\mathsf{cls}}\bullet\tuple{\Sigma^\wedge_\mathsf{cls},\vdash^\wedge_\mathsf{cls}}\bullet\tuple{\Sigma^\vee_\mathsf{cls},\vdash^\vee_\mathsf{cls}}\bullet\tuple{\Sigma^\to_\mathsf{cls},\vdash^\to_\mathsf{cls}}.$$
 
A suitable semantics for the resulting combined logic is not obvious.\hfill$\triangle$
\end{example}

The previous example, though very familiar, was useful for illustrating some of the notions and notations used in this paper, and in particular the differences between the single and multiple-conclusion settings, in particular when combining logics. However, even in a two-valued scenario, non-determinism allows for characterizing some interesting non-classical connectives and logics.

\begin{example}\label{platypuses}(Some new two-valued connectives).\\
Consider the signature  $\Sigma$ such that $\Sigma^{(0)}=\{\botop\}$, $\Sigma^{(1)}=\{\square\}$, 
$\Sigma^{(2)}=\{\rightsquigarrow,\pl\}$, and $\Sigma^{(n)}=\emptyset$ for $n>2$, and the Boolean-like $\Sigma$-Nmatrix
$\mathbbm{2'}=\tuple{\{0,1\},\{1\},\cdot_\mathbbm{2'}}$
defined by the following tables.

\begin{center}
\begin{tabular}{c}
$\botop_{\TWO'}$ \\%&      \\
\hline
$0,1$\\
\quad
\end{tabular} \,\,
  \begin{tabular}{c | c c c }
& $\square_{\TWO'}$   \\
\hline
$0$&  $0,1$  \\
$1$ &$1$ 
\end{tabular} 
\,\,
  \begin{tabular}{c | c c c }
$%\mathsf{mp}
 \rightsquigarrow_{\TWO'}$ & $0$ & $1$  \\
\hline
$0$&  $0,1$ & $ 0,1 $ \\
$1$ &$0$ & $0,1$ 
\end{tabular} \,\,
 \begin{tabular}{c | c c c }
$\pl_{\TWO'}$ & $0$ & $1$  \\
\hline
$0$&  $0$ & $0,1$ \\
$1$ &$0,1$ & $1$ 
\end{tabular} 
\end{center}

The  $0$-place connective $\botop$ is known as {botop} in~\cite{mar:09a}. Since the interpretation of $\botop$ is fully non-deterministic the  connective is unrestricted, in the sense that its logic is characterized by the empty calculus, without any rules.

Then, $\square$ is a $1$-place modal-like box operator whose logic is characterized by the usual (global) \emph{necessitation} rule:
$$\frac{p}{\square p}.$$

The $2$-place implication-like connective $\rightsquigarrow$ is characterized solely by the rule of \emph{modus ponens}:
$$\frac{p\,,\, p \rightsquigarrow q}{q}.$$

All rules above serve both the single and multiple-conclusion settings. \smallskip

Lastly, the $2$-place connective $\pl$, was named \emph{platypus} in~\cite{platypus}, and features some properties of classical conjunction mingled with classical disjunction. Its multiple-conclusion logic is characterized by the following two (familiar) rules.

$$\frac{p\,,\,q}{p\pl q}\qquad\qquad\frac{p\pl q}{p\,,\,q}$$

As shown in~\cite{platypus}, the single-conclusion logic of platypus cannot be finitely axiomatized.\hfill$\triangle$
\end{example}

Of course, the added richness provided by non-determinism goes well beyond the two-valued cases seen above.

\begin{example}\label{processor}(Information sources).\\
In~\cite{avr:bn:kon:genFVS}, the authors introduce a logic for modelling the reasoning of a processor which collects information from different classical sources. Each source may provide information that a certain formula of classical logic is true, or false, or no information at all. This situation gives rise to four possible situation as, for each formula, $(t)$ the processor may have information that it is true and no information that it is false, or  $(f)$ information that it is false and no information that it is true, or $(\top)$ have information that it is true and also information that it is false, or $(\bot)$ have no information at all about the formula. This situation is easily understandable if one reads the situations as collecting the available classical truth-values (0,1) for each formula as $$t=\{1\},\quad f=\{0\},\quad \top=\{0,1\},\quad \bot=\emptyset.$$

As originally presented, the logic is defined over the signature $\Sigma_S=\Sigma_\mathsf{cls}^{\wedge,\vee,\neg}$ and characterized by the Nmatrix 
 $\mathbb{S} = \tuple{\{f,\bot,\top,t\},\{\top,t\},\cdot_{\mathbb{S}}}$ defined by the tables below.

 \begin{center} 
  \begin{tabular}{c | c c c c }
    %\hline
    $\wedge_{\mathbb{S}}$&  $f$ & $\bot$ & $\top$ & $t$   \\ 
    \hline
      $f$ & $f$ &  $f$ & $f$ &  $f$  \\  
$\bot$ & $f$ &  $f,\bot$ & $f$ &  $f,\bot$  \\  
$\top$ & $f$ &  $f$ & $\top$ &  $\top$  \\  
     $t$ & $f$ &  $f,\bot$ & $\top$ &  $t,\top$  \\  
  \end{tabular}
  \quad
  \begin{tabular}{c | c c c c }
    %\hline
    $\vee_{\mathbb{S}}$&  $f$ & $\bot$ & $\top$ & $t$   \\ 
    \hline
      $f$ & $f,\top$ &  $t,\bot$ & $\top$ &  $t$  \\  
$\bot$ & $t,\bot$ &  $t,\bot$ & $t$ &  $t$  \\  
$\top$ & $\top$ &  $t$ & $\top$ &  $t$  \\  
     $t$ & $t$ &  $t$ & $t$ &  $t$  \\  
  \end{tabular}
 \quad
  \begin{tabular}{c | c  }
    %\hline
    &$\neg_{\mathbb{S}}$      \\ 
    \hline
    $f$ & $t$    \\  
$\bot$ & $\bot$    \\  
$\top$ & $\top$   \\  
     $t$ & $f$  
  \end{tabular}
  \end{center} 
 
 Clearly, both conjunction and disjunction have non-deterministic interpretations. For instance, note that $t\wedge_{\mathbb{S}}t=\{t,\top\}$. Thus,  a valuation $v\in\Val(\mathbb{S})$ may be such that there exist formulas $A,B$ with $v(A)=v(B)=v(A\wedge B)=t$ and $v(B\wedge A)=\top$.\smallskip

From~\cite{SYNTH,wollic19} we know that both $\der_{\mathbb{S}}$ and $\vdash_{\mathbb{S}}$ are axiomatized by the following 
  rules.
   $$\frac{p\,,\, q}{\;p\e q\;}\quad 
 \frac{\;p\e q\;}{p}\quad \frac{\;p\e q\;}{q} 
 \quad \frac{\neg p}{\;\neg(p\e q)\;}  \quad \frac{\neg q}{\;\neg(p\e q)\;}$$
 %\quad\frac{\neg(p\e p_2)}{\neg p,\neg p_2}\ _{r6}$$ 

   $$\frac{p}{\;p\ou q\;}\quad \frac{q}{\;p\ou q\;}\quad \frac{\;\neg(p\ou q)\;}{\neg p}
   \quad \frac{\;\neg(p\ou q)\;}{\neg q}
  \quad   
   \frac{\neg p\,,\,\neg q}{\;\neg(p\ou q)\;} 
 %  \frac{p\e p_2}{p_2}\ _{.9} 
 %\quad \frac{\neg p}{\neg(p\e p_2)}\ _{.10} 
  $$ 
  
  $$\quad\frac{p}{\;\neg \neg p\;}\quad \frac{\;\neg \neg p\;}{p} $$
 
The fact that this single-conclusion calculus characterizes the multiple-conclusion consequence $\der_{\mathbb{S}}$ is remarkable, implying that 
$\posit(\tuple{\Sigma_S,\vdash_\mathbb{S}}){=}\tuple{\Sigma_S,\der_\mathbb{S}}$.  
Particularly, given any set of formulas $\Gamma$, note that $\Gamma\not\der_\mathbb{S}\emptyset$, as the function such that $v(A)=\top$ for all formulas $A\in L_{\Sigma_S}(P)$, is such that $v\in\Val(\mathbb{S})$.\smallskip

Note that these logics could not possibly be characterized by finite matrices. Let $p\in P$ be a variable, and define $A_0=p$ and $A_{k+1}=p\lor A_k$ for $k\in\nats_0$. Note that $\Gamma=\{A_k:k\in\nats_0\}$ satisfies $\sub(\Gamma)\subseteq\Gamma$, and consider for each $i\in\nats_0$ the prevaluation $f_i:\Gamma\to\{f,\bot,\top,t\}$ of $\mathbb{S}$ defined by $$f_i(A_k)=
\begin{cases}
 f & \mbox{ if }k<i \\
 \top & \mbox{ otherwise. }
\end{cases}
$$

Each $f_i$ thus extends to a valuation $v_i\in\Val(\mathbb{S})$ showing that $A_i\not\vdash_\mathbb{S} A_k$ for $k<i$. Hence, $\vdash_\mathbb{S}$ fails to be locally finite and therefore cannot be characterized by a finite set of finite matrices~\cite{finval}. The same applies to $\der_\mathbb{S}$, simply because $\sing(\tuple{\Sigma_S,\der_\mathbb{S}}){=}\tuple{\Sigma_S,\vdash_\mathbb{S}}$.
\hfill$\triangle$
\end{example}

Besides allowing for a great amount of compression of the number of necessary truth-values, non-determinism has another outstanding property regarding modularity with respect to syntax.

\begin{example}\label{langext}(Language extensions).\\
Above, given a logic $\tuple{\Sigma_0,\propto_0}$ and $\Sigma_0\subseteq\Sigma$, 
we have defined its extension to the larger signature as 
$\tuple{\Sigma,\propto_0^{\Sigma}}$. If $\tuple{\Sigma_0,\propto_0}$ has an associated calculus $\tuple{\Sigma_0,R}$, then it is clear that the extension $\tuple{\Sigma,\propto_0^{\Sigma}}$ is associated to the exact same rules via the calculus $\tuple{\Sigma,R}$.\smallskip

Take, for instance, single-conclusion classical propositional logic $\tuple{\Sigma_\mathsf{cls},\vdash_\mathsf{cls}}$ as defined in Example~\ref{classical1}. Suppose that we wish to consider its extension to a larger signature $\Sigma\supseteq\Sigma_\mathsf{cls}$ containing a new connective $c\in\Sigma^{(2)}$. We know that the exact same calculus presented above is associated to the extended logic 
$\tuple{\Sigma,\vdash^\Sigma_\mathsf{cls}}$. What is more, we also have that the extended logic is characterized by the $\Sigma$-Nmatrix 
$\mathbbm{2}^\Sigma$ which extends the $\Sigma_\mathsf{cls}$-matrix $\mathbbm{2}$ by letting 
 \begin{center}
   \begin{tabular}{c | c c}
    %\hline
    $c_{\mathbbm{2}^\Sigma}$&  $0$ & $1$   \\ 
    \hline
    $0$ & $0,1$ &  $0,1$ \\ %\hline
     $1$  &$0,1$  &  $0,1$  %\hline
   % $-1$ & $0$   \\
 %   \hline
  \end{tabular}
  \end{center}
where the unrestricted connective is interpreted in a fully non-deterministic manner. 
This extension could not be characterized by a finite matrix~\cite{Avron01062005}.
\hfill$\triangle$
\end{example}

At last, we should illustrate the advantages and intricacies that result from introducing partiality.

\begin{example}\label{kleene1}(Kleene's strong three-valued logic).\\
We consider the (single-conclusion) implication-free fragment of Kleene's strong three-valued logic as defined in~\cite{Belfont}. This logic is defined over the signature $\Sigma=\Sigma_\mathsf{cls}^{\wedge,\vee,\neg}$ and is usually characterized by means of two three-valued $\Sigma$-matrices. Equivalently, the logic is given just by the four-valued Pmatrix $\KSMt=\tuple{\{0,a,b,1\},\{b,1\},\cdot_\KSMt}$ defined by the following truth-tables. 
\setlength{\tabcolsep}{5pt}
  
 \begin{center}
          \begin{tabular}{c | c c c c}
    %\hline
   $\e_\KSMt$ &  $0$ & $a$ & $b$ & $1$   \\ 
    \hline
    $0$ & $0$ & $0$ & $0$& $0$ \\ %\hline
     $a$  &$0$ & $a$ & $\emptyset$ & $a$\\
     $b$ &$0$ & $\emptyset$ & $b$ & $b$\\
     $1$  &$0$ & $a$ & $b$ & $1$    
     %\hline
   % $-1$ & $0$   \\
 %   \hline
  \end{tabular}	
  \quad
       \begin{tabular}{c | c c c c}
    %\hline
   $\ou_\KSMt$ &  $0$ & $a$ & $b$ & $1$   \\ 
    \hline
    $0$ & $0$ & $a$ & $b$& $1$ \\ %\hline
     $a$  &$a$ & $a$ & $\emptyset$ & $1$\\
     $b$ &$b$ & $\emptyset$ & $b$ & $1$\\
     $1$  &$1$ & $1$ & $1$ & $1$    
     %\hline
   % $-1$ & $0$   \\
 %   \hline
  \end{tabular}
  \quad
    \begin{tabular}{c | c }
    %\hline
    & $\neg_\KSMt$   \\ 
    \hline
    $0$ & $1$  \\
     $a$  &$a$ \\
     $b$ &$b$ \\
     $1$  &$0$
     %\hline
   % $-1$ & $0$   \\
 %   \hline
  \end{tabular}	
 \end{center}
%with $D=\{1,b\}$

{Note that several entries of the tables are empty, namely $a\wedge_{\KSMt}b=b\wedge_{\KSMt}a=a\vee_{\KSMt}b=b\vee_{\KSMt}a=\emptyset$. As such, a valuation in $\Val(\KSMt)$ cannot both use the values $a$ and $b$. 
Further, a valuation must either use both $0,1$, or none, because $\neg_\KSMt(0)=1$ and $\neg_\KSMt(1)=0$.
As a consequence, we have
$\mathsf{Tot}_{\KSMt}=\{\{a\},\{b\},\{0,1\},\{0,a,1\},\{0,b,1\}\}.$ The two three-valued matrices mentioned above correspond to 
the restrictions $\KSMt_X$ to the maximal total components $X=\{0,a,1\}$ and $X=\{0,b,1\}$.}\smallskip

The single-conclusion logic $\tuple{\Sigma,\vdash_{\KSMt}}$ is known to be associated with the calculus

    $$\frac{\;p\e q\;}{\;\; p\;\;} \quad
    \frac{\;p\e q\;}{\;\; q\;\;}\quad
    \frac{\;\;p\quad  q\;\;}{p\e q}$$
    $$\frac{p}{\;p\ou q\;}\quad
    \frac{\;p\ou p\;}{p}\quad
    \frac{\;p\ou q\;}{q\ou p}\quad
    \frac{\;p\ou (q \ou r)\;}{(p \ou q)\ou r}$$ 
    
    $$\frac{p\lor (q\land r)}{(p\lor q)\land (p\lor r)}\qquad \frac{(p\lor q)\land (p\lor r)}{p\lor (q\land r)}$$
   
   $$\frac{p \lor r}{\neg \neg p \lor r}\qquad \frac{\neg \neg p \lor r}{p \lor r} \qquad \frac{\neg (p \lor q) \lor r}{ (\neg p \land \neg q) \lor r}$$
   
    $$\frac{(\neg p \land \neg q) \lor r}{\neg (p \lor q) \lor r}\qquad\frac{\neg (p \land q) \lor r}{(\neg p \lor \neg q) \lor r}\qquad\frac{(\neg p \lor \neg q) \lor r}{\neg (p \land q) \lor r} \qquad \frac{(p\land \neg p)\lor r}{(q\lor \neg q) \lor r}$$
     
The Pmatrix ${\KSMt}$ also characterizes the multiple-conclusion logic $\tuple{\Sigma,\der_{\KSMt}}$, which is known from~\cite{wollic19} to be associated with the calculus

$$\frac{p\,,\, q}{\;p\e q\;}\quad 
 \frac{\;p\e q\;}{p}\quad \frac{\;p\e q\;}{q}
 \quad \frac{\neg p}{\;\neg(p\e q)\;}  \quad \frac{\neg q}{\;\neg(p\e q)\;}
 \quad  \frac{\;\neg(p\e q)\;}{\;\neg p\,,\, \neg q\;}$$
 %\quad\frac{\neg(p\e p_2)}{\neg p,\neg p_2}\ _{r6}$$ 

   $$\frac{p}{\;p\ou q\;}\quad \frac{q}{\;p\ou q\;}\quad 
   \frac{\;\neg(p\ou q)\;}{\neg p}
   \quad \frac{\;\neg(p\ou q)\;}{\neg q}
  \quad   
   \frac{\neg p\,,\,\neg q}{\;\neg(p\ou q)\;} \quad
   \frac{\;p\ou q\;}{\;p\,,\, q\;}
 %  \frac{p\e p_2}{p_2}\ _{.9} 
 %\quad \frac{\neg p}{\neg(p\e p_2)}\ _{.10} 
  $$ 
  
  $$\quad\frac{p}{\;\neg \neg p\;}\quad \frac{\;\neg \neg p\;}{p}\quad \frac{\;p\,,\, \neg p\;}{q\,,\,\neg q} $$

This latter multiple-conclusion calculus is clearly more natural. As it contains genuinely multiple-conclusion rules, in this case, 
$\posit(\tuple{\Sigma,\vdash_\KSMt}){\subsetneq}\tuple{\Sigma,\der_\KSMt}$. Note, for instance, that $p\vee q\der_\KSMt p,q$, but that 
$p\vee q\not\vdash_\KSMt p$ and $p\vee q\not\vdash_\KSMt q$.\hfill$\triangle$
   
\end{example}

Partiality can indeed be used in almost all cases to provide a single PNmatrix for a logic characterized by a set of Nmatrices, as we will see later.

\section{Semantics for combined logics}\label{sect:semcomblog}

In this section, we are seeking semantic characterizations for combined logics, using bivaluations, and then PNmatrices. The multiple-conclusion abstraction, which we will analyze first, is important for its purity with respect to combination. As we will show, the step toward the single-conclusion case is then a matter of controlling, semantically, the relationship with the multiple-conclusion companions. 
In both cases, non-determinism and partiality play fundamental roles, as we know that combining two logics, each given by a single finite matrix, can result in a logic that cannot even be characterized by a single matrix~\cite{finval}.

\subsection{The multiple-conclusion case}\label{sec:withstrictproduct}

In the multiple-conclusion case, the characterization we need is really very simple, if one considers bivaluations.
The result stems directly from the known fact (see~\cite{SS}) that for every multiple-conclusion logic $\tuple{\Sigma,\der}$ the set $\calB=\{b\in\BVal(\Sigma):b^{-1}(1)\not\der b^{-1}(0)\}$ is the only one set of bivaluations such that ${\der}={\der_{\calB}}$.
Note that this also implies that ${\der_{\calB}}\subseteq{\der_{\calB'}}$ if and only if ${\calB'}\subseteq{\calB}$.

We fix signatures $\Sigma_1,\Sigma_2$, and sets of bivaluations $\calB_1\subseteq\BVal(\Sigma_1)$ and $\calB_2\subseteq\BVal(\Sigma_2)$, both closed under substitutions.

\begin{proposition}\label{bivmultiple}
We have  that
$\tuple{\Sigma_1,\der_{\calB_1}}\bullet\tuple{\Sigma_2,\der_{\calB_2}}=\tuple{\Sigma_1\cup\Sigma_2,\der_{\calB^{\mult}_{12}}}$ with $\calB^{\mult}_{12}={\calB_1^{\Sigma_1\cup\Sigma_2}\cap \calB_2^{\Sigma_1\cup\Sigma_2}}$.
\end{proposition}
\proof{Let $\tuple{\Sigma_1\cup\Sigma_2,\der_{\calB}}=\tuple{\Sigma_1,\der_{\calB_1}}\bullet\tuple{\Sigma_2,\der_{\calB_2}}$, meaning that $\der_{\calB}$ is the least consequence with ${\der_{\calB_1},\der_{\calB_2}}\subseteq{\der_{\calB}}$.

As $\calB^{\mult}_{12}\subseteq{\calB_1^{\Sigma_1\cup\Sigma_2},\calB_2^{\Sigma_1\cup\Sigma_2}}$ it follows that ${\der_{\calB_1}}\subseteq{\der_{\calB_1^{\Sigma_1\cup\Sigma_2}}}\subseteq{\der_{\calB^{\mult}_{12}}}$ and ${\der_{\calB_2}}\subseteq{\der_{\calB_2^{\Sigma_1\cup\Sigma_2}}}\subseteq{\der_{\calB^{\mult}_{12}}}$. Since $\der_{\calB}$ is the least 
consequence with this property we can conclude that ${\der_{\calB}}\subseteq{\der_{\calB^{\mult}_{12}}}$, and therefore ${\calB^{\mult}_{12}}\subseteq{\calB}$.

Conversely, as ${\der_{\calB_1}}\subseteq{\der_{\calB_1^{\Sigma_1\cup\Sigma_2}}}\subseteq{\der_{\calB}}$ and ${\der_{\calB_2}}\subseteq{\der_{\calB_2^{\Sigma_1\cup\Sigma_2}}}\subseteq{\der_{\calB}}$ it  follows that ${{\calB}}\subseteq{\calB_1^{\Sigma_1\cup\Sigma_2}}$ and ${{\calB}}\subseteq{\calB_2^{\Sigma_1\cup\Sigma_2}}$, and so ${{\calB}}\subseteq{\calB^{\mult}_{12}}$.\qed}\\

This result allows us to obtain a clean abstract characterization of the combination of multiple-conclusion logics. One just needs to note that, given a signature $\Sigma$, a partition $\tuple{\underline{\Omega},\overline{\Omega}}$ of $L_\Sigma(P)$ is essentially the same thing as a bivaluation $b\in\BVal(\Sigma)$ with $b^{-1}(1)=\underline{\Omega}$ and $b^{-1}(0)=\overline{\Omega}$.

\begin{corollary}\label{multiplechar}
Let $\tuple{\Sigma_1,\der_1}$, $\tuple{\Sigma_2,\der_2}$ be multiple-conclusion logics, and consider their combination $\tuple{\Sigma_1\cup\Sigma_2,\der_{12}}=\tuple{\Sigma_1,\der_1}\bullet\tuple{\Sigma_2,\der_2}$. For every $\Gamma,\Delta\subseteq L_{\Sigma_1\cup\Sigma_2}(P)$, we have:
$$\Gamma\der_{12}\Delta$$
$$\textrm{if and only if}$$ 
$$\textrm{for each partition }\tuple{\underline{\Omega},\overline{\Omega}}\textrm{ of }L_{\Sigma_1\cup\Sigma_2}(P)\textrm{, there is }k\in\{1,2\}\textrm{ such that }$$
$$
\Gamma\cup\underline{\Omega}\der_k^{\Sigma_1\cup\Sigma_2}\overline{\Omega}\cup\Delta.$$
\end{corollary}
\proof{Using Proposition~\ref{bivmultiple}, and letting ${\der_1}={\der_{\calB_1}}$ and ${\der_2}={\der_{\calB_2}}$, we have
$\Gamma\not\der_{12}\Delta$ if and only if
there exists $b\in\calB^{\mult}_{12}$ such that $b(\Gamma)\subseteq\{1\}$ and $b(\Delta)\subseteq\{0\}$ if and only if
there exists $b\in\calB_{1}^{\Sigma_1\cup\Sigma_2},\calB_{2}^{\Sigma_1\cup\Sigma_2}$ such that $b(\Gamma)\subseteq\{1\}$ and $b(\Delta)\subseteq\{0\}$ if and only if
there is a partition $\tuple{\underline{\Omega},\overline{\Omega}}$ of $L_{\Sigma_1\cup\Sigma_2}(P)$ such that $\Gamma\cup\underline{\Omega}\not\der_1^{\Sigma_1\cup\Sigma_2}\overline{\Omega}\cup\Delta$ and $\Gamma\cup\underline{\Omega}\not\der_2^{\Sigma_1\cup\Sigma_2}\overline{\Omega}\cup\Delta$.
\qed}\\

The question is now whether we can mimic this simplicity at the level of PNmatrices. 
It turns out that one can define a
very simple but powerful operation (already studied in~\cite{wollic17} with respect to Nmatrices, but only in the disjoint case) in order to combine PNmatrices.

We fix PNmatrices $\tuple{\Sigma_1,\Mt_1}$ and $\tuple{\Sigma_2,\Mt_2}$, with $\Mt_1=\tuple{V_1,D_1,\cdot_{\Mt_1}}$ and $\Mt_2=\tuple{V_2,D_2,\cdot_{\Mt_2}}$. 

\begin{definition}\label{strictproduct}
The \emph{strict product} of $\tuple{\Sigma_1,\Mt_1}$ and $\tuple{\Sigma_2,\Mt_2}$ is the PNmatrix $\tuple{\Sigma_1\cup\Sigma_2,\Mt_1\ast \Mt_2}$ such that $\Mt_1\ast \Mt_2=\tuple{V_\ast,D_\ast,\cdot_{\ast}}$ where 
\begin{itemize}
\item $V_\ast=\{(x,y)\in V_1\times V_2: x\in D_1\textrm{ if and only if }y\in D_2\}$, 
\item $D_\ast=D_1\times D_2$, and 
\item for every $n\in\nats_0$, $c\in(\Sigma_1\cup\Sigma_2)^{(n)}$ and $(x_1,y_1),\dots,(x_n,y_n)\in V_\ast$, we define $c_\ast((x_1,y_1),\dots,(x_n,y_n))\subseteq V_\ast$ by letting $(x,y)\in c_\ast((x_1,y_1),\dots,(x_n,y_n))$ if and only if the following conditions hold:

\begin{itemize}
\item if $c\in\Sigma_1$ then $x\in c_{\Mt_1}(x_1,\dots,x_n)$, and
\item if $c\in\Sigma_2$ then $y\in c_{\Mt_2}(y_1,\dots,y_n)$.
\end{itemize}
\end{itemize}
\end{definition}

Note that $V_\ast$ contains all pairs of truth-values which are \emph{compatible}, that is, either both designated, or both undesignated. The pairs where both values are designated constitute $D_\ast$. The truth-table of a shared connective $c\in\Sigma_1\cap\Sigma_2$ comprises all the pairs of values in the truth-tables of $\Mt_1,\Mt_2$ which are compatible. The truth-table of a non-shared connective, say $c\in\Sigma_1\setminus\Sigma_2$, has each possible value in the truth-table of $\Mt_1$ paired with all compatible values in $V_2$. Clearly, the resulting PNmatrix is fundamentally non-deterministic for non-shared connectives whenever the given PNmatrices have several designated values, and several undesignated values.  For shared connectives, non-determinism may appear if inherited from some of the given PNmatrices. Partiality, on its turn, is crucial to the interpretation of shared connectives, showing up whenever the values given by the truth-tables of $\Mt_1,\Mt_2$ cannot be paired compatibly. For non-shared connectives, partiality may still appear if inherited from the given PNmatrices, or in pathological cases where the given PNmatrices do not have designated values, or do not have undesignated values.\smallskip

Before characterizing the exact scope of this construction, we should note that valuations in $\Mt_1\ast\Mt_2$ are always suitable combinations of valuations in $\Mt_1$ and $\Mt_2$. For $k\in\{1,2\}$, we will use $\pi_k:V_\ast\to V_k$ to denote the obvious projection functions, i.e., $\pi_1(x,y)=x$ and $\pi_2(x,y)=y$. 
Easily, if $v\in\Val({\Mt_1\ast\Mt_2})$ then $(\pi_k\circ v)\in\Val({\Mt_k}^{\Sigma_1\cup\Sigma_2})$. Also, it is clear that $v(A)$, $(\pi_1\circ v)(A)$, $(\pi_2\circ v)(A)$ are all compatible for every formula $A$ in the combined language. Further, if $v_1\in\Val(\Mt_1^{\Sigma_1\cup\Sigma_2})$, $v_2\in\Val(\Mt_2^{\Sigma_1\cup\Sigma_2})$, and $v_1(A)$ is compatible with $v_2(A)$ for every $A\in L_{\Sigma_1\cup\Sigma_2}(P)$, then $(v_1\ast v_2)\in\Val(\Mt_1\ast\Mt_2)$ with $(v_1\ast v_2)(A)=(v_1(A),v_2(A))$ for each $A$. These properties apply also to prevaluations over any set of formulas closed under taking subformulas.\smallskip

\begin{lemma}\label{strictvals}
$\BVal(\Mt_1\ast\Mt_2)=\BVal(\Mt_1^{\Sigma_1\cup\Sigma_2})\cap\BVal(\Mt_2^{\Sigma_1\cup\Sigma_2})$.
\end{lemma}
\proof{Pick $k\in\{1,2\}$. If $v\in \Val(M_1\ast M_2)$ then $({\pi_k}\circ {v})\in\Val(M_k^{\Sigma_1\cup\Sigma_2})$. By definition of $M_1\ast M_2$, it follows that 
$({\pi_k}\circ {v})$ and 
$v$ are compatible for every formula in $L_{\Sigma_1\cup\Sigma_2}(P)$, and thus induce the same bivaluation. We then have that $\BVal(M_1\ast M_2)\subseteq\BVal(\Mt_1^{\Sigma_1\cup\Sigma_2})\cap\BVal(\Mt_2^{\Sigma_1\cup\Sigma_2})$.

Conversely, given valuations $v_1\in\Val(\Mt_1^{\Sigma_1\cup\Sigma_2})$ and $v_2\in\Val(\Mt_2^{\Sigma_1\cup\Sigma_2})$ inducing the same bivaluation $b$, it turns out that they must be compatible for all formulas, $v_1(A)\in D_1$ if and only if $v_2(A)\in D_2$. Therefore, $(v_1\ast v_2)\in\Val(\Mt_1\ast \Mt_2)$ is compatible with both, as $(v_1\ast v_2)(A)\in D_\ast$ if and only if $A\in b^{-1}(1)$, and induces the exact same bivaluation $b$. We conclude that ${\BVal(\Mt_1^{\Sigma_1\cup\Sigma_2})\cap\BVal(\Mt_2^{\Sigma_1\cup\Sigma_2})}\subseteq\BVal(\Mt_1\ast\Mt_2)$.
\qed}\\

Now, we can show that the multiple-conclusion logic characterized by a strict product is the corresponding combined logic. 
\begin{theorem}\label{multiplefibring}
The combination of multiple-conclusion logics characterized by PNmatrices is the multiple-conclusion logic characterized by their strict product, that is, 
$\tuple{\Sigma_1,\der_{\Mt_1}}\bullet\tuple{\Sigma_2,\der_{\Mt_2}}=\tuple{\Sigma_1\cup\Sigma_2,\der_{\Mt_1\ast \Mt_2}}$.
\end{theorem}
\proof{The result follows directly from Proposition~\ref{bivmultiple} and Lemma~\ref{strictvals}. 

With $\calB_1=\BVal(\Mt_1)$ and $\calB_2=\BVal(\Mt_2)$, just note that 
$\BVal(\Mt_1\ast\Mt_2)=\BVal(\Mt_1^{\Sigma_1\cup\Sigma_2})\cap\BVal(\Mt_2^{\Sigma_1\cup\Sigma_2})=\calB_1^{\Sigma_1\cup\Sigma_2}\cap\calB_2^{\Sigma_1\cup\Sigma_2}=\calB_{12}^{\mult}$.\qed}\\

This result, besides being based on a simple operation on PNmatrices, has a very useful feature: it provides a finite-valued semantics for the combined logic whenever we are given finite-valued semantics for both given multiple-conclusion logics.\smallskip

Let us look at some examples, starting with the two-valued case.

\begin{example}\label{multi2}(Multiple-conclusion two-valued combinations).\\
Take two signatures $\Sigma_1,\Sigma_2$ and consider two-valued PNmatrices $\tuple{\Sigma_1,\Mt_1}, \tuple{\Sigma_2,\Mt_2}$ with 
$\Mt_1=\tuple{\{0,1\},\{1\},\cdot_{\Mt_1}}$ and $\Mt_2=\tuple{\{0,1\},\{1\},\cdot_{\Mt_2}}$, such as those in Examples~\ref{classical1} or~\ref{platypuses}.

Note that according to Definition~\ref{strictproduct}, $\Mt_1\ast\Mt_2$ is also two-valued, having values $(0,0)$ and $(1,1)$, with the latter designated. In the following discussion, for simplicity, we shall rename the two values to just 0 and 1, respectively, and assume that the strict product is $\Mt_1\ast\Mt_2=\tuple{\{0,1\},\{1\},\cdot_\ast}$.

If $c\in \Sigma_1\cap\Sigma_2$ is a shared connective, and assuming for the sake of exposition that $c$ is a $2$-place connective, we have that $c_\ast$ behaves as depicted below.

\begin{center}
  \begin{tabular}{c | c c c }
$c_{\Mt_1}$ & $0$ & $1$  \\
\hline
$0$&  $X_{00}$ & $ X_{01} $ \\
$1$ &$X_{10}$ & $X_{11}$ 
\end{tabular} %\quad 
$\star$ %\quad
  \begin{tabular}{c | c c c }
$c_{\Mt_2}$ & $0$ & $1$  \\
\hline
$0$&   $Y_{00}$ & $ Y_{01} $ \\
$1$ & $Y_{10}$ & $ Y_{11} $ 
\end{tabular}\,\,
$\approx$ %\quad  
\begin{tabular}{c | c c c }
$c_\ast$ & $0$ & $1$  \\
\hline
$0$&   $X_{00}\cap Y_{00}$ & $X_{01}\cap Y_{01} $ \\
$1$ & $X_{10}\cap Y_{10}$ & $X_{11}\cap Y_{11} $ 
\end{tabular}
\end{center}

This behaviour is absolutely similar for connectives with any number of places, that is, $c_\ast(x_1,\ldots,x_k)=c_{\Mt_1}(x_1,\ldots,x_k)\cap c_{\Mt_2}(x_1,\ldots,x_k)$ for any $k\in\nats_0$ and $c\in(\Sigma_1\cap\Sigma_2)^{(k)}$. Consequently, $c_\ast(x_1,\ldots,x_k)=c_{\Mt_1}(x_1,\ldots,x_k)= c_{\Mt_2}(x_1,\ldots,x_k)$ whenever $c_{\Mt_1}(x_1,\ldots,x_k)= c_{\Mt_2}(x_1,\ldots,x_k)$.

If $c\in \Sigma_1\cup\Sigma_2$ is not a shared connective, and assuming without loss of generality that $c\in\Sigma_1$, $c\notin\Sigma_2$, then in the strict product we have $c_\ast=c_{\Mt_1}$. This implies that one could simply consider the extended PNmatrices $\Mt_1^{\Sigma_1\cup\Sigma_2},\Mt_2^{\Sigma_1\cup\Sigma_2}$ and just use the equation above, as $c_{\Mt_1}(x_1,\ldots,x_k)\cap\{0,1\}=c_{\Mt_1}(x_1,\ldots,x_k)$.

It is clear that the resulting PNmatrix may be genuinely partial in case some of the sets are disjoint, and in general will be a Pmatrix whenever starting with two matrices.\smallskip

Recalling Example~\ref{classical1}, we see that given any two sets of classical connectives $\Theta_1,\Theta_2$, if $\Mt_1$ and $\Mt_2$ are the 
corresponding fragments of the classical matrix $\TWO$, 
it follows that $\Mt_1\ast\Mt_2$ is the fragment of $\TWO$ corresponding to 
$\Theta_1\cup\Theta_2$. According to Theorem~\ref{multiplefibring}, this implies that 
$\tuple{\Sigma_\mathsf{cls}^{\Theta_1\cup\Theta_2},\der_\mathsf{cls}^{\Theta_1\cup\Theta_2}}=\tuple{\Sigma_\mathsf{cls}^{\Theta_1},\der^{\Theta_1}_\mathsf{cls}}\bullet\tuple{\Sigma^{\Theta_2}_\mathsf{cls},\der^{\Theta_2}_\mathsf{cls}}$ and thus we have
$$\tuple{\Sigma_\mathsf{cls},\der_\mathsf{cls}}=\tuple{\Sigma_\mathsf{cls}^\top,\der^\top_\mathsf{cls}}\bullet\tuple{\Sigma^\neg_\mathsf{cls},\der^\neg_\mathsf{cls}}\bullet\tuple{\Sigma^\wedge_\mathsf{cls},\der^\wedge_\mathsf{cls}}\bullet\tuple{\Sigma^\vee_\mathsf{cls},\der^\vee_\mathsf{cls}}\bullet\tuple{\Sigma^\to_\mathsf{cls},\der^\to_\mathsf{cls}}.$$

We will see later that in the single-conclusion case the situation can be dramatically different.\hfill$\triangle$
\end{example}

Let us now consider a richer example.

\begin{example}\label{lk}(Combining the three-valued implications of Kleene and  \L ukasiewicz).\\
Consider the signature $\Sigma$ with $\Sigma^{(2)}=\{\to\}$ and $\Sigma^{{n}}=\emptyset$ for $n\neq 2$. The three-valued implications of Kleene and  \L ukasiewicz are defined by the well known matrices $\mathbb{K}=\tuple{\{0,\frac{1}{2},1\},\{1\},\cdot_\mathbb{K}}$ and 
$\mathbb{L}=\tuple{\{0,\frac{1}{2},1\},\{1\},\cdot_\mathbb{L}}$ defined below.

 \begin{center}
 \begin{tabular}{c | c c c}
$\to_{\mathbb{K}}$& $0$ & $\frac{1}{2}$ & $1$\\
\hline
$0$& $1$ & $1$ & $1$\\
$\frac{1}{2}$ & $\frac{1}{2}$ & $\frac{1}{2}$ & $1$ \\
$1$ & $0$ & $\frac{1}{2}$ & $1$
\end{tabular}
\qquad
 \begin{tabular}{c | c c c}
$\to_{\mathbb{L}}$& $0$ & $\frac{1}{2}$ & $1$\\
\hline
$0$& $1$ & $1$ & $1$\\
$\frac{1}{2}$ & $\frac{1}{2}$ & $1$ & $1$ \\
$1$ & $0$ & $\frac{1}{2}$ & $1$
\end{tabular}  
\end{center}

According to Theorem~\ref{multiplefibring}, we know that $\tuple{\Sigma,\der_\mathbb{K}}\bullet\tuple{\Sigma,\der_\mathbb{L}}$ is characterized by the five-valued $\Sigma$-Pmatrix $\mathbb{K}\ast\mathbb{L}=\tuple{\{00,0\frac{1}{2},\frac{1}{2}0,\frac{1}{2}\frac{1}{2},11\},\{11\},\cdot_\ast}$ given by the truth-table below, where for simplicity we have renamed each truth-value $(x,y)$ to $xy$.
 
 \begin{center}
 \begin{tabular}{c | c c ccc}
$\to_{\ast}$& $00$ & $0\frac{1}{2}$& $\frac{1}{2}0$& $\frac{1}{2}\frac{1}{2}$ & $11$\\
\hline
$00$& $11$ & $11$ & $11$& $11$ & $11$\\
$0\frac{1}{2}$ & $\emptyset$ & $11$ &$\emptyset$& $11$ & $11$\\ 
$\frac{1}{2}0$ & $\emptyset$ & $\emptyset$ & $\emptyset$& $\emptyset$ & $11$\\ 
$\frac{1}{2}\frac{1}{2}$ & $\frac{1}{2}\frac{1}{2}$ & $\emptyset$ & $\frac{1}{2}\frac{1}{2}$& $\emptyset$ & $11$\\ 
$11$ & $00$ & $0\frac{1}{2}$ & $\frac{1}{2}0$& $\frac{1}{2}\frac{1}{2}$ & $11$

\end{tabular}
 \end{center}

Analyzing $\mathbb{K}\ast\mathbb{L}$ it is clear that the values $\frac{1}{2}0$ and $\frac{1}{2}\frac{1}{2}$ are spurious, in the sense that they cannot be used by any valuation. This happens because the two corresponding diagonal entries of the table are empty, that is, $(\frac{1}{2}0\to_\ast\frac{1}{2}0)=(\frac{1}{2}\frac{1}{2}\to_\ast\frac{1}{2}\frac{1}{2})=\emptyset$.

Removing these spurious elements we obtain the following table.
 \begin{center}
 \begin{tabular}{c | c c c}
$\to_{\ast}$& $00$ & $0\frac{1}{2}$ & $11$\\
\hline
$00$& $11$ & $11$ & $11$\\
$0\frac{1}{2}$ & $\emptyset$ & $11$ & $11$ \\
$11$ & $00$ & $0\frac{1}{2}$ & $11$
\end{tabular}
 \end{center}

{ Further, we also have that $(0\frac{1}{2}\to_\ast00)=\emptyset$. Hence, taking into account the fact that $x \to_{\ast} x=11$ we conclude that the total components of the Pmatrix are $\mathsf{Tot}_{\mathbb{K}\ast\mathbb{L}}=\{\{11\},\{00,11\},\{0\frac{1}{2},11\}\}$. This fact implies that $\BVal(\mathbb{K}\ast\mathbb{L})$ is precisely the set of all classical interpretations of classical implication, simply because the restrictions of the Pmatrix to its maximal total components $X=\{00,11\}$ and $X=\{0\frac{1}{2},11\}$ are both isomorphic copies of the two-valued truth-table of classical implication.
We can conclude, thus, that 
$\tuple{\Sigma,\der_\mathbb{K}}\bullet\tuple{\Sigma,\der_\mathbb{L}}=\tuple{\Sigma,\der_{\mathbb{K}\ast\mathbb{L}}}=
\tuple{\Sigma_\mathsf{cls}^\to,\der_\mathsf{cls}^\to}$ coincides with the multiple-conclusion implication-only fragment of classical logic, as defined in Example~\ref{classical1}.}\smallskip

Confirmation of this interesting fact can be obtained by putting together calculi for the logics, according to Proposition~\ref{fibringcalculi}. These multiple-conclusion versions of the logics are not very well known, but a calculus for $\der_\mathbb{K}$ can be readily obtained using the technique of~\cite{SYNTH,wollic19}, whereas a calculus for $\der_\mathbb{L}$ can be found in~\cite{AvronNatural}. Of course, all the rules in these calculi are classically valid, because $\der_\mathbb{K},\der_\mathbb{L}{\subseteq} \der_\mathsf{cls}^\to$ (all classical valuations are permitted in both matrices). Further, as we know the rules of classical implication from Example~\ref{classical1}, it suffices to note that: (1) ${q}\der_\mathbb{K}{p\to q}$, and also ${q}\der_\mathbb{L}{p\to q}$; (2) ${p,p\to q}\der_\mathbb{K}q$, and also ${p,p\to q}\der_\mathbb{L}q$; and (3) 
$\emptyset\der_\mathbb{L}p\to p$ and ${p\to p}\der_\mathbb{K}p,p\to q$ which implies that $\emptyset\der_{\mathbb{K}\ast\mathbb{L}}p,p\to q$.\smallskip

We will return to this example in the more familiar single-conclusion setting.
\hfill $\triangle$
\end{example}

We shall now illustrate how different it is to combine logics with or without shared connectives.

\begin{example}\label{lkdisj}(The three-valued implications of Kleene and  \L ukasiewicz, disjointly).\\
We shall revisited Example~\ref{lk} above, combining the (multiple-conclusion) logics of Kleene's and  \L ukasiewicz's three-valued implications, but now assuming that the connectives are syntactically different, i.e., that Kleene's implication comes from a signature $\Sigma_K$ with $\Sigma_K^{(2)}=\{\to^K\}$ and $\Sigma_K^{{n}}=\emptyset$ for $n\neq 2$, whereas  \L ukasiewicz's implication comes from a disjoint signature $\Sigma_L$ with $\Sigma_L^{(2)}=\{\to^L\}$ and $\Sigma_L^{{n}}=\emptyset$ for $n\neq 2$. The corresponding matrices are now $\tuple{\Sigma_K,\mathbb{K'}}$ with $\mathbb{K'}=\tuple{\{0,\frac{1}{2},1\},\{1\},\cdot_\mathbb{K'}}$ and $\tuple{\Sigma_L,\mathbb{L'}}$ with
$\mathbb{L'}=\tuple{\{0,\frac{1}{2},1\},\{1\},\cdot_\mathbb{L'}}$ defined below.

 \begin{center}
 \begin{tabular}{c | c c c}
$\to^K_{\mathbb{K'}}$& $0$ & $\frac{1}{2}$ & $1$\\
\hline
$0$& $1$ & $1$ & $1$\\
$\frac{1}{2}$ & $\frac{1}{2}$ & $\frac{1}{2}$ & $1$ \\
$1$ & $0$ & $\frac{1}{2}$ & $1$
\end{tabular}
\qquad
 \begin{tabular}{c | c c c}
$\to^L_{\mathbb{L'}}$& $0$ & $\frac{1}{2}$ & $1$\\
\hline
$0$& $1$ & $1$ & $1$\\
$\frac{1}{2}$ & $\frac{1}{2}$ & $1$ & $1$ \\
$1$ & $0$ & $\frac{1}{2}$ & $1$
\end{tabular}  
\end{center}

According to Theorem~\ref{multiplefibring}, we know that $\tuple{\Sigma_K,\der_\mathbb{K'}}\bullet\tuple{\Sigma_L,\der_\mathbb{L'}}$ is characterized by the five-valued $(\Sigma_K\cup\Sigma_L)$-Nmatrix $\mathbb{K'}\ast\mathbb{L'}=\tuple{\{00,0\frac{1}{2},\frac{1}{2}0,\frac{1}{2}\frac{1}{2},11\},\{11\},\cdot_\ast}$ given by the truth-tables below.

 \begin{center}
 \begin{tabular}{c | c c ccc}
$\to^K_{\ast}$& $00$ & $0\frac{1}{2}$& $\frac{1}{2}0$& $\frac{1}{2}\frac{1}{2}$ & $11$\\
\hline
$00$& $11$ & $11$ & $11$& $11$ & $11$\\
$0\frac{1}{2}$ & $11$ & $11$ & $11$& $11$ & $11$\\
$\frac{1}{2}0$ & $\frac{1}{2}0,\frac{1}{2}\frac{1}{2}$ & $\frac{1}{2}0,\frac{1}{2}\frac{1}{2}$ & $\frac{1}{2}0,\frac{1}{2}\frac{1}{2}$& $\frac{1}{2}0,\frac{1}{2}\frac{1}{2}$ & $11$\\ 
$\frac{1}{2}\frac{1}{2}$ & $\frac{1}{2}0,\frac{1}{2}\frac{1}{2}$ & $\frac{1}{2}0,\frac{1}{2}\frac{1}{2}$ & $\frac{1}{2}0,\frac{1}{2}\frac{1}{2}$& $\frac{1}{2}0,\frac{1}{2}\frac{1}{2}$ & $11$\\ 
$11$ & $00,0\frac{1}{2}$ & $00,0\frac{1}{2}$ & $\frac{1}{2}0,\frac{1}{2}\frac{1}{2}$& $\frac{1}{2}0,\frac{1}{2}\frac{1}{2}$ & $11$
\end{tabular}
 \end{center}
 \begin{center}
  \begin{tabular}{c | c c ccc}
$\to^L_{\ast}$& $00$ & $0\frac{1}{2}$& $\frac{1}{2}0$& $\frac{1}{2}\frac{1}{2}$ & $11$\\
\hline
$00$& $11$ & $11$ & $11$& $11$ & $11$\\
$0\frac{1}{2}$ & $0\frac{1}{2},\frac{1}{2}\frac{1}{2}$ & $11$ &$0\frac{1}{2},\frac{1}{2}\frac{1}{2}$& $11$ & $11$\\ 
$\frac{1}{2}0$ & $11$ & $11$ & $11$& $11$ & $11$\\
$\frac{1}{2}\frac{1}{2}$ & $0\frac{1}{2},\frac{1}{2}\frac{1}{2}$ & $11$ & $0\frac{1}{2},\frac{1}{2}\frac{1}{2}$& $11$ & $11$\\ 
$11$ & $00,\frac{1}{2}0$ & $0\frac{1}{2},\frac{1}{2}\frac{1}{2}$ & $00,\frac{1}{2}0$& $0\frac{1}{2},\frac{1}{2}\frac{1}{2}$ & $11$
\end{tabular}
 \end{center}

Of course, this is quite different from simply considering the given three-valued interpretations of each connective.
Note that since all designated entries of $\to^K_{\mathbb{K'}}$ are also designated in $\to^L_{\mathbb{L'}}$ one could perhaps wrongly expect to have that $p\to^Kq\der_{\mathbb{K'}\ast\mathbb{L'}}p\to^Lq$. This assertion can be shown to fail by considering, for instance, any valuation $v\in\Val(\mathbb{K'}\ast\mathbb{L'})$ with $v(p)=0\frac{1}{2}$ and $v(q)=00$, simply because 
$(0\frac{1}{2}\to^K_{\ast}00)=\{11\}$ is necessarily designated but $(0\frac{1}{2}\to^L_{\ast}00)=\{0\frac{1}{2},\frac{1}{2}\frac{1}{2}\}$ contains no designated value.\hfill $\triangle$
\end{example}

\subsection{The single-conclusion case}

The results above really illustrate the simplifying power of using multiple conclusions. Of course, things are not so simple in the single-conclusion scenario. Still, we know enough to be able to take nice conclusions from the same line of reasoning. Let us start with bivaluations. 

If we analyze the proof of Proposition~\ref{bivmultiple}, it is simple to understand why it cannot be simply replicated in the single-conclusion case.
The same line of reasoning would allows us to conclude that ${\vdash_{\calB_1},\vdash_{\calB_2}}\subseteq {\vdash_{\calB_{12}^{\mult}}}$ with $\calB^{\mult}_{12}={\calB_1^{\Sigma_1\cup\Sigma_2}\cap \calB_2^{\Sigma_1\cup\Sigma_2}}$. However, in general, 
$\tuple{\Sigma_1\cup\Sigma_2,\vdash_{\calB_{12}^{\mult}}}$ may very well not be the least such single-conclusion logic. Given $\tuple{\Sigma_1\cup\Sigma_2,\vdash_{\calB}}$ with ${\vdash_{\calB_1},\vdash_{\calB_2}}\subseteq {\vdash_{\calB}}$, we cannot guarantee that $\calB\subseteq{\calB_1^{\Sigma_1\cup\Sigma_2}}$ and $\calB\subseteq{\calB_2^{\Sigma_1\cup\Sigma_2}}$. 
Concretely, assuming that $\Gamma\not\vdash_{\calB} A$ we still know that $\Gamma\not\vdash_{\calB_1}^{\Sigma_1\cup\Sigma_2} A$ and $\Gamma\not\vdash_{\calB_2}^{\Sigma_1\cup\Sigma_2} A$. 
Thus, we still have bivaluations  $b_k\in\calB_k^{\Sigma_1\cup\Sigma_2}$ such that $b_k(\Gamma)\subseteq \{1\}$ and $b_k(A)=0$, for $k\in\{1,2\}$. However, now, despite the fact that $b_1$ and $b_2$ coincide for all formulas in $\Gamma\cup\{A\}$, we have no way of making sure that $b_1=b_2$. 

This problem lies in a crucial difference from the multiple-conclusion case, as the same single-conclusion logic can be characterized by distinct sets of bivaluations. 
Given $\calB\subseteq\BVal(\Sigma)$, let its \emph{meet-closure} be the set 
$\calB^\cap=\{b_X:X\subseteq \calB\}\subseteq\BVal(\Sigma)$ with each meet defined by $b_X(A)=1$ precisely if $b(A)=1$ for every $b\in X$, for each $A\in L_\Sigma(P)$. It is well known that 
${\vdash_{\calB}}\subseteq{\vdash_{\calB'}}$ if and only if ${\calB'}\subseteq{\calB}^\cap$, and thus that two sets of bivaluations characterize the same single-conclusion logic precisely when their meet-closures coincide. Indeed, every single-conclusion logic $\tuple{\Sigma,\vdash}$ is such that ${\vdash}={\vdash_{\calB}}$ provided that ${{\calB^\cap}}=\{b\in\BVal(\Sigma):b^{-1}(1)\textrm{ is a theory of }\tuple{\Sigma,\vdash}\}$, which is the largest set of bivaluations that characterizes $\tuple{\Sigma,\vdash}$.

The trick is then to work with meet-closed sets of bivaluations, that is, $\calB=\calB^\cap$. 
In that case, it is worth noting that $\tuple{\Sigma,\der_{\calB}}=\posit(\tuple{\Sigma,\vdash_{\calB}})$.

\begin{proposition}\label{bivsingle}
We have that  
$\tuple{\Sigma_1,\vdash_{\calB_1}}\bullet\tuple{\Sigma_2,\vdash_{\calB_2}}=\tuple{\Sigma_1\cup\Sigma_2,\vdash_{\calB^{\sing}_{12}}}$ with $\calB^{\sing}_{12}={(\calB_1^{\Sigma_1\cup\Sigma_2})^{\cap}\cap (\calB_2^{\Sigma_1\cup\Sigma_2})^{\cap}}$.
\end{proposition}
\proof{Let $\tuple{\Sigma_1\cup\Sigma_2,\vdash_{\calB}}=\tuple{\Sigma_1,\vdash_{\calB_1}}\bullet\tuple{\Sigma_2,\vdash_{\calB_2}}$, meaning that $\vdash_{\calB}$ is the least consequence with ${\vdash_{\calB_1},\vdash_{\calB_2}}\subseteq{\vdash_{\calB}}$.

As $\calB^{\sing}_{12}\subseteq{(\calB_1^{\Sigma_1\cup\Sigma_2})^\cap}$ and $\calB^{\sing}_{12}\subseteq{(\calB_2^{\Sigma_1\cup\Sigma_2})^\cap}$ it follows that 
${\vdash_{\calB_1}}\subseteq{\vdash_{(\calB_1^{\Sigma_1\cup\Sigma_2})^\cap}}\subseteq{\vdash_{\calB^{\sing}_{12}}}$ and ${\vdash_{\calB_2}}\subseteq{\vdash_{(\calB_2^{\Sigma_1\cup\Sigma_2})^\cap}}\subseteq{\vdash_{\calB^{\sing}_{12}}}$. Since $\vdash_{\calB}$ is the least 
consequence with this property we can conclude that ${\vdash_{\calB}}\subseteq{\vdash_{\calB^{\sing}_{12}}}$, and therefore ${\calB^{\sing}_{12}}\subseteq{\calB}^\cap$.

Conversely, as ${\vdash_{\calB_1}}\subseteq{\vdash_{\calB_1^{\Sigma_1\cup\Sigma_2}}}\subseteq{\vdash_{\calB}}$ and ${\vdash_{\calB_2}}\subseteq{\vdash_{\calB_2^{\Sigma_1\cup\Sigma_2}}}\subseteq{\vdash_{\calB}}$ it  follows that ${{\calB}}\subseteq{(\calB_1^{\Sigma_1\cup\Sigma_2})^\cap}$ and ${{\calB}}\subseteq{(\calB_2^{\Sigma_1\cup\Sigma_2})^\cap}$, and so ${{\calB}}\subseteq{\calB^{\sing}_{12}}$.\qed}\\

Note that $\calB^{\sing}_{12}$ is necessarily meet-closed, as it is the intersection of meet-closed sets. Hence, not only $\posit(\tuple{\Sigma_1,\vdash_{\calB_1}})=\tuple{\Sigma_1,\der_{\calB_1^\cap}}$ and $\posit(\tuple{\Sigma_2,\vdash_{\calB_2}})=\tuple{\Sigma_2,\der_{\calB_2^\cap}}$, but also $\posit(\tuple{\Sigma_1,\vdash_{\calB_1}}\bullet\tuple{\Sigma_2,\vdash_{\calB_2}})=\tuple{\Sigma_1,\der_{\calB_1^\cap}}\bullet\tuple{\Sigma_2,\der_{\calB_2^\cap}}$.
\smallskip

The next abstract characterization of combined single-conclusion logics follows. Just note that given $b\in\calB\subseteq\BVal(\Sigma)$ closed under substitutions, it is always the case that $b^{-1}(1)$ is a theory of $\tuple{\Sigma,\vdash_{\cal B}}$. The converse is in general not true, unless $\calB$ is meet-closed. Concretely, if 
$\Gamma\subseteq L_\Sigma(P)$ is a theory of $\tuple{\Sigma,\vdash_{\cal B}}$ then there always exists $b\in\calB^\cap$ such that $b^{-1}(1)=\Gamma$.

\begin{corollary}\label{singlechar}
Let $\tuple{\Sigma_1,\vdash_1}$, $\tuple{\Sigma_2,\vdash_2}$ be single-conclusion logics, and consider their combination $\tuple{\Sigma_1\cup\Sigma_2,\vdash_{12}}=\tuple{\Sigma_1,\vdash_1}\bullet\tuple{\Sigma_2,\vdash_2}$. For every $\Gamma\cup\{A\}\subseteq L_{\Sigma_1\cup\Sigma_2}(P)$, we have:
$$\Gamma\vdash_{12}A$$
$$\textrm{if and only if}$$
$$\textrm{ for each }\Gamma\subseteq\Omega\subseteq L_{\Sigma_1\cup\Sigma_2}(P),$$
$$\textrm{if }
{\Omega}\textrm{ is  a theory of both }\tuple{\Sigma_1\cup\Sigma_2,\vdash_1^{\Sigma_1\cup\Sigma_2}}\textrm{ and }
\tuple{\Sigma_1\cup\Sigma_2,\vdash_2^{\Sigma_1\cup\Sigma_2}}\textrm{ then }A\in\Omega.$$
\end{corollary}
\proof{Using Proposition~\ref{bivsingle}, and letting ${\vdash_1}={\vdash_{\calB_1}}$ and ${\vdash_2}={\vdash_{\calB_2}}$, we have
$\Gamma\not\vdash_{12}A$ if and only if
there exists $b\in\calB^{\sing}_{12}$ such that $b(\Gamma)\subseteq\{1\}$ and $b(A)=0$ if and only if
there exists $b\in(\calB_{1}^{\Sigma_1\cup\Sigma_2})^\cap\cap(\calB_{2}^{\Sigma_1\cup\Sigma_2})^\cap$ such that $b(\Gamma)\subseteq\{1\}$ and $b(A)=0$ if and only if
there is $\Omega\subseteq L_{\Sigma_1\cup\Sigma_2}(P)$ such that $\Omega$ is a theory of both
$\tuple{\Sigma_1\cup\Sigma_2,\vdash_1^{\Sigma_1\cup\Sigma_2}}$ and
$\tuple{\Sigma_1\cup\Sigma_2,\vdash_2^{\Sigma_1\cup\Sigma_2}}$ with $\Gamma\subseteq\Omega\not\owns A$.
\qed}\\

This result captures neatly the intuition one already had from Proposition~\ref{fibringcalculi}, amounting to the fact that $\Gamma\vdash_{12}A$ precisely when $A$ is obtained by closing 
$\Gamma$ with respect to both $\vdash_1$ and $\vdash_2$.\smallskip

Turning to PNmatrices, expectedly, Theorem~\ref{multiplefibring} does not immediately apply to single-conclusion logics, as shown by the next counterexample.

\begin{example}\label{negand1}(Limitations of the strict product operation).\\
We present an example showing that, contrarily to the multiple-conclusion setting, the single-conclusion logic characterized by the strict product of two PNmatrices may fail to be the combination of the single-conclusion logics characterized by each of the given PNmatrices.\smallskip

Consider the $\neg$-only and $\wedge$-only fragments of classical logic, $\tuple{\Sigma^\neg_\mathsf{cls},\vdash^\neg_\mathsf{cls}}$ and $\tuple{\Sigma^\wedge_\mathsf{cls},\vdash^\wedge_\mathsf{cls}}$ respectively, as defined in Example~\ref{classical1}. As we have seen in Example~\ref{multi2}, the strict product of their two-valued matrices is simply the classical two-valued matrix for both connectives $\neg,\wedge$ which characterizes the $\neg,\wedge$-fragment of classical logic.

However, letting $\tuple{\Sigma^\neg_\mathsf{cls},\vdash^\neg_\mathsf{cls}}\bullet\tuple{\Sigma^\wedge_\mathsf{cls},\vdash^\wedge_\mathsf{cls}}=\tuple{\Sigma^{\neg,\wedge}_\mathsf{cls},\vdash}$, and taking into account the single-conclusion calculi for each of the fragments, shown also in Example~\ref{classical1} and also Proposition~\ref{fibringcalculi}, it is not at all clear that $\vdash{=}\vdash^{\neg,\wedge}_\mathsf{cls}$. As a matter of fact, this is not the case, as $\vdash$ is a strictly weaker consequence, as shown in~\cite{softcomp,wollic17}. For instance, $\neg p\not\vdash\neg(p\wedge p)$. We will show exactly why in a forthcoming example.\smallskip

Clearly, the set of classical bivaluations for conjunction is meet-closed, and therefore the phenomenon above can only be justified by the fact that the set of classical bivaluations for negation is not meet-closed. Concretely, it is clear that there is a single classical (bi)valuation $b_0:L_{\Sigma^\neg_\mathsf{cls}}(P)\to\{0,1\}$ such that $b_0(p)=0$ for all $p\in P$, and thus $b_0(\neg p)=1$ and so on, and also a single classical (bi)valuation $b_1:L_{\Sigma^\neg_\mathsf{cls}}(P)\to\{0,1\}$ such that $b_1(p)=1$ for all $p\in P$. Letting $X=\{b_0,b_1\}$, it is immediate that $b_X(A)=0$ for $A\in L_{\Sigma^\neg_\mathsf{cls}}(P)$, but it is also clear that $b_X$ is not a classical bivaluation.
\hfill $\triangle$
\end{example}

As before, the problem lies with the fact that, in general, the set of bivaluations induced by a PNmatrix does not have to be meet-closed. In order to cope with this possibility, we consider the following property.

\begin{definition}
A PNmatrix $\tuple{\Sigma,\Mt}$ with $\Mt=\tuple{V,D,\cdot_{\Mt}}$ is \emph{saturated} if for each consistent theory $\Gamma$ of $\tuple{\Sigma,\vdash_\Mt}$ there exists $v\in\Val(\Mt)$ such that 
$\Gamma=v^{-1}(D)$.
\end{definition}

If $\tuple{\Sigma_0,\Mt_0}$ with $\Mt_0=\tuple{V_0,D_0,\cdot_{\Mt_0}}$ is a saturated PNmatrix and $\Sigma_0\subseteq\Sigma$ then it is worth noting that $\tuple{\Sigma,\Mt_0^\Sigma}$ is also saturated. 
Indeed, it suffices to observe, first, that $\Gamma$ is a consistent theory of  $\tuple{\Sigma,\der_{\Mt_0^\Sigma}}$ if and only if $\skel(\Gamma)$ is a consistent theory of $\tuple{\Sigma_0,\der_{\Mt_0}}$, and second, that if $v^{-1}(D_0)=\skel(\Gamma)$ for some valuation $v\in\Val(\Mt_0)$ then $(v\circ\skel)\in\Val(\Mt_0^\Sigma)$ with 
$(v\circ\skel)^{-1}(D_0)=\skel^{-1}(v^{-1}(D_0))=\skel^{-1}(\skel(\Gamma))=\Gamma$.\smallskip

Of course, a PNmatrix that has a meet-closed set of bivaluations is necessarily saturated.
It turns out, however, that saturation does not necessarily imply a meet-closed set of bivaluations, but it gets sufficiently close. Let ${\tt 1}:L_\Sigma(P)\to\{0,1\}$ be the trivial bivaluation such that ${\tt 1}(A)=1$ for every $A\in L_\Sigma(P)$. It is well known that ${\vdash_{\calB}}={\vdash_{\calB'}}$ for $\calB'=\{{\tt 1}\}\cup\calB$.

\begin{lemma}\label{theone}
$\tuple{\Sigma,\Mt}$ is saturated if and only if $\BVal(\Mt)^\cap=\{{\tt 1}\}\cup\BVal(\Mt)$.
\end{lemma}
\proof{Let $\tuple{\Sigma,\Mt}$ be saturated and $X\subseteq\BVal(\Mt)$. If $X=\emptyset$ then $b_X={\tt 1}$. Otherwise, $\Gamma=b_X^{-1}(1)$ is a theory of $\tuple{\Sigma,\vdash_\Mt}$: if $\Gamma$ is inconsistent then again $b_X={\tt 1}$; if $\Gamma$ is consistent then saturation implies that $b_X\in\BVal(\Mt)$. We conclude that $\BVal(\Mt)^\cap\subseteq\{{\tt 1}\}\cup\BVal(\Mt)$.
To conclude the proof note that $\{{\tt 1}\}\cup\BVal(\Mt)\subseteq\BVal(\Mt)^\cap$, simply because ${\tt 1}=b_\emptyset$, and any bivaluation $b\in\BVal(\Mt)$ is such that $b=b_X$ with $X=\{b\}$.

Conversely, assume that $\BVal(\Mt)^\cap=\{{\tt 1}\}\cup\BVal(\Mt)$. If $\Gamma$ is a consistent theory of $\tuple{\Sigma,\vdash_\Mt}$ then, for each formula $A\notin\Gamma$ there is a bivaluation $b_A\in\BVal(\Mt)$ such that $b_A(\Gamma)\subseteq\{1\}$ and $b_A(A)=0$. Take $X=\{b_A:A\notin\Gamma\}\subseteq\BVal(\Mt)$. The meet bivaluation $b_X\in\BVal(\Mt)^\cap$ and $b_X^{-1}(1)=\Gamma$, which implies that $b_X\neq\{{\tt 1}\}$ as $\Gamma$ is consistent. We conclude that $b_X\in\BVal(\Mt)$ and thus that $\Mt$ is saturated.\qed}\\

Under the assumption that the PNmatrices are saturated, Theorem~\ref{multiplefibring} can now be adapted to the single-conclusion case too.

\begin{theorem}\label{saturatedfibring}
The combination of single-conclusion logics characterized by saturated PNmatrices is the single-conclusion logic characterized by their strict product, that is, if both 
$\tuple{\Sigma_1,\Mt_1}$ and $\tuple{\Sigma_2,\Mt_2}$ are saturated then we have
$\tuple{\Sigma_1,\vdash_{\Mt_1}}\bullet\tuple{\Sigma_2,\vdash_{\Mt_2}}=\tuple{\Sigma_1\cup\Sigma_2,\vdash_{\Mt_1\ast \Mt_2}}$.
\end{theorem}
\proof{The result follows directly from Proposition~\ref{bivsingle}, and Lemmas~\ref{strictvals}~and~\ref{theone}. 

With $\calB_1=\BVal(\Mt_1)$ and $\calB_2=\BVal(\Mt_2)$,   note that we necessarily have
$\{{\tt 1}\}\cup\BVal(\Mt_1\ast\Mt_2)=\{{\tt 1}\}\cup(\calB_1^{\Sigma_1\cup\Sigma_2}\cap\calB_2^{\Sigma_1\cup\Sigma_2})
=(\{{\tt 1}\}\cup\calB_1^{\Sigma_1\cup\Sigma_2})\cap(\{{\tt 1}\}\cup\calB_2^{\Sigma_1\cup\Sigma_2}\})={(\calB_1^{\Sigma_1\cup\Sigma_2})^{\cap}\cap (\calB_2^{\Sigma_1\cup\Sigma_2})^{\cap}}=\calB_{12}^{\sing}$.\qed}\\

As a corollary of this proof it also results that $\tuple{\Sigma_1\cup\Sigma_2,{\Mt_1\ast \Mt_2}}$ is a saturated PNmatrix, and thus that saturation is preserved by the strict product operation. 

\begin{example}\label{negand2}(Strict product and saturation).\\
Recall the discussion in Example~\ref{negand1} about the single-conclusion combination of the $\neg$-only and $\wedge$-only fragments of classical logic, $\tuple{\Sigma^\neg_\mathsf{cls},\vdash^\neg_\mathsf{cls}}$ and $\tuple{\Sigma^\wedge_\mathsf{cls},\vdash^\wedge_\mathsf{cls}}$. 

Since the set of all classical bivaluations for the $\wedge$-only language is meet-closed, the classical $\Sigma^\wedge_\mathsf{cls}$-matrix $\Mt_1=\tuple{\{0,1\},\{1\},\cdot_{\Mt_1}}$ with $\wedge_{\Mt_1}=\wedge_{\TWO}$, whose truth-table is shown below, is saturated.

On the other hand, as seen before, the set of all classical bivaluations for the $\neg$-only language is not meet-closed, and the corresponding two-valued matrix is not saturated. Instead, let us consider
the three-valued $\Sigma^\neg_\mathsf{cls}$-matrix $\Mt_2=\tuple{\{0,\frac{1}{2},1\},\{1\},\cdot_{\Mt_2}}$ as defined below.

 \begin{center}
    \begin{tabular}{c c c c}
    %\hline
    $\land_{\Mt_1}$&\vline&  $0$ & $1$   \\ 
    \hline
    $0$ & \vline&$0$ &  $0$ \\ %\hline
     $1$  &\vline& $0$  &  $1$\\  %\hline
   \;
   % $-1$ & $0$   \\
 %   \hline
  \end{tabular}
\qquad
 \begin{tabular}{c | c c}
    %\hline
    &$\neg_{\Mt_2}$     \\ 
    \hline
    $0$ & $1$   \\ %\hline
    $\frac{1}{2}$ & $\frac{1}{2}$\\
     $1$  &$0$    %\hline
   % $-1$ & $0$   \\
 %   \hline
  \end{tabular}
  \end{center} 

It turns out that $\vdash_{\Mt_2}{=}\vdash^\neg_\mathsf{cls}$. Further, $\Mt_2$ is saturated. To see this, given a consistent theory 
$\Gamma$ of $\tuple{\Sigma^\neg_\mathsf{cls},\vdash^\neg_\mathsf{cls}}$, it suffices to consider the valuation $v\in\Val(\Mt_2)$ such that 

$$v(A)=
\begin{cases}
1 & \mbox{ if }\;A\in\Gamma \\
0 &  \mbox{ if }\;\neg A\in\Gamma\\
\frac{1}{2} & \mbox{ if }\;A,\neg A\notin\Gamma.
\end{cases}
$$

By the consistency of $\Gamma$ we know that $v$ is well defined, because $A$ and $\neg A$ cannot both be in $\Gamma$ due to rule  
$\frac{\;\;p\quad  \neg p\;\;}{q}$. Additionally, rules  $\frac{p}{\;\;\neg\neg p\;\;}$ and $\frac{\;\;\neg \neg p\;\;}{p}$ guarantee that $v$ is indeed a valuation of $\Mt_2$.\smallskip

Given the saturation of the two matrices, according to Theorem~\ref{saturatedfibring}, the resulting combined logic $\tuple{\Sigma^\neg_\mathsf{cls},\vdash^\neg_\mathsf{cls}}\bullet\tuple{\Sigma^\wedge_\mathsf{cls},\vdash^\wedge_\mathsf{cls}}{=}\tuple{\Sigma^{\neg,\wedge}_\mathsf{cls},\vdash}$ is characterized by the strict product $\Mt_1\ast\Mt_2$. By Definition~\ref{strictproduct} the resulting truth-values are $00,0\frac{1}{2},11$ (where, again, we are simplifying each pair $(x,y)$ to simply $xy$). Renaming these values to simply $0,\frac{1}{2},1$, respectively, we have that $\Mt_1\ast\Mt_2=\tuple{\{0,\frac{1}{2},1\},\{1\},\cdot_\ast}$ is defined by the truth-tables below.

 \begin{center}
    \begin{tabular}{c |c c c}
    %\hline
    $\land_{\ast}$ &  $0$ & $\frac{1}{2}$& $1$   \\ 
    \hline
    $0$ & $0,\frac{1}{2}$ &  $0,\frac{1}{2}$ & $0,\frac{1}{2}$ \\ %\hline
    $\frac{1}{2}$ & $0,\frac{1}{2}$ &  $0,\frac{1}{2}$ & $0,\frac{1}{2}$ \\ %\hline
     $1$  & $0,\frac{1}{2}$  & $0,\frac{1}{2}$ &  $1$   % $-1$ & $0$   \\
 %   \hline
  \end{tabular}
\qquad
 \begin{tabular}{c | c c}
    %\hline
    &$\neg_{\ast}$     \\ 
    \hline
    $0$ & $1$   \\ %\hline
    $\frac{1}{2}$ & $\frac{1}{2}$\\
     $1$  &$0$    %\hline
   % $-1$ & $0$   \\
 %   \hline
  \end{tabular}
  \end{center} 

A valuation $v\in\Val(\Mt_1\ast\Mt_2)$ such that $v(p)=0$, $v(\neg p)=1$ and $v(p\wedge p)=v(\neg(p\wedge p))=\frac{1}{2}$ shows that 
$\neg p\not\vdash\neg(p\wedge p)$. 

Of course, this implies that putting together calculi for each of the connectives is not enough to fully characterize the way they interact in classical logic, as can be visually confirmed from the rules in Example~\ref{classical1}.
\hfill $\triangle$
\end{example}

As this point, we should question the scope of applicability of Theorem~\ref{saturatedfibring}, and how often we can expect to find naturally saturated PNmatrices as in the example above. What should we do when given PNmatrices that are not saturated? Fortunately, there is a simple way of transforming a PNmatrix into a saturated PNmatrix both characterizing the same single-conclusion logic~\footnote{Note that the operation works precisely by weakening the associated multiple-conclusion logic, a subject to which we will return later on.}.

\begin{definition}
Let $\tuple{\Sigma,\Mt}$ with $\Mt=\tuple{V,D,\cdot_{\Mt}}$ be a PNmatrix. The \emph{$\omega$-power of $\tuple{\Sigma,\Mt}$} is the PNmatrix $\tuple{\Sigma,\Mt^\omega}$ with $\Mt^\omega=\tuple{V^\omega,D^\omega,\cdot_{\Mt^\omega}}$ such that
\begin{itemize}
\item $V^\omega=\{\tuple{x_i}_{i\in\nats}: x_i\in V\textrm{ for every }i\in \nats\}$, 
\item $D^\omega=\{\tuple{x_i}_{i\in\nats}: x_i\in D\textrm{ for every }i\in \nats\}$, and 
\item for every $n\in\nats_0$, $c\in(\Sigma_1\cup\Sigma_2)^{(n)}$ and $\tuple{x_{1,i}}_{i\in\nats},\dots,\tuple{x_{n,i}}_{i\in\nats}\in V^\omega$, we let 
$\tuple{y_i}_{i\in\nats}\in c_{\Mt^\omega}(\tuple{x_{1,i}}_{i\in\nats},\dots,\tuple{x_{n,i}}_{i\in\nats})$ if and only if the following condition holds:
$$y_i\in c_\Mt(x_{1,i},\dots,x_{n,i})\textrm{ for every }i\in \nats.$$
\end{itemize}
\end{definition}

For each $k\in\nats$ we use $\pi_k:V^\omega\to V$ to denote the obvious projection function, i.e., $\pi_k(\tuple{x_i}_{i\in\nats})=x_k$. Clearly, we have that $v\in\Val(\Mt^\omega)$ if and only if $(\pi_k\circ v)\in\Val(\Mt)$ for every $k\in\nats$.

\begin{lemma}\label{wpower}
If $\tuple{\Sigma,\Mt}$ is a PNmatrix then $\tuple{\Sigma,\Mt^\omega}$ is saturated, and further we have ${\vdash_\Mt}={\vdash_{\Mt^\omega}}$.
\end{lemma}
\proof{We first show that ${\vdash_\Mt}={\vdash_{\Mt^\omega}}$. 

Suppose that $\Gamma\vdash_\Mt A$, and let $v\in\Val(\Mt^\omega)$ be such that $v(\Gamma)\subseteq D^\omega$. 
For every $k\in\nats$, this means that $(\pi_k\circ v)(\Gamma)\subseteq D$ and thus that $(\pi_k\circ v)(A)\in D$.
We conclude that $v(A)\in D^\omega$ and so $\Gamma\vdash_{\Mt^\omega} A$.

Suppose now that $\Gamma\vdash_{\Mt^\omega} A$, and let $v\in\Val(\Mt)$ be such that $v(\Gamma)\subseteq D$.
Easily, $v^\omega$ defined by $\pi_k(v^\omega(A))=v(A)$ for every $k\in\nats$ is such that $v^\omega\in \Val(\Mt^\omega)$, and $v^\omega(\Gamma)\subseteq D^\omega$.
Therefore, we have that $v^\omega(A)\in D^\omega$ and therefore, for any $k\in\nats$, $\pi_k(v^\omega(A))=v(A)\in D$.
We conclude that $\Gamma\vdash_{\Mt} A$.\smallskip

To show that $\tuple{\Sigma,\Mt^\omega}$ is saturated, let $\Gamma$ be a consistent theory of $\tuple{\Sigma,\vdash_\Mt}$, and fix $A\notin \Gamma$.
For each formula $B\notin\Gamma$ we know that there exists a valuation $v_B\in\Val(\Mt)$ such that $v_B(\Gamma)\subseteq D$ and $v_B(B)\notin D$.
Fix an enumeration $\eta:\nats\to L_\Sigma(P)$ and define, for each $k\in\nats$, the valuation $v_k\in\Val(\Mt)$ such that

$$v_k=\left\{\begin{array}{ll}
v_{\eta(k)} & \textrm{if }\eta(k)\notin\Gamma\\
v_A & \textrm{if }\eta(k)\in\Gamma
\end{array}\right..$$

Let $v\in\Val(\Mt^\omega)$ be such that $(\pi_k\circ v)=v_k$ for every $k\in\nats$.
We have that $v(\Gamma)\subseteq D^\omega$ because $v_k(\Gamma)\subseteq D$ for every $k\in\nats$, 
and for every $B\notin\Gamma$, we have that $v(B)\notin D^\omega$ because with $k=\eta^{-1}(B)$ we have $v_k(B)=v_B(B)\notin D$.
We conclude that $\tuple{\Sigma,\Mt^\omega}$ is saturated\footnote{The proof of Lemma~\ref{wpower} actually shows that the bivaluations of $\Mt^\omega$ are closed under non-empty countable meets. This is enough, in our case, since the empty meet corresponds to the irrelevant ${\tt 1}$ bivaluation, and also because non-countable meets are not necessary given that we always work with denumerable languages.}.\qed}\smallskip

The following result follows from Theorem~\ref{saturatedfibring}~and~Lemma~\ref{wpower}.

\begin{theorem}\label{singlefibring}
The combination of single-conclusion logics characterized by PNmatrices is the multiple-conclusion logic characterized by the strict product of their $\omega$-powers, that is, 
$\tuple{\Sigma_1,\vdash_{\Mt_1}}\bullet\tuple{\Sigma_2,\vdash_{\Mt_2}}=\tuple{\Sigma_1\cup\Sigma_2,\vdash_{\Mt_1^\omega\ast \Mt_2^\omega}}$.
\end{theorem}

Though still usable, this result does not provide us with a finite-valued semantics for the combined logic, even if departing from finite-valued PNmatrices. Indeed, in the relevant cases, a $\omega$-power PNmatrix is always non-countable, the exceptions being $\omega$-powers of PNmatrices with at most one truth-value, which are actually saturated to begin with. Note also that the situation is not unexpected, as we know that some combinations cannot be endowed with finite-valued semantics~\cite{charfinval}, and it is playing the role of meet-closure in the case of bivaluations. Still, there are cases when saturation is not necessary (we will see a relevant example of this phenomenon in Subsection~\ref{sub:axioms} below), or when saturation can be achieved by a finite power of the given PNmatrix. We present some examples below, and further discuss this question in the conclusion of the paper.

\begin{example}\label{lk2}(Three-valued implications of Kleene and  \L ukasiewicz).\\
Recall Example~\ref{lk}, where we discussed the combination of the three-valued implications of Kleene and  \L ukasiewicz in the multiple-conclusion setting, with shared implication, and we concluded that the resulting logic coincided with the multiple-conclusion version of the implication-fragment of classical logic.\smallskip

Now, we shall see that, incidentally, the corresponding single-conclusion combination also coincides with classical implication, albeit for distinct reasons. 
For easier readability, we recall here the signature $\Sigma$ with $\Sigma^{(2)}=\{\to\}$ and $\Sigma^{{n}}=\emptyset$ for $n\neq 2$, and the three-valued implication matrices of Kleene and  \L ukasiewicz, respectively $\mathbb{K}=\tuple{\{0,\frac{1}{2},1\},\{1\},\cdot_\mathbb{K}}$ and 
$\mathbb{L}=\tuple{\{0,\frac{1}{2},1\},\{1\},\cdot_\mathbb{L}}$, defined below.

 \begin{center}
 \begin{tabular}{c | c c c}
$\to_{\mathbb{K}}$& $0$ & $\frac{1}{2}$ & $1$\\
\hline
$0$& $1$ & $1$ & $1$\\
$\frac{1}{2}$ & $\frac{1}{2}$ & $\frac{1}{2}$ & $1$ \\
$1$ & $0$ & $\frac{1}{2}$ & $1$
\end{tabular}
\qquad
 \begin{tabular}{c | c c c}
$\to_{\mathbb{L}}$& $0$ & $\frac{1}{2}$ & $1$\\
\hline
$0$& $1$ & $1$ & $1$\\
$\frac{1}{2}$ & $\frac{1}{2}$ & $1$ & $1$ \\
$1$ & $0$ & $\frac{1}{2}$ & $1$
\end{tabular}  
\end{center}

Concerning the combined logic $\tuple{\Sigma,\vdash}=\tuple{\Sigma,\vdash_\mathbb{K}}\bullet\tuple{\Sigma,\vdash_\mathbb{L}}$, 
Theorem~\ref{saturatedfibring} cannot be applied directly as the matrices are not saturated, which renders the strict product obtained in Example~\ref{lk} useless as our envisaged semantics for the combined logic.\smallskip 

Concerning $\mathbb{K}$, let $\Gamma=\{p\to p\}^{\vdash_\mathbb{K}}$. The theory is consistent, as in particular we have 
$p\to p \not\vdash_{\mathbb{K}} p$ (as can be confirmed by any valuation $v\in\Val(\mathbb{K})$ with $v(p)=0$) and also 
$p\to p \not\vdash_\mathbb{K} p\to q$ (as witnessed by any valuation $v\in\Val(\mathbb{K})$ with $v(p)=1$ and $v(q)=0$). Therefore, 
$p,p\to q\notin\Gamma$. However, there is no valuation $v\in\Val(\mathbb{K})$ such that $\Gamma=v^{-1}(\{1\})$, simply because if $v(p\to p)=1$ and $v(p)\neq 1$ then by just inspecting the truth-table we conclude that $v(p)=0$ and $v(p\to q)=1$.

Concerning $\mathbb{L}$, consider $\Gamma=\emptyset^{\vdash_\mathbb{L}}$. The theory is consistent, as in particular we have 
$\emptyset \not\vdash_{\mathbb{L}} p$ (as can be confirmed by any valuation $v\in\Val(\mathbb{L})$ with $v(p)=0$),  
$\emptyset \not\vdash_\mathbb{L} p\to q$ (as witnessed by any valuation $v\in\Val(\mathbb{L})$ with $v(p)=1$ and $v(q)=0$), and also 
$\emptyset \not\vdash_\mathbb{L} q\to r$ (as witnessed by any valuation $v\in\Val(\mathbb{L})$ with $v(q)=1$ and $v(r)=0$). Therefore, 
$p,p\to q,q\to r \notin\Gamma$. However, there is no valuation $v\in\Val(\mathbb{L})$ such that $\Gamma=v^{-1}(\{1\})$, simply because if $v(p),v(p\to q)\neq 1$ then by just inspecting the truth-table we conclude that $v(p)=\frac{1}{2}$, $v(q)=0$ and $v(q\to r)=1$.\smallskip

We must therefore consider the $\omega$-power of each of the matrices. Note that in both cases, the resulting truth-values are denumerable tuples 
$\tuple{x_i}_{i\in\nats}$ with each $x_i\in\{0,\frac{1}{2},1\}$. For simplicity, 
we will rename each $\tuple{x_i}_{i\in\nats}$ as the pair of sets $\tuple{X_0,X_1}$ with $X_0=\{i\in\nats:x_i=0\}$ and $X_1=\{i\in\nats:x_i=1\}$. 
Of course, then, $X_{\frac{1}{2}}=\{i\in\nats:x_i=\frac{1}{2}\}=\nats\setminus(X_0\cup X_1)$. The resulting set of truth-values corresponds to 
$V^\omega=\{\tuple{X_0,X_1}:X_0,X_1\subseteq\nats, X_0\cap X_1=\emptyset\}$. As both matrices have $1$ as their only designated value, we also get in both cases that $D^\omega=\{\tuple{\emptyset,\nats}\}$. The resulting matrices are thus 
$\mathbb{K}^\omega=\tuple{V^\omega,D^\omega,\cdot_{\mathbb{K}^\omega}}$ and 
$\mathbb{L}^\omega=\tuple{V^\omega,D^\omega,\cdot_{\mathbb{L}^\omega}}$, defined according to the tables below. 

 \begin{center}
 \begin{tabular}{c | c c c}
$\to_{\mathbb{K}^\omega}$& $\tuple{Y_0,Y_1}$\\
\hline
$\tuple{X_0,X_1}$& $\tuple{X_1\cap Y_0,X_0\cup Y_1}$
\end{tabular}
\quad
 \begin{tabular}{c | c c c}
$\to_{\mathbb{L}^\omega}$& $\tuple{Y_0,Y_1}$\\
\hline
$\tuple{X_0,X_1}$& $\tuple{X_1\cap Y_0,X_0\cup Y_1\cup (X_{\frac{1}{2}}\cap Y_{\frac{1}{2}})}$
\end{tabular}
\end{center}

According to Theorem~\ref{singlefibring}, we know that $\tuple{\Sigma,\vdash_\mathbb{K}}\bullet\tuple{\Sigma,\vdash_\mathbb{L}}$ is characterized by the infinite $\Sigma$-Pmatrix $\mathbb{K}^\omega\ast\mathbb{L}^\omega$. Representing pairs of pairs as four-tuples we get the set of truth-values 
$V_\ast=\{\tuple{X_0,X_1,Y_0,Y_1}:X_0,X_1,Y_0,Y_1\subseteq\nats, X_0\cap X_1=Y_0\cap Y_1=\emptyset\textrm{, and }X_1=\nats\textrm{ iff }Y_1=\nats\}$, and designated values 
$D_\ast=\{\tuple{\emptyset,\nats,\emptyset,\nats}\}$. The resulting strict product is then 
$\mathbb{K}^\omega\ast\mathbb{L}^\omega=\tuple{V_\ast,D_\ast,\cdot_\ast}$ as defined by the table below.

 \begin{center}
 \begin{tabular}{c  c c ccc}
$\to_{\ast}$&\vline& $\tuple{Z_0,Z_1,W_0,W_1}$\\
\hline
$\tuple{X_0,X_1,Y_0,Y_1}$&\vline& $\tuple{X_1\cap Z_0,X_0\cup Z_1,Y_1\cap W_0,Y_0\cup W_1\cup(Y_{\frac{1}{2}}\cap W_{\frac{1}{2}})}$\\
&\vline& $\textrm{ provided }X_0\cup Z_1=\nats\textrm{ iff }Y_0\cup W_1\cup(Y_{\frac{1}{2}}\cap W_{\frac{1}{2}})=\nats$\\ \hline
$\tuple{X_0,X_1,Y_0,Y_1}$&\vline& $\emptyset$\\
&\vline& $\textrm{ otherwise }$
%&\vline& $\textrm{ otherwise }
\end{tabular}
 \end{center}

This semantics has a non-trivial look, but it turns out that one can still draw valuable conclusions from it. Due to the partiality of $\mathbb{K}^\omega\ast\mathbb{L}^\omega$ it is worth noting that the Pmatrix has a lot of spurious values. 
For a truth-value $\tuple{X_0,X_1,Y_0,Y_1}$, if $\tuple{X_0,X_1,Y_0,Y_1}\to_\ast\tuple{X_0,X_1,Y_0,Y_1}$ is non-empty then  it must contain a single element $\tuple{
X_1\cap X_0,X_0\cup X_1,Y_1\cap Y_0,Y_0\cup Y_1\cup Y_{\frac{1}{2}}}=\tuple{
\emptyset,X_0\cup X_1,\emptyset,\nats}=\tuple{\emptyset,\nats,\emptyset,\nats}$. Of course, this happens only when $X_{\frac{1}{2}}=\emptyset$, or otherwise $\tuple{X_0,X_1,Y_0,Y_1}$ is spurious.

This observation implies that if $v\in\Val(\mathbb{K}^\omega\ast\mathbb{L}^\omega)$ then $\pi_1\circ v\in \Val(\mathbb{K}^\omega)$ never uses the $\frac{1}{2}$ value of matrix $\mathbb{K}$. As $\mathbb{K}$ behaves classically for the $0,1$ values, we conclude that all bivaluations in $\BVal(\mathbb{K}^\omega\ast\mathbb{L}^\omega)$ are classical. On the other hand, we know from Lemmas~\ref{strictvals}, \ref{theone} and~\ref{wpower} that $\BVal(\mathbb{K}^\omega\ast\mathbb{L}^\omega)=\BVal(\mathbb{K}^\omega)\cap\BVal(\mathbb{L}^\omega)=
\BVal(\mathbb{K})^\cap\cap\BVal(\mathbb{L})^\cap$, because obviously ${\tt 1}\in\BVal(\mathbb{K}),\BVal(\mathbb{L})$. Since not just $\mathbb{K}$ but also $\mathbb{L}$ behaves classically for the $0,1$ values, we conclude that all classical bivaluations are in $\BVal(\mathbb{K}^\omega\ast\mathbb{L}^\omega)$. Thus, we have that $\tuple{\Sigma,\vdash_\mathbb{K}}\bullet\tuple{\Sigma,\vdash_\mathbb{L}}=\tuple{\Sigma,\vdash_{\mathbb{K}^\omega\ast\mathbb{L}^\omega}}=\tuple{\Sigma^\to_\mathsf{cls},\vdash^\to_\mathsf{cls}}$ is precisely the implication fragment of classical logic.\smallskip

Again, we can confirm this fact by putting together calculi for the logics, according to Proposition~\ref{fibringcalculi}. Even without listing the rules, and since it is well known that $\vdash_{\mathbb{K}},\vdash_{\mathbb{L}}{\subseteq} \vdash^\to_\mathsf{cls}$, it suffices to check that all rules of the calculus for $\vdash^\to_\mathsf{cls}$ (see Example~\ref{classical1}) obtain when we join them. The rule of \emph{modus ponens} is unproblematic, as we have both $p,p\to q\der_{\mathbb{K}}q$ and $p,p\to q\der_{\mathbb{L}}q$. Concerning the classical axioms $p\to(q\to p),(p\to(q\to r))\to ((p\to q)\to (p\to r)),(((p\to q)\to p)\to p)$, or actually any other classical tautology $A$, note that $\{p\to p:p\in\var(A)\}\vdash_{\mathbb{K}}A$, simply because a valuation $v\in\Val(\mathbb{K})$ such that $v(p\to p)=1$ must be classical, as $v(p)\neq\frac{1}{2}$. 
On the side of  \L ukasiewicz, it is clear that $\emptyset\vdash_{\mathbb{K}}B\to B$ for any formula $B$. Hence, in the combined logic, we know that $\emptyset\vdash p\to p$ for each $p\in\var(A)$, and also that $\{p\to p:p\in\var(A)\}\vdash A$ which implies that $\emptyset\vdash A$ for any classical tautology $A$.
\hfill $\triangle$
\end{example}

One must, of course, look again at the combination of fragments of classical logic, now in the single-conclusion setting.

\begin{example}\label{classical2}(Combining fragments of classical logic).\\
Recall Example~\ref{classical1}, as well as the multiple-conclusion scenario that we have explored in Example~\ref{multi2}. We will now see how much challenging and interesting it is to combine fragments of classical logic in the single-conclusion setting, by means of four distinct cases.\smallskip

{\bf(1)}~Let us first consider the combination of the fragments of classical logic corresponding to $\Sigma^\wedge_\mathsf{cls}$ and $\Sigma^\top_\mathsf{cls}$. It is easy to see that the corresponding two-valued classical matrices are both saturated. Hence, the combined logic $\tuple{\Sigma^\wedge_\mathsf{cls},\vdash^\wedge_\mathsf{cls}}\bullet\tuple{\Sigma^\top_\mathsf{cls},\vdash^\top_\mathsf{cls}}$ is directly characterized by the corresponding strict product, which coincides with the
 two-valued classical matrix for the fragment $\Sigma^{\wedge,\top}_\mathsf{cls}$. We conclude thus, that $\tuple{\Sigma^\wedge_\mathsf{cls},\vdash^\wedge_\mathsf{cls}}\bullet\tuple{\Sigma^\top_\mathsf{cls},\vdash^\top_\mathsf{cls}}=\tuple{\Sigma^{\wedge,\top}_\mathsf{cls},\vdash^{\wedge,\top}_\mathsf{cls}}$, which is compatible with our knowledge that joining single-conclusion calculi for $\vdash^\wedge_\mathsf{cls}$ and $\vdash^\top_\mathsf{cls}$ yields a calculus for $\vdash^{\wedge,\top}_\mathsf{cls}$.\smallskip
 
 {\bf(2)}~Now, let us reanalyze the combination of $\tuple{\Sigma^\neg_\mathsf{cls},\vdash^\neg_\mathsf{cls}}$ and $\tuple{\Sigma^\wedge_\mathsf{cls},\vdash^\wedge_\mathsf{cls}}$. We already know from Examples~\ref{negand1} and~\ref{negand2} that the combination $\tuple{\Sigma^\neg_\mathsf{cls},\vdash^\neg_\mathsf{cls}}\bullet\tuple{\Sigma^\wedge_\mathsf{cls},\vdash^\wedge_\mathsf{cls}}$ is weaker than $\vdash^{\neg,\wedge}_\mathsf{cls}$. However, since the two-valued classical matrix for negation is not saturated, we have instead considered an adequate three-valued saturated matrix. Of course, we did not need to know that such a matrix existed, and instead we could have blindly obtained the $\omega$-power of the two-valued classical matrix for negation, and then obtained its strict product with the saturated two-valued classical matrix for negation. The semantics thus obtained would be less amiable, but still characterizes the combination. Alternatively, we could have noted that the four-valued $2$-power of the matrix of negation would already be saturated (the three-valued matrix used before is a simplification of it).\smallskip
 
{\bf(3)}~As a last interesting example of disjoint combination let us consider the fragments $\tuple{\Sigma^\wedge_\mathsf{cls},\vdash^\wedge_\mathsf{cls}}$ and $\tuple{\Sigma^\vee_\mathsf{cls},\vdash^\vee_\mathsf{cls}}$. We already know from~\cite{softcomp} that the combination $\tuple{\Sigma^{\wedge,\vee}_\mathsf{cls},\vdash}=\tuple{\Sigma^\wedge_\mathsf{cls},\vdash^\wedge_\mathsf{cls}}\bullet\tuple{\Sigma^\vee_\mathsf{cls},\vdash^\vee_\mathsf{cls}}$ is weaker than $\vdash^{\wedge,\vee}_\mathsf{cls}$. Indeed, we have that $p\vee(p\wedge p)\not\vdash p$, as can be confirmed by any valuation with 
$\emptyset\subsetneq v(p)\subsetneq\nats$ and $v(p\wedge p)=\nats\setminus v(p)\neq \nats$, which then implies that $v(p\vee(p\wedge p))=\nats$, on the Nmatrix resulting from the strict product of the two-valued classical matrix for conjunction and the $\omega$-power of the two-valued classical matrix for disjunction, corresponding (after renaming) to $\Mt_\ast=\tuple{\wp(\nats),\{\nats\},\cdot_\ast}$ as defined by the tables below.
 
 \begin{center}
 \begin{tabular}{c  c c ccc}
$\wedge_{\ast}$&\vline& $Y$\\
\hline
$X$&\vline& $\nats$\\
&\vline& $\textrm{ provided }X=Y=\nats$\\ \hline
$X$&\vline& $\{Z\subseteq\nats:Z\neq\nats\}$\\
&\vline& $\textrm{ otherwise }$
\end{tabular}
\quad
  \begin{tabular}{c  c c c}
$\vee_\ast$&\vline& $Y$\\ \hline
$X$&\vline& $X\cup Y$\\
\\ \\
\\
\end{tabular}

 \end{center}

Of course a simpler semantics would be desirable, but we know that the two-valued classical matrix for disjunction is not saturated. For instance, note that $\{p\vee q\}\not\vdash^\vee_\mathsf{cls}p$ and $\{p\vee q\}\not\vdash^\vee_\mathsf{cls}q$ but any valuation $v$ of the two-valued matrix that sets $v(p\vee q)=1$ necessarily must have $v(p)=1$ or $v(q)=1$.\smallskip

{\bf(4)}~Finally, we shall look at a non-disjoint example, the combination of the fragments corresponding to $\Sigma^{\neg,\to}_\mathsf{cls}$ and $\Sigma^{\vee,\to}_\mathsf{cls}$. It seems clear, from the calculus for $\vdash_\mathsf{cls}$ shown in Example~\ref{classical1}, that all rules involving negation or disjunction only additionally use the shared connective of implication. We may therefore conjecture that 
$\tuple{\Sigma^{\neg,\to}_\mathsf{cls},\vdash^{\neg,\to}_\mathsf{cls}}\bullet\tuple{\Sigma^{\vee,\to}_\mathsf{cls},\vdash^{\vee,\to}_\mathsf{cls}}=\tuple{\Sigma^{\neg,\vee,\to}_\mathsf{cls},\vdash^{\neg,\vee,\to}_\mathsf{cls}}$. To confirm this, we first note that none of the two-valued classical matrices for the fragments is saturated. This is simply because both include implication, $\emptyset\not\vdash_\mathsf{cls} p$ and 
$\emptyset\not\vdash_\mathsf{cls} p\to q$ but any classical valuation will have $v(p)=1$ or $v(p\to q)=1$. Thus, we need to consider the strict product of the $\omega$-power of each of the matrices, which results (after renaming) in the PNmatrix
$\Mt_\ast=\tuple{\wp(\nats)\times\wp(\nats),\{(\nats,\nats)\},\cdot_\ast}$ as defined by the tables below, where we use $\overline{X}=\nats\setminus X$.

 \begin{center}
 \begin{tabular}{c  c c ccc}
&\vline&$\neg_{\ast}$\\
\hline
$\tuple{X,Y}$&\vline& $\tuple{\nats,\nats}$\\
&\vline& $\textrm{ provided }X=\emptyset$\\ \hline
$\tuple{X,Y}$&\vline& $\{\tuple{\overline{X},U}:U\neq\nats\}$\\
&\vline& $\textrm{ otherwise }$
\end{tabular}
\qquad
 \begin{tabular}{c  c c ccc}
$\vee_{\ast}$&\vline& $\tuple{Z,W}$\\
\hline
$\tuple{X,Y}$&\vline& $\tuple{\nats,\nats}$\\
&\vline& $\textrm{ provided }Y\cup W=\nats$\\ \hline
$\tuple{X,Y}$&\vline& $\{\tuple{U,Y\cup W}:U\neq\nats\}$\\
&\vline& $\textrm{ otherwise }$
\end{tabular}
\qquad
 \begin{tabular}{c  c c ccc}
$\to_{\ast}$&\vline& $\tuple{Z,W}$\\
\hline
$\tuple{X,Y}$&\vline& $\tuple{\overline{X}\cup Z,\overline{Y}\cup W}$\\
&\vline& $\textrm{ provided }\overline{X}\cup Z=\nats\textrm{ iff }\overline{Y}\cup W=\nats$\\ \hline
$\tuple{X,Y}$&\vline& $\emptyset$\\
&\vline& $\textrm{ otherwise }$
\end{tabular}
 \end{center}

All classical bivaluations can be simply obtained in $\Mt_*$ by taking only the two values $\tuple{\emptyset,\emptyset}$ and $\tuple{\nats,\nats}$.
Seeing that all valuations of $\Mt_\ast$ are classical (i.e., they respect the operations in $\TWO^\omega$) is slightly harder. Let $v\in\Val(\Mt_\ast)$. 

We have $v(p\to p)=\tuple{\nats,\nats}$ and hence $v(\neg(p\to p))=\tuple{\emptyset,U_0}$ for some $U_0\subsetneq\nats$. Consequently, for any formula $A$, if $v(A)=\tuple{Z,W}$ then $v(\neg(p\to p)\to A)\in(\tuple{\emptyset,U_0}\to_\ast\tuple{Z,W})$ and thus $v(\neg(p\to p)\to A)=\tuple{\overline{\emptyset}\cup Z,\overline{U_0}\cup W}=\tuple{\nats\cup Z,\overline{U_0}\cup W}=\tuple{\nats,\nats}$, which implies $\overline{U_0}\cup W=\nats$, or equivalently $U_0\subseteq W$. 
Further, $v(\neg A)=\tuple{\overline{Z},T}$ for some $T$ such that $U_0\subseteq T\subseteq\nats$. We can see that $v(A\to(\neg A\to \neg(p\to p)))=
\tuple{\overline{Z}\cup Z\cup \emptyset,\overline{W}\cup\overline{T}\cup U_0}=\tuple{\nats,\nats}$, which means that $\overline{W}\cup\overline{T}\cup U_0=\nats$ and thus that $W\cap T=U_0$. Also, we can see that $v((\neg A\to A)\to A)=
\tuple{\overline{Z}\cup Z,(T\cap\overline{W})\cup W}=\tuple{\nats,\nats}$, which means that $(T\cap\overline{W})\cup W=\nats$ and thus that $\overline{W}\subseteq T$. From these observations, we conclude that $T=\overline{W}\cup U_0$.

Let $f:\nats\to\overline{U_0}$ be any surjective function. We claim that the function defined by $v'(B)=f^{-1}(\pi_2(v(B))\setminus U_0)$ is a valuation over $\TWO^\omega$ (on the relevant connectives), which is compatible with $v$, as whenever $U_0\subseteq X$ we have $f^{-1}(X\setminus U_0)=\nats$ if and only if $X=\nats$. It is straightforward to check that $v'(B\vee C)=v'(B)\cup v'(C)$ and $v'(B\to C)=\overline{v'(B)}\cup v'(C)$, and also easy to see that $v'(\neg B)=\overline{v'(B)}$ using the observations in the previous paragraph.
\hfill $\triangle$
\end{example}

\section{Universal properties}\label{secuniv}

According to the very successful mathematical approach to \emph{General Systems Theory} initiated by J. Goguen in~\cite{goguen1,goguen2}, composition operations should always be explained by universal properties, in the sense of category theory~\cite{maclane,adamek}. This approach has been used also with respect to combinations of logics, for instance in~\cite{ccal:car:jfr:css:04,acs:css:ccal:98a}. Herein, we briefly show how our results are explained by simple universal constructions.\smallskip

\subsection{PNmatrices and rexpansions}

In order to explore the relationships among PNmatrices, let us consider a suitable notion of homomorphism. If $\Sigma_0\subseteq\Sigma$, $\tuple{\Sigma_0,\Mt_0}$ and $\tuple{\Sigma,\Mt}$ are PNmatrices, with $\Mt_0=\tuple{V_0,D_0,\cdot_{\Mt_0}}$ and $\Mt=\tuple{V,D,\cdot_{\Mt}}$, a \emph{strict homomorphism} $h:\tuple{\Sigma,\Mt}\to\tuple{\Sigma_0,\Mt_0}$ is a function $h:V\to V_0$ such that $h^{-1}(D_0)=D$, and for every $n\in\nats_0$, $c\in\Sigma_0^{(n)}$, and $x_1,\dots x_n\in V$, $h(c_\Mt(x_1,\dots,x_n))\subseteq c_{\Mt_0}(h(x_1),\dots,h(x_n))$.
If $v\in\Val(\Mt)$ then $(h\circ v)\in\Val({\Mt_0^\Sigma})$. Further, the strictness condition $h^{-1}(D_0)=D$  implies that $v(A)$ and $(h\circ v)(A)$ are compatible for every formula $A\in L_\Sigma(P)$. Consequently, we have that
$\BVal(\Mt)\subseteq\BVal(\Mt_0^\Sigma)$, and therefore ${\propto_{\Mt_0}}\subseteq{\propto_{\Mt}}$ both in the single and the multiple-conclusion cases.  PNmatrices with strict homomorphisms constitute a category~$\cat{PNMatr}$.\smallskip

Expectedly, the strict product of PNmatrices as introduced in Definition~\ref{strictproduct} enjoys a universal characterization, whose proof is straightforward.

\begin{proposition}\label{productstrict}
Let $\tuple{\Sigma_1,\Mt_1}$ and $\tuple{\Sigma_2,\Mt_2}$ be PNmatrices. 
Their strict product $\tuple{\Sigma_1\cup\Sigma_2,\Mt_1\ast\Mt_2}$ is a product in $\cat{PNMatr}$, with projection homomorphisms  
$\pi_i:\tuple{\Sigma_1\cup\Sigma_2,\Mt_1\ast\Mt_2}\to\tuple{\Sigma_i,\Mt_i}$ for each $i\in\{1,2\}$.
\end{proposition}

It goes without saying that the relationships between valuations in the PNmatrices $\tuple{\Sigma_1,\Mt_1},\tuple{\Sigma_2,\Mt_2}$ and valuations in $\tuple{\Sigma_1\cup\Sigma_2,\Mt_1\ast\Mt_2}$ mediated by the projection homomorphisms that were presented in Section~\ref{sec:withstrictproduct} are simple consequences of the universal property enjoyed by products in a category.\smallskip

A dual construction is also possible. The key idea is that one can take crucial advantage of partiality to blend together PNmatrices. Let $\cM=\{\tuple{\Sigma,\Mt_i}:i\in I\}$ be a set of PNmatrices, each $\Mt_i=\tuple{V_i,D_i,\cdot_{\Mt_i}}$.
The \emph{sum} of $\cM$ is the PNmatrix $(\Sigma,\oplus\cM)$ where
$\oplus\cM=\tuple{V,D,\cdot_\oplus}$ is defined by $V=\bigcup_{i\in I} (\{i\}\times V_i)$, 
$D=\bigcup_{i\in I} (\{i\}\times D_i)$, and 
for $n\in\nats_0$ and $c\in\Sigma^{(n)}$,
$c_\oplus((i_1,x_1),\dots,(i_n,x_n))=\bigcup_{\stackrel{i\in I : i=i_1=\dots=i_n}{}}(\{i\}\times c_{\Mt_i}(x_1,\dots,x_n))$. 

It is clear that $(\Sigma,\oplus\cM)$ is a \emph{coproduct} of $\cM$ in~$\cat{PNMatr}$, 
with inclusion homomorphisms $\iota_i:\tuple{\Sigma,\Mt_i}\to\tuple{\Sigma,\oplus\cM}$ defined, for each $i\in I$ and each $x\in V_i$, by $\iota_i(x)=(i,x)$. Therefore, we have $\bigcup_{i\in I}\BVal(\Mt_i)\subseteq\BVal(\oplus\cM)$. Perhaps surprisingly, however, it may happen that 
$\BVal(\oplus\cM)\neq\bigcup_{i\in I}\BVal(\Mt_i)$.

\begin{example}\label{nobin}(A badly-behaved sum).\\
Let $\Sigma$ be the signature whose only connectives are $@\in\Sigma^{(0)}$ and $f\in\Sigma^{(1)}$, and consider the matrices $\tuple{\Sigma,\Mt_0}$ with $\Mt_0=\tuple{\{0\},\emptyset,\cdot_{\Mt_0}}$ and $@_{\Mt_0}=f_{\Mt_0}(0)=\{0\}$, and $\tuple{\Sigma,\Mt_2}$ with $\Mt_2=\tuple{\{0,1\},\{1\},\cdot_{\Mt_2}}$ and $@_{\Mt_2}=f_{\Mt_2}(0)=f_{\Mt_2}(1)=\{1\}$.

We have that $\BVal(\Mt_0)=\{{\tt 0}\}$ with ${\tt 0}:L_\Sigma(P)\to\{0,1\}$ such that ${\tt 0}(A)=0$ for every formula $A\in L_\Sigma(P)$. We also have that $\BVal(\Mt_2)=\{b\in\BVal(\Sigma):b(A)=1\textrm{ for every }A\in L_\Sigma(P)\setminus P\}$.

However, it is clear that $\BVal(\oplus\{\Mt_0,\Mt_2\})$ contains bivaluations which are neither in $\BVal(\Mt_0)$ nor in $\BVal(\Mt_2)$. For instance, given a variable $p\in P$, it is clear that 
$v:L_\Sigma(P)\to\{(0,0),(2,0),(2,1)\}$ such that $v(A)=(0,0)$ if $p$ occurs in A, and $v(A)=(2,1)$ otherwise, defines a valuation $v\in\Val(\oplus\{\Mt_0,\Mt_2\})$. Hence, $\BVal(\oplus\{\Mt_0,\Mt_2\})$ contains the bivaluation $b$ such that $b(A)=1$ if and only if $p$ does not occur in $A$, which is clearly not in $\BVal(\Mt_0)\cup\BVal(\Mt_2)$.\hfill$\triangle$
\end{example}

\subsubsection{The good, the bad, and the ugly}

Now, we will prove the (good, very good) fact that `almost' every logic (in the single or multiple-conclusion sense, it does not matter) can be characterized by a single PNmatrix (actually, a Pmatrix). This property is striking, as it is well known to fail for matrices, or even Nmatrices (see~\cite{finval,SS,Wojcicki88}), and at the same time reinforces the wide range of applicability of our results. The bad and the ugly are actually both related to the `almost' part of our statement. On the one hand, we need to understand the reason for the exception, which actually lies on the (bad) fact that valuations on PNmatrices cannot be locally assessed for a sublanguage (or subsignature). On the other hand, the nature of the exception is related to a rather annoying (and mathematically ugly) syntactic property. Both are suitably illustrated by Example~\ref{nobin}, as the the crucial reason for the behaviour shown in this example is the absence of a 2-place connective (actually, of any connective with at least two places).

\begin{lemma}\label{badugly}
 Let $\cM=\{\tuple{\Sigma,\Mt_i}:i\in I\}$ be a set of PNmatrices. If $\Sigma^{(n)}\neq\emptyset$ for some $n>1$
then $\BVal(\oplus\cM)=\bigcup_{i\in I}\BVal(\Mt_i)$.
\end{lemma}
\proof{Observe that if $v\in\Val(\oplus\cM)$ and $A,B\in L_\Sigma(P)$ are such that $v(A)=(i,x)$ and $v(B)=(j,y)$ then it must be the case that $i=j$. This is a consequence of the existence of a connective $c\in \Sigma^{(n)}$ for some $n>1$, as $v(c(A,B,\dots,B))\in c_\oplus(v(A),v(B),\dots,v(B))=c_\oplus((i,x),(j,y),\dots,(j,y))$ and, by definition, we have $c_\oplus((i,x),(j,y),\dots,(j,y))=\emptyset$ if $i\neq j$.

Hence, it is clear that if $v\in\Val(\oplus\cM)$ then letting $v_i(A)=x$ if $v(A)=(i,x)$ defines a valuation $v_i\in\Val(\Mt_i)$. As the compatibility of these values is granted by the definition of designated values in $\oplus\cM$, it follows that $\BVal(\oplus\cM)=\bigcup_{i\in I}\BVal(\Mt_i)$.\qed}\\

We can now take advantage of well-known results about logical matrices. 
Recall that a \emph{Lindenbaum matrix} over $\Sigma$ is a matrix of the form $\Mt_\Gamma=\tuple{L_\Sigma(P),\Gamma,\cdot_\Gamma}$ for some $\Gamma\subseteq L_\Sigma(P)$, with $c_\Gamma(A_1,\dots,A_n)=\{c(A_1,\dots,A_n)\}$ for every $n\in\nats_0$, $c\in\Sigma^{{n}}$, and $A_1,\dots,A_n\in L_\Sigma(P)$. Given a bivaluation $b:L_\Sigma(P)\to\{0,1\}$, we will use $\Mt_b$ to denote the Lindenbaum matrix 
$\Mt_{b^{-1}(1)}$.

In the single-conclusion case, we know that a logic $\tuple{\Sigma,\vdash}$ is precisely characterized by the set of its theories~\cite{Wojcicki88}. Hence, the \emph{Lindenbaum bundle} $\Lind(\tuple{\Sigma,\vdash})$ consists of the Lindenbaum matrices $\Mt_\Gamma$ over $\Sigma$ such that $\Gamma$ is a theory of $\tuple{\Sigma,\vdash}$ (when the logic is compact, we could as well consider only its relatively maximal theories).

In the multiple-conclusion case, a logic $\tuple{\Sigma,\der}$ is known to be precisely characterized by its maximal theory-pairs~\cite{SS,czelak2,Zygmunt}. Thus, the \emph{Lindenbaum bundle} $\Lind(\tuple{\Sigma,\der})$ contains precisely the Lindenbaum matrices $\Mt_\Gamma$ over $\Sigma$ such that $\Gamma\not\der (L_\Sigma(P)\setminus \Gamma)$ is a maximal theory-pair of $\tuple{\Sigma,\der}$ (see~\cite{blasio}).

Given a set of bivaluation $\cB\subseteq\BVal(\Sigma)$ closed under substitutions, we also define $\Lind(\tuple{\Sigma,\cB})$ as consisting of the Lindenbaum matrices $\Mt_b$ for every $b\in\cB$.\smallskip

The following result is an immediate consequence of Lemma~\ref{badugly}, taking into account the Pmatrix corresponding to summing the Lindenbaum bundle into consideration.

\begin{corollary}
Let $\tuple{\Sigma,\propto}$ be a logic. 

If $\Sigma^{(n)}\neq\emptyset$ for some $n>1$ then
 $\propto{=}\propto_{\oplus\Lind(\tuple{\Sigma,\propto})}$.
\end{corollary}

We should note that when considering signatures with no connectives with two or more places, it really may happen that the logic characterized by a set of PNmatrices does not coincide with the logic characterized by their coproduct.

\begin{example}\label{nobin2}(A badly-behaved sum, continued).\\
Recall Example~\ref{nobin}, with a signature $\Sigma$ whose only connectives are $@\in\Sigma^{(0)}$ and $f\in\Sigma^{(1)}$, and matrices $\tuple{\Sigma,\Mt_0}$ with $\Mt_0=\tuple{\{0\},\emptyset,\cdot_{\Mt_0}}$ such that $@_{\Mt_0}=f_{\Mt_0}(0)=\{0\}$, and $\tuple{\Sigma,\Mt_2}$ with $\Mt_2=\tuple{\{0,1\},\{1\},\cdot_{\Mt_2}}$ such that $@_{\Mt_2}=f_{\Mt_2}(0)=f_{\Mt_2}(1)=\{1\}$. Let $\cB=\BVal(\Mt_0)\cup\BVal(\Mt_2)=\{{\tt 0}\}\cup\{b\in\BVal(\Sigma):b(A)=1\textrm{ for every }A\in L_\Sigma(P)\setminus P\}$. Recall also that $\BVal(\oplus\{\Mt_0,\Mt_2\})$ contains the bivaluation $b$ such that $b(A)=1$ if and only if $p$ does not occur in $A$, which is clearly not in $\cB$.

It is not hard to see that given any PNmatrix $\tuple{\Sigma,\Mt}$ with $\Mt=\tuple{V,D,\cdot_{\Mt}}$ such that $\cB\subseteq\BVal(\Mt)$ it must also be the case that $b\in\BVal(\Mt)$. Indeed, since we have $\{{\tt 0},{\tt 1}\}\subseteq\cB$, there exist valuations $v_0,v_1\in\Val(\Mt)$ such that $v_0(A)\notin D$ and $v_1(A)\in D$ for every $A\in L_\Sigma(P)$. Therefore, the valuation such that 
$$v(A)= \begin{cases}
    v_0(A) & \textrm{if $p$ occurs in $A$} \\
    v_1(A) & \textrm{otherwise}
      \end{cases}
$$
is also in $\Val(\Mt)$. Just note that $v(@)=v_1(@)\in @_\Mt$, and that  $v(f(A))\in f_\Mt(v(A))$, because both $v_0$ and $v_1$ are in $\Val(\Mt)$ and $p$ occurs in $A$ if and only if $p$ occurs in $f(A)$. We conclude that $b\in\BVal(\Mt)$.

It turns out that it is impossible for a PNmatrix to have $\BVal(\Mt)=\cB$, and thus necessarily $\der_\Mt{\neq}\der_\cB$ for every PNmatrix $\tuple{\Sigma,\Mt}$.

Moreover, it is easy to see that $\cB$ is meet-closed, and that if $\cB'\subseteq\cB$ is such that $\cB'$ is closed under substitutions and $\cB'^\cap=\cB$ then one must have $\cB'=\cB$. Thus,  it is impossible for a PNmatrix to have $\BVal(\Mt)^\cap=\cB$, and one can also conclude that $\vdash_\Mt{\neq}\vdash_\cB$ for every PNmatrix $\tuple{\Sigma,\Mt}$.

For the sake of closure, we should note that $\der_\cB{=}\der_R$ and $\vdash_\cB{=}\vdash_R$ where $R$ contains the two rules
$$\frac{\;p\;}{\;f(q)\;}\qquad\qquad\frac{\;p\;}{\;\;\;@\;\;\;}$$
\noindent which are both single-conclusioned rules simply because $\cB$ is meet-closed. Indeed, we have that ${\Gamma}\der_\cB{\Delta}$ if and only if $\Gamma\vdash_\cB A$ for some $A\in \Delta$, and additionally $\Gamma\vdash_\cB A$ if and only if $\Gamma\neq\emptyset$ and $A\in L_\Sigma(P)\setminus P$, i.e., $A$ is not a variable.
\hfill$\triangle$
\end{example}

\subsection{Multiple-conclusion combination}

The results of Section~\ref{sect:semcomblog}, in particular Proposition~\ref{bivmultiple}, Lemma~\ref{strictvals} and Theorem~\ref{multiplefibring}, allowed us to characterize the combination of the multiple-conclusion logics defined by two PNmatrices as the logic defined by the intersection of their induced bivaluations, or equivalently as the logic of their strict product. Our aim is to provide a categorial explanation for these results.\smallskip

\begin{center}
\begin{tikzpicture}[-latex,semithick,
%estado/.style ={circle, top color=white, draw, text=blue, minimum width=2mm},
inicio/.style={fill=none,draw=none}
]
\node (pnmatr)    at    (0,0)    {$\cat{PNMatr}$};
\node (rexp)    at    (3,0)    {$\cat{Rexp}$};
\node (biv)    at    (6,0)    {$\cat{Biv}$};
\node (mult)  at (9,0) {$\cat{Mult}^{\cat{op}}$};
\node (adj) at (4.5,0) {$\top$};
\path (pnmatr) edge node[above] {Q} (rexp);
\path (rexp) edge [bend left] node[above] {\BVal} (biv);
\path (biv) edge [bend left] node[below] {$\Lind_\oplus$} (rexp);
\path (biv) edge node[below] {$\cong$} node[above]{$\textrm{Mult}$} (mult);

\end{tikzpicture}
\end{center}

In the terminology of~\cite{rexpansions}, that we extend here to PNmatrices, 
the existence of a {strict homomorphism} $h:\tuple{\Sigma,\Mt}\to\tuple{\Sigma_0,\Mt_0}$ means that $\tuple{\Sigma,\Mt}$ is a \emph{rexpansion} of $\tuple{\Sigma_0,\Mt_0}$. For simplicity, we will consider the quotient of 
$\cat{PNMatr}$ obtained by identifying all strict homomorphisms between the same two PNmatrices, thus obtaining a thin category that we will dub $\cat{Rexp}$. Equivalently, $\cat{Rexp}$ consists of the preordered class of all PNmatrices under the rexpansion relation, and we write $\tuple{\Sigma,\Mt}\sqsubseteq\tuple{\Sigma_0,\Mt_0}$ precisely when $\tuple{\Sigma,\Mt}$ is a {rexpansion} of $\tuple{\Sigma_0,\Mt_0}$. The obvious quotient functor $\textrm{Q}:\cat{PNMatr}\to\cat{Rexp}$ that sends each PNmatrix to itself is continuous and cocontinuous. Therefore, 
it follows that $\tuple{\Sigma,\oplus\cM}=\bigsqcup_{\tuple{\Sigma,\Mt}\in\cM}\tuple{\Sigma,\Mt}$ is a (non-unique) join in $\cat{Rexp}$, 
and we can say that $\oplus\cM$ is the least PNMatrix of which all PNmatrices in $\cM$ are rexpansions.
Dually, $\tuple{\Sigma_1\cup\Sigma_2,\Mt_1\ast\Mt_2}=\tuple{\Sigma_1,\Mt_1}\sqcap \tuple{\Sigma_2,\Mt_2}$ is a (non-unique) meet in $\cat{Rexp}$, and we can say that the strict product of two PNmatrices is the largest PNmatrix that is a rexpansion of both $\tuple{\Sigma_1,\Mt_1}$ and $\tuple{\Sigma_2,\Mt_2}$.\smallskip

We need also to consider the posetal category $\cat{Biv}$, consisting of all pairs $\tuple{\Sigma,\cB}$ where $\Sigma$ is a signature and $\cB\subseteq\BVal(\Sigma)$ is closed under substitutions, partially ordered by the relation defined as $\tuple{\Sigma,\cB}\sqsubseteq\tuple{\Sigma_0,\cB_0}$ if $\Sigma_0\subseteq\Sigma$ and $\cB\subseteq\cB_0^\Sigma$. The mapping $\BVal$ extends to a (order-preserving) functor $\BVal:\cat{Rexp}\to \cat{Biv}$ such that $\BVal(\tuple{\Sigma,\Mt})=\tuple{\Sigma,\BVal(\Mt)}$. If we further define $\Lind_\oplus: \cat{Biv}\to\cat{Rexp}$ by $\Lind_\oplus(\tuple{\Sigma,\cB})=\tuple{\Sigma,\oplus\Lind(\cB)}$ we also obtain a (order-preserving) functor, which is actually left-adjoint to $\BVal$.

\begin{proposition}\label{propmult}
The functors $\Lind_\oplus,\BVal$ constitute a Galois connection, that is, for every  $\tuple{\Sigma,\cB}$ in $\cat{Biv}$ and every $\tuple{\Sigma_0,\Mt_0}$ in $\cat{Rexp}$, the following conditions are equivalent:
\begin{itemize}
\item $\Lind_\oplus(\tuple{\Sigma,\cB})\sqsubseteq\tuple{\Sigma_0,\Mt_0}$,
\item $\tuple{\Sigma,\cB}\sqsubseteq\BVal(\tuple{\Sigma_0,\Mt_0})$.
\end{itemize}
\end{proposition}
\proof{First assume that $\Lind_\oplus(\tuple{\Sigma,\cB})\sqsubseteq\tuple{\Sigma_0,\Mt_0}$, i.e., there exists a strict homomorphism $h:\tuple{\Sigma,\oplus\{\Mt_b:b\in\cB\}}\to\tuple{\Sigma_0,\Mt_0}$. Obviously, 
we have $\Sigma_0\subseteq\Sigma$. Given $b\in \cB$, we have the inclusion homomorphism $\iota_b:\tuple{\Sigma,\Mt_b}\to\tuple{\Sigma,\oplus\{\Mt_b:b\in\cB\}}$. Composing, we have $\tuple{\Sigma,\Mt_b}\sqsubseteq\tuple{\Sigma_0,\Mt_0}$ and therefore $b\in\BVal(\Mt_b)\subseteq\BVal(\Mt_0^\Sigma)=\BVal(\Mt_0)^\Sigma$. We conclude that 
$\cB\subseteq\BVal(\Mt_0)^\Sigma$, and so
$\tuple{\Sigma,\cB}\sqsubseteq\BVal(\tuple{\Sigma_0,\Mt_0})$.

Conversely, assume that $\tuple{\Sigma,\cB}\sqsubseteq\BVal(\tuple{\Sigma_0,\Mt_0})$. This means that not only 
$\Sigma_0\subseteq\Sigma$ but also $\cB\subseteq\BVal(\Mt_0)^\Sigma$. The latter inclusion implies that for each $b\in\cB$ there exists $v_b\in\Val(\Mt_0)$ such that $b=t\circ v_b\circ \skel$, and thus $(v_b\circ \skel):\tuple{\Sigma,\Mt_b}\to\tuple{\Sigma_0,\Mt_0}$ is a strict homomorphism. The universal property of $\oplus\{\Mt_b:b\in\cB\}$ then implies that there exists an homomorphism $h:\tuple{\Sigma,\oplus\{\Mt_b:b\in\cB\}}\to\tuple{\Sigma_0,\Mt_0}$ (actually a single one with the property that $h\circ\iota_b=v_b\circ\skel$ for every $b\in \cB$). We conclude that $\Lind_\oplus(\tuple{\Sigma,\cB})\sqsubseteq\tuple{\Sigma_0,\Mt_0}$.\qed}\\

Consequently, $\BVal$ preserves meets (limits) and we rediscover Lemma~\ref{strictvals}, as then 
$\tuple{\Sigma_1\cup\Sigma_2,\BVal(\Mt_1\ast\Mt_2)}=\BVal(\tuple{\Sigma_1\cup\Sigma_2,\Mt_1\ast\Mt_2})=\BVal(\tuple{\Sigma_1,\Mt_1}\sqcap\tuple{\Sigma_2,\Mt_2})=\BVal(\tuple{\Sigma_1,\Mt_1})\sqcap\BVal(\tuple{\Sigma_2,\Mt_2}){=}\tuple{\Sigma_1,\BVal(\Mt_1)}\sqcap\tuple{\Sigma_1,\BVal(\Mt_2)}$, which equals $\tuple{\Sigma_1\cup\Sigma_2,\BVal(\Mt_1)^{\Sigma_1\cup\Sigma_2}\cap\BVal(\Mt_2)^{\Sigma_1\cup\Sigma_2}}=
\tuple{\Sigma_1\cup\Sigma_2,\calB^{\mult}_{12}}$.\smallskip

Finally, we consider another posetal category $\cat{Mult}$, of multiple-conclusion logics ordered by inclusion, that is, $\tuple{\Sigma_0,\der_0}\sqsubseteq\tuple{\Sigma,\der}$ if $\Sigma_0\subseteq\Sigma$ and $\der_0{\subseteq}\der$. It is clear that a combined logic $\tuple{\Sigma_1,\der_1}\bullet\tuple{\Sigma_2,\der_2}$ is a join $\tuple{\Sigma_1,\der_1}\sqcup\tuple{\Sigma_2,\der_2}$ in $\cat{Mult}$. The mapping $\textrm{Mult}:\cat{Biv}\to\cat{Mult}$ such that $\textrm{Mult}(\tuple{\Sigma,\cB})=\tuple{\Sigma,\der_\cB}$ establishes an order isomorphism between $\cat{Biv}$ and $\cat{Mult}^{\cat{op}}$, as we show next.

\begin{proposition}
$\textrm{Mult}:\cat{Biv}\to\cat{Mult}$ is a dual order isomorphism, that is:
\begin{itemize}
\item $\textrm{Mult}$ is bijective, and
\item $\tuple{\Sigma,\cB}\sqsubseteq\tuple{\Sigma_0,\cB_0}$ if and only if $\tuple{\Sigma_0,\der_{\cB_0}}\sqsubseteq\tuple{\Sigma,\der_\cB}$.
\end{itemize}
\end{proposition}
\proof{$\textrm{Mult}$ is bijective precisely because for every multiple-conclusion logic $\tuple{\Sigma,\der}$ there exists only one set $\cB\subseteq\BVal(\Sigma)$ such that $\der{=}\der_\cB$. We know from~\cite{SS} that $\cB=\{b\in \BVal(\Sigma): b^{-1}(1)\not\der b^{-1}(0)\}$, which relates to the set of all maximal theory-pairs of the logic. 

As a consequence, just note that $\der_{\cB_0}{\subseteq}\der_\cB$ is equivalent to having that, for every $\Gamma,\Delta\subseteq L_\Sigma(P)$, $\Gamma\not\der_\cB\Delta$ implies $\skel(\Gamma)\not\der_{\cB_0}\skel(\Delta)$. On its turn, the later is equivalent to having, for every $b\in\cB$, that $\skel(b^{-1}(1))\not\der_{\cB_0}\skel(b^{-1}(0))$, which actually means that for every $b\in\cB$ there exists $b_0\in\cB_0$ such that $b=b_0\circ \skel$, or simply that $\cB\subseteq\cB_0^\Sigma$.\qed}\\

As a consequence, we can  recover Proposition~\ref{bivmultiple}, because $\tuple{\Sigma_1\cup\Sigma_2,\der_{\cB_{12}^{\mult}}}{=}$ $\textrm{Mult}(\tuple{\Sigma_1\cup\Sigma_2,\cB_{12}^{\mult}}){=}
\textrm{Mult}(\tuple{\Sigma_1,\cB_1}\sqcap\tuple{\Sigma_2,\cB_2})$ which, by duality, is equal to $\textrm{Mult}(\tuple{\Sigma_1,\cB_1})\sqcup\textrm{Mult}(\tuple{\Sigma_2,\cB_2}){=}\tuple{\Sigma_1,\der_{\cB_1}}\sqcup\tuple{\Sigma_2,\der_{\cB_2}}{=}\tuple{\Sigma_1,\der_{\cB_1}}\bullet\tuple{\Sigma_2,\der_{\cB_2}}$.\smallskip

Theorem~\ref{multiplefibring} can also be recovered, simply, by noting that $\tuple{\Sigma_1\cup\Sigma_2,\der_{\Mt_1\ast \Mt_2}} {=}$ $ \textrm{Mult}(\tuple{\Sigma_1\cup\Sigma_2,\BVal(\Mt_1\ast\Mt_2)})
{=}\textrm{Mult}(\BVal(\tuple{\Sigma_1,\Mt_1}) {\;\sqcap\;} \BVal(\tuple{\Sigma_2,\Mt_2}))$, which is equal to 
$\textrm{Mult}(\BVal(\tuple{\Sigma_1,\Mt_1}))\sqcup\textrm{Mult}(\BVal(\tuple{\Sigma_2,\Mt_2})){=}
\tuple{\Sigma_1,\der_{\Mt_1}}\sqcup\tuple{\Sigma_2,\der_{\Mt_2}}{=}$ $
\tuple{\Sigma_1,\der_{\Mt_1}}\bullet\tuple{\Sigma_2,\der_{\Mt_2}}$.

\subsection{Single-conclusion combination}

The results of Section~\ref{sect:semcomblog}, particularly Proposition~\ref{bivsingle} and 
Theorems~\ref{saturatedfibring} and~\ref{singlefibring}, allowed us to characterize the combination of the single-conclusion logics defined by two PNmatrices as the logic defined by the intersection of the meet-closure of their induced bivaluations, or equivalently as the logic of the strict product of their $\omega$-powers (or the PNmatrices themselves, when saturated). In order to explain these results categorially, we can adopt a similar strategy.

\begin{center}
\begin{tikzpicture}[-latex,semithick,
inicio/.style={fill=none,draw=none}
]
\node (pnmatr)    at    (0,0)    {$\cat{SPNMatr}$};
\node (rexp)    at    (3,0)    {$\cat{SRexp}$};
\node (biv)    at    (6,0)    {$\cat{Biv}^\cap$};
\node (mult)  at (9,0) {$\cat{Sing}^{\cat{op}}$};
\node (adj) at (4.5,0) {$\top$};
\path (pnmatr) edge node[above] {Q} (rexp);
\path (rexp) edge [bend left] node[above] {$\BVal^{+{\tt 1}}$} (biv);
\path (biv) edge [bend left] node[below] {$\Lind_\oplus^{-{\tt 1}}$} (rexp);
\path (biv) edge node[below] {$\cong$} node[above]{$\textrm{Sing}$} (mult);

\end{tikzpicture}
\end{center}

Given the properties of $\omega$-powers, as stated in Lemma~\ref{wpower}, we can restrict our attention to the full subcategory $\cat{SPNMatr}$ of $\cat{PNMatr}$ whose objects are just the saturated PNmatrices. We apply the same restriction to obtain the full subcategory $\cat{SRexp}$ of $\cat{Rexp}$. As we know that the strict product of saturated PNmatrices is still saturated, the operation will still correspond to a (non-unique) meet in $\cat{SRexp}$. Interestingly, however, coproducts of saturated matrices need not be saturated.

\begin{example}\label{nosat}(Sums and saturation).\\
Let $\Sigma$ be the signature whose only connectives are $@\in\Sigma^{(0)}$ and $f\in\Sigma^{(1)}$, and consider the PNmatrix $\tuple{\Sigma,\Mt_\emptyset}$ with $\Mt_\emptyset=\tuple{\{0,1\},\{1\},\cdot_{\Mt_\emptyset}}$ such that $@_{\Mt_\emptyset}=\emptyset$, $f_{\Mt_\emptyset}(0)=\{1\}$, $f_{\Mt_\emptyset}(1)=\{0\}$, and the matrix $\tuple{\Sigma,\Mt_1}$ with $\Mt_1=\tuple{\{1\},\{1\},\cdot_{\Mt_1}}$ such that $@_{\Mt_1}=f_{\Mt_1}(1)=\{1\}$.

We have that $\BVal(\Mt_\emptyset)=\emptyset$, and 
$\BVal(\Mt_1)=\{{\tt 1}\}$. Since $\emptyset^\cap=\{{\tt 1}\}^\cap=\{{\tt 1}\}$ we can conclude from Lemma~\ref{theone} that both PNmatrices are saturated.

However, the PNmatrix $\oplus\{\Mt_\emptyset,\Mt_1\}$ is not saturated. In fact, it suffices to show that 
$\BVal(\oplus\{\Mt_\emptyset,\Mt_1\})$, which clearly contains the bivaluation ${\tt 1}$, is not meet-closed. 
To see this, fix a variable $p\in P$ and note that $v_1,v_2:L_\Sigma(P)\to\{(\emptyset,0),(\emptyset,1),(1,1)\}$ such that $$v_1(A)=\begin{cases}
    (\emptyset,n\,\text{mod }2) & \textrm{if $A=f^n(p)$} \\
    (1,1) & \textrm{otherwise}
      \end{cases}
      \textrm{, }
      v_2(A)=\begin{cases}
    (\emptyset,n+1\,\text{mod }2) & \textrm{if $A=f^n(p)$} \\
    (1,1) & \textrm{otherwise}
      \end{cases}$$

\noindent are valuations $v_1,v_2\in\Val(\oplus\{\Mt_\emptyset,\Mt_1\})$. Hence, $b_1,b_2:L_\Sigma(P)\to\{0,1\}$ such that $b_1(A)=0$ if and only if $A=f^n(p)$ with $n$ even, and $b_2(A)=0$ if and only if $A=f^n(p)$ with $n$ odd, both are bivaluations $b_1,b_2\in\BVal(\oplus\{\Mt_\emptyset,\Mt_1\})$. Letting $X=\{b_1,b_2\}$ we have $b_X$ such that $b_X(A)=0$ if and only if $p$ occurs in $A$. It is not difficult to conclude that $b_X\notin \BVal(\oplus\{\Mt_\emptyset,\Mt_1\})$, noting that $(\emptyset,0)$ is the only value not designated and that
any valuation $v\in\Val(\oplus\{\Mt_\emptyset,\Mt_1\})$ such that 
$v(p)=(\emptyset,0)$ must have $v(f(p))\in f_\oplus(v(p))=f_\oplus(\emptyset,0)=\{(\emptyset,1)\}$.
\hfill$\triangle$
\end{example}

Following the lines developed in Section~\ref{sect:semcomblog}, we shall also consider 
the full subcategory $\cat{Biv}^\cap$ of $\cat{Biv}$ whose objects are meet-closed sets of bivaluations. 
The functor $\BVal^{+{\tt 1}}:\cat{SRexp}\to \cat{Biv}^\cap$ is such that $\BVal^{+{\tt 1}}(\tuple{\Sigma,\Mt})=\tuple{\Sigma,\BVal(\Mt)\cup\{{\tt 1}\}}$. The functor $\Lind^{-{\tt 1}}_\oplus: \cat{Biv}^\cap\to\cat{SRexp}$ 
is defined by 
$\Lind^{-{\tt 1}}_\oplus(\tuple{\Sigma,\cB})=\tuple{\Sigma,\oplus\Lind(\cB\setminus\{{\tt 1}\})}$. The following result shows that, despite of Example~\ref{nosat}, $\Lind^{-{\tt 1}}_\oplus$ is well-defined.

\begin{lemma}
If $\Sigma$ is a signature and $\cB\subseteq\BVal(\Sigma)$ is closed under substitutions and meet-closed then $\Lind^{-{\tt 1}}_\oplus(\tuple{\Sigma,\cB})$ is saturated.
\end{lemma}
\proof{For simplicity, let $\Lind^{-{\tt 1}}_\oplus(\tuple{\Sigma,\cB})=\tuple{\Sigma,\Mt}$.

The result is immediate if $\Sigma$ contains a connective with two or more places. When that is the case, Lemma~\ref{badugly} guarantees that $\BVal(\Mt)=\bigcup_{b\in \cB\setminus\{{\tt 1}\}}\BVal(\Mt_b)$. As one can easily check, for Lindenbaum matrices, 
$\BVal(\Mt_b)$ is simply $\{b\}$ closed under substitutions,
and thus $b\in\BVal(\Mt_b)\subseteq\cB$ since $\cB$ is itself closed under substitutions. Therefore, we have that $ \cB\setminus\{{\tt 1}\}\subseteq\bigcup_{b\in \cB\setminus\{{\tt 1}\}}\BVal(\Mt_b)\subseteq\cB$, and by Lemma~\ref{theone}, since $\cB$ is meet-closed, we can conclude that $\Mt$ is saturated.

Things are less straightforward when $\Sigma^{(n)}=\emptyset$ for every $n>1$, since we know that the set $\BVal(\Mt)$ may include bivaluations not in $\cB$. Indeed, 
letting
$\atm(A)\in \Atm=(P\cup\Sigma^{(0)})$ be the only atomic subformula of $A$, one has $b\in\BVal(\Mt)$ if and only if there exists a function $\gamma:\Atm\to(\cB\setminus\{{\tt 1}\})$ such that $b(A)=\gamma(\atm(A))(A)$ for every formula $A\in L_\Sigma(P)$. For convenience, we use $b_\gamma$ instead of just $b$ to denote each such substitution
\footnote{To make it more explicit, this means that each valuation $v\in\Val(\Mt)$ can be made to correspond with 
choosing a bivaluation $\gamma(t)=b_t\in\cB\setminus\{{\tt 1}\}$ for each $t\in\Atm$, and setting $v(A)=(b_{\atm(A)},A)$ 
(which is designated precisely when $b_{\atm(A)}(A)=1$) for each formula $A$. 

A similar characterization would apply, \emph{mutatis mutandis}, to arbitrary coproducts of (total) Nmatrices over a signature without connectives with two or more places.}. 

To establish that $\tuple{\Sigma,\Mt}$ is saturated, using again Lemma~\ref{theone}, it is sufficient to show that given $X=\{b_{\gamma_i}:i\in I\}\subseteq\BVal(\Mt)$ if the meet $b_X\neq{\tt 1}$ then $b_X\in\BVal(\Mt)$. First note that if $b_X\neq{\tt 1}$ then it is the case that $I\neq\emptyset$, and we can fix $j\in I$. For each $t\in \Atm$ consider $X_t=\{\gamma_i(t):i\in I\}\subseteq\cB\setminus\{{\tt 1}\}$. As $\cB$ is meet-closed, the meet $b_{X_t}\in\cB$, but it may still happen that $b_{X_t}={\tt 1}$. Define $\gamma:\Atm\to(\cB\setminus\{{\tt 1}\})$ by
$$\gamma(t)=\begin{cases}
b_{X_t} & \textrm{if $b_{X_t}\neq{\tt 1}$}\\
\gamma_j(t) & \textrm{otherwise}
\end{cases}.$$
We claim that $b_X=b_\gamma$. To see this, consider a formula $A$ and let $t=\atm(A)$.

If $b_{X_t}\neq{\tt 1}$ then $b_\gamma(A)=\gamma(t)(A)=b_{X_t}(A)=1$ if and only if $b_{\gamma_i}(A)=\gamma_i(t)(A)=1$ for every $i\in I$ if and only if $b_X(A)=1$.

If, on the contrary, we have $b_{X_t}={\tt 1}$ this means that $\gamma_i(t)={\tt 1}$ for every $i\in I$. Therefore, it follows that $b_\gamma(A)=\gamma(t)(A)=\gamma_j(t)(A)=1=b_X(A)$, since $b_{\gamma_i}(A)=\gamma_i(t)(A)=1$ for every $i\in I$.\qed}\\

Similarly, we have that  $\Lind^{-{\tt 1}}_\oplus$ is left-adjoint to $\BVal^{+{\tt 1}}$.

\begin{proposition}
The functors $\Lind^{-{\tt 1}}_\oplus,\BVal^{+{\tt 1}}$ constitute a Galois connection, that is, for every  $\tuple{\Sigma,\cB}$ in $\cat{Biv}^\cap$ and every $\tuple{\Sigma_0,\Mt_0}$ in $\cat{SRexp}$, the following conditions are equivalent:
\begin{itemize}
\item $\Lind^{-{\tt 1}}_\oplus(\tuple{\Sigma,\cB})\sqsubseteq\tuple{\Sigma_0,\Mt_0}$,
\item $\tuple{\Sigma,\cB}\sqsubseteq\BVal^{+{\tt 1}}(\tuple{\Sigma_0,\Mt_0})$.
\end{itemize}
\end{proposition}
\proof{Let $\cB\subseteq\BVal(\Sigma)$ be a meet-closed set of bivaluations closed under substitutions, and 
$\tuple{\Sigma_0,\Mt_0}$ be a saturated PNmatrix.
The result follows from Proposition~\ref{propmult}. Indeed, $\Lind^{-{\tt 1}}_\oplus(\tuple{\Sigma,\cB})=
\Lind_\oplus(\tuple{\Sigma,\cB\setminus\{{\tt 1}\}})\sqsubseteq\tuple{\Sigma_0,\Mt_0}$ if and only if 
$\tuple{\Sigma,\cB\setminus\{{\tt 1}\}}\sqsubseteq\BVal(\tuple{\Sigma_0,\Mt_0})$ in $\cat{Biv}$, which is equivalent to having $\tuple{\Sigma,\cB}\sqsubseteq\tuple{\Sigma,\BVal(\Mt_0)\cup\{{\tt 1}\}}=\BVal^{+{\tt 1}}(\tuple{\Sigma_0,\Mt_0})$ in $\cat{Biv}^\cap$.\qed}\\

Consequently, $\BVal^{+{\tt 1}}$ preserves meets (limits). Note also that, in the category $\cat{Biv}^\cap$, we have that $\tuple{\Sigma_1,\cB_1}\sqcap\tuple{\Sigma_2,\cB_2}=\tuple{\Sigma_1\cup\Sigma_2,\cB_1\cap \cB_2}$.\smallskip

Finally, we consider another posetal category $\cat{Sing}$, of single-conclusion logics ordered by inclusion, that is, $\tuple{\Sigma_0,\vdash_0}\sqsubseteq\tuple{\Sigma,\vdash}$ if $\Sigma_0\subseteq\Sigma$ and $\vdash_0{\subseteq}\vdash$. It is clear that a combined logic $\tuple{\Sigma_1,\vdash_1}\bullet\tuple{\Sigma_2,\vdash_2}$ is a join $\tuple{\Sigma_1,\vdash_1}\sqcup\tuple{\Sigma_2,\vdash_2}$ in $\cat{Sing}$. The mapping $\textrm{Sing}:\cat{Biv}^\cap\to\cat{Sing}$ such that $\textrm{Sing}(\tuple{\Sigma,\cB})=\tuple{\Sigma,\vdash_\cB}$ is again an order isomorphism, now between $\cat{Biv}^\cap$ and $\cat{Sing}^{\cat{op}}$.

\begin{proposition}
$\textrm{Sing}:\cat{Biv}^\cap\to\cat{Sing}$ is a dual order isomorphism, that is:
\begin{itemize}
\item $\textrm{Sing}$ is bijective, and
\item $\tuple{\Sigma,\cB}\sqsubseteq\tuple{\Sigma_0,\cB_0}$ if and only if $\tuple{\Sigma_0,\vdash_{\cB_0}}\sqsubseteq\tuple{\Sigma,\vdash_\cB}$.
\end{itemize}
\end{proposition}
\proof{$\textrm{Sing}$ is bijective precisely because for every single-conclusion logic $\tuple{\Sigma,\vdash}$ there exists a single meet-closed set of bivaluations $\cB\subseteq\BVal(\Sigma)$ such that $\der{=}\der_\cB$. We know from~\cite{Wojcicki88} that $\cB=\{b\in \BVal(\Sigma): b^{-1}(1)\textrm{ is a theory of }\tuple{\Sigma,\vdash}\}$. 

As a consequence, just note that $\vdash_{\cB_0}{\subseteq}\vdash_\cB$ 
if and only if $\cB\subseteq\cB_0^\cap=\cB_0$, just because $\cB_0$ is meet-closed and thus, equivalently, 
$\tuple{\Sigma,\cB}\sqsubseteq\tuple{\Sigma_0,\cB_0}$.\qed}\\

As a consequence, we  can recover Proposition~\ref{bivsingle} too, just because we have 
$\tuple{\Sigma_1\cup\Sigma_2,\vdash_{\cB_{12}^{\sing}}}{=}\textrm{Sing}(\tuple{\Sigma_1\cup\Sigma_2,\cB_{12}^{\sing}}){=}
\textrm{Sing}(\tuple{\Sigma_1\cup\Sigma_2,(\cB_{1}^{\Sigma_1\cup\Sigma_2})^\cap\cap(\cB_{2}^{\Sigma_1\cup\Sigma_2})^\cap}){=}$  

\noindent $\textrm{Sing}(\tuple{\Sigma_1,\cB_1^\cap}\sqcap\tuple{\Sigma_2,\cB_2^\cap})$ which, by duality, equals $\textrm{Sing}(\tuple{\Sigma_1,\cB_1^\cap})\,\sqcup\,\textrm{Sing}(\tuple{\Sigma_2,\cB_2^\cap}){=}$ 
$\tuple{\Sigma_1,\vdash_{\cB_1^\cap}}\sqcup\tuple{\Sigma_2,\vdash_{\cB_2^\cap}}{=}$ $\tuple{\Sigma_1,\vdash_{\cB_1}}\sqcup\tuple{\Sigma_2,\vdash_{\cB_2}}{=}\tuple{\Sigma_1,\vdash_{\cB_1}}\bullet\tuple{\Sigma_2,\vdash_{\cB_2}}$.\smallskip

Theorem~\ref{saturatedfibring} can also be recovered. If $\tuple{\Sigma_1,\Mt_1}$ and $\tuple{\Sigma_2,\Mt_2}$ are saturated PNmatrices, then note that $\tuple{\Sigma_1\cup\Sigma_2,\vdash_{\Mt_1\ast \Mt_2}} {=}
 \textrm{Sing}(\BVal^{+{\tt 1}}(\tuple{\Sigma_1\cup\Sigma_2,\Mt_1\ast\Mt_2}))%{=}\der(\BVal(\tuple{\Sigma_1\,\Mt_1}\sqcap \tuple{\Sigma_2\,\Mt_2}))
{=}$ $\textrm{Sing}(\BVal^{+{\tt 1}}(\tuple{\Sigma_1,\Mt_1}) {\;\sqcap\;} \BVal^{+{\tt 1}}(\tuple{\Sigma_2,\Mt_2}))$, which equals 
$\textrm{Sing}(\BVal^{+{\tt 1}}(\tuple{\Sigma_1,\Mt_1}))\sqcup\textrm{Sing}(\BVal^{+{\tt 1}}(\tuple{\Sigma_2,\Mt_2})){=}
\tuple{\Sigma_1,\vdash_{\Mt_1}}\sqcup\tuple{\Sigma_2,\vdash_{\Mt_2}}{=}$ $
\tuple{\Sigma_1,\vdash_{\Mt_1}}\bullet\tuple{\Sigma_2,\vdash_{\Mt_2}}$. Theorem~\ref{singlefibring} also follows, by the same argument, using Lemma~\ref{wpower}.

\section{Examples and applications}\label{secapps}

Besides illustrating examples, of which we have already shown a few, our aim is now to present concrete ways of applying the tools we have previously defined to problems pertaining to the modular conception and analysis of logics. 

\subsection{Adding axioms}\label{sub:axioms}

Often, one works with a single-conclusion logic which can then be streghtened by the addition of new axioms~\cite{rexpansions,taming,addax} (and possibly also new syntax). Concretely, let $\tuple{\Sigma_1,\vdash_1}$ be a single-conclusion logic, $\Sigma_2$ be a signature, $\Ax\subseteq L_{\Sigma_2}(P)$ be a set of \emph{axiom schemata}, and define $\Ax^{\inst}=\{A^\sigma:A\in\Ax\textrm{ and } \sigma:P\to L_{\Sigma_1\cup\Sigma_2}(P)\}$. The \emph{strengthening of $\tuple{\Sigma_1,\vdash_1}$ with the schema axioms $\Ax$} is the single-conclusion logic $\tuple{\Sigma_1\cup\Sigma_2,\vdash_1^{\Ax}}$ 
defined by $\Gamma\vdash_1^{\Ax}A$ if and only if $\Gamma\cup\Ax^{\inst}\vdash_1^{\Sigma_1\cup\Sigma_2} A$, for every $\Gamma\cup\{A\}\subseteq L_{\Sigma_1\cup\Sigma_2}(P)$. 
Our aim is to apply the ideas developed above in order to obtain a semantic characterization of $\tuple{\Sigma_1\cup\Sigma_2,\vdash_1^{\Ax}}$ from a given semantic characterization of $\tuple{\Sigma_1,\vdash_1}$.\smallskip

In terms of bivaluations, the following simple result from~\cite{softcomp} is instrumental.

\begin{proposition}\label{axiomsandvals}
Let $\Sigma_1,\Sigma_2$ be signatures, $\calB_1\subseteq\BVal(\Sigma_1)$ closed under substitutions, and $\Ax\subseteq L_{\Sigma_2}(P)$.
The single-conclusion logic $\tuple{\Sigma_1\cup\Sigma_2,\vdash_{\calB_1}^{\Ax}}$ is characterized by $\calB_1^{\Ax}=\{b\in\calB_1^{\Sigma_1\cup\Sigma_2}:b(\Ax^{\inst})\subseteq\{1\}\}$.
\end{proposition}
\proof{Given $\Gamma\cup\{A\}\subseteq L_{\Sigma_1\cup\Sigma_2}(P)$, note that $\Gamma\not\vdash_{\calB_1^{\Ax}}A$ if and only if there exists $b\in\calB_1^{\Ax}$ such that $b(\Gamma)\subseteq\{1\}$ and $b(A)=0$ if and only if there exists $b\in\calB_1^{\Sigma_1\cup\Sigma_2}$ such that 
$b(\Ax^{\inst}),b(\Gamma)\subseteq\{1\}$ and $b(A)=0$ if and only if 
$\Gamma\cup\Ax^{\inst}\not\vdash_{\calB_1}^{\Sigma_1\cup\Sigma_2}A$ if and only if $\Gamma\not\vdash_{\calB_1}^{\Ax}A$.
\qed}\\

Clearly, $\calB_1^{\Ax}=\calB_1^{\Sigma_1\cup\Sigma_2}\cap\{b\in\BVal(\Sigma_1\cup\Sigma_2):b(\Ax^{\inst})\subseteq\{1\}\}$. The set of bivaluations 
$\{b\in\BVal(\Sigma_1\cup\Sigma_2):b(\Ax^{\inst})\subseteq\{1\}\}$ is easily seen to be meet-closed. 
The result follows, though, even if $\calB_1$ is not meet-closed.\smallskip

With respect to PNmatrices, we can take advantage of non-determinism for building a simple semantics for the logic associated to the calculus whose rules are precisely $\frac{\;\emptyset\;}{A}$ for each $A\in\Ax$, also equivalently defined by the set of bivaluations $\{b\in\BVal(\Sigma_1\cup\Sigma_2):b(\Ax^{\inst})\subseteq\{1\}\}$. 

\begin{definition}\label{axioms}
The Nmatrix $\tuple{\Sigma_1\cup\Sigma_2,\Mt_{\Ax}}$ with $\Mt_{\Ax}=\tuple{V_{\Ax},D_{\Ax},\cdot_{\Ax}}$ is defined by

\begin{itemize}
\item $V_{\Ax}=\{(A,0):A\notin\Ax^{\inst}\}\cup\{(A,1):A\in L_{\Sigma_1\cup\Sigma_2}(P)\}$, 
\item $D_{\Ax}=\{(A,1):A\in L_{\Sigma_1\cup\Sigma_2}(P)\}$, and 
\item for every $n\in\nats_0$, $c\in({\Sigma_1\cup\Sigma_2})^{(n)}$ and $(A_1,{x_1}),\dots,(A_n,{x_n}) \in V_{\Ax}$, we let 
$${\small c_{\Ax}((A_1,{x_1}),\dots,(A_n,{x_n}))=\left\{\begin{array}{ll}
\{(A,1)\} & \textrm{if }A=c(A_1,\dots,A_n)\in\Ax^{\inst}\\
\{(A,0),(A,1)\} & \textrm{if }A=c(A_1,\dots,A_n)\notin\Ax^{\inst}
\end{array}\right..}$$
\end{itemize}
\end{definition}

It is easy to characterize the properties of this construction.

\begin{lemma}\label{bivsaxioms}
We have that $\BVal(\Mt_{\Ax})=\{b\in\BVal(\Sigma_1\cup\Sigma_2):b(\Ax^{\inst})\subseteq\{1\}\}$.
\end{lemma}
\proof{If $v\in\Val(\Mt_{\Ax})$ and $A\in\Ax^{\inst}$ then necessarily $v(A)=(A,1)\in D_{\Ax}$. Therefore, we have $\BVal(\Mt_{\Ax})\subseteq\{b\in\BVal(\Sigma_1\cup\Sigma_2):b(\Ax^{\inst})\subseteq\{1\}\}$.

If $b\in\BVal(\Sigma_1\cup\Sigma_2)$ is such that $b(\Ax^{\inst})\subseteq\{1\}$, then we can build a compatible valuation $v\in\Val(\Mt_{\Ax})$ by letting 
$v(A)=(A,{b(A)})$ for each $A\in L_{\Sigma_1\cup\Sigma_2}(P)$. We conclude that $\{b\in\BVal(\Sigma_1\cup\Sigma_2):b(\Ax^{\inst})\subseteq\{1\}\}\subseteq
\BVal(\Mt_{\Ax})$.\qed}\\

Of course, it is simple to check that $\Gamma\vdash_{\Mt_{\Ax}} A$ if and only if $A\in\Gamma\cup\Ax^{\inst}$, and thus that 
$\tuple{\Sigma_1\cup\Sigma_2,\vdash_1^{\Ax}}=\tuple{\Sigma_1,\vdash_1}\bullet
\tuple{\Sigma_1\cup\Sigma_2,\vdash_{\Mt_{\Ax}}}$. It is worth noting too that $\tuple{\Sigma_1\cup\Sigma_2,\vdash_{\Mt_{\Ax}}}$ is saturated.

\begin{theorem}\label{adding}
The strengthening with $\Ax$ of the single-conclusion logic characterized by a PNmatrix is the single-conclusion logic characterized by its strict product with $\Mt_{\Ax}$, that is, 
given a PNmatrix $\tuple{\Sigma_1,\Mt_1}$,  we have that
$\tuple{\Sigma_1\cup\Sigma_2,\vdash_{\Mt_1}^{\Ax}}=\tuple{\Sigma_1\cup\Sigma_2,\vdash_{\Mt_1\ast \Mt_{\Ax}}}$.
\end{theorem}
\proof{The result follows directly from Proposition~\ref{axiomsandvals}, and Lemma~\ref{strictvals}~and~\ref{bivsaxioms}. 

With $\calB_1=\BVal(\Mt_1)$,  note that we necessarily have
$\BVal(\Mt_1\ast\Mt_{\Ax})=\calB_1^{\Sigma_1\cup\Sigma_2}\cap\BVal(\Mt_{\Ax})
=\calB_1^{\Ax}$.\qed}\\

Note that this construction, though similar to the ones presented in the previous section, has a remarkable difference: the PNmatrix $\Mt_1$ does not need to be saturated. 
The proof technique is essentially the same, but takes advantage of the fact that any set extending a theory of $\tuple{\Sigma_1\cup\Sigma_2,\vdash_{\Mt_{\Ax}}}$ is still a theory, as it retains all the instances of axioms. The end result, though, is still infinite-valued (denumerable, provided $\Mt_1$ is countable), but has really interesting consequences. 

For instance, we have seen in Example~\ref{platypuses} that the logics associated to the rules of \emph{modus ponens} and/or \emph{necessitation},

$$\frac{p\,,\, p \rightsquigarrow q}{q}\qquad\qquad\frac{p}{\square p}$$

\noindent can be given very simple non-deterministic two-valued semantics. Hence, according to Theorem~\ref{adding}, every logic obtained from these by the addition of axioms can be characterized by a denumerable PNmatrix, which of course includes every (global) modal logic, normal or not. A most interesting consequence of this idea is shown in the following example.

\begin{example}\label{intuition}(Intuitionistic logic).\\
Fix a signature $\Sigma_\mathsf{int}$ containing the desired set of intuitionistic propositional connectives, including implication, i.e., ${\to}\in\Sigma^{(2)}_\mathsf{int}$. Also, fix a set of axioms $\Int$ which together with the rule of \emph{modus ponens} constitutes a calculus for intuitionistic propositional logic $\tuple{\Sigma_\mathsf{int},\vdash_\mathsf{int}}$ ($\Int$ could consist of just the first two axioms in the axiomatization of $\vdash^\to_\mathsf{cls}$, in case implication is the only connective of $\Sigma_\mathsf{int}$).

Take a signature $\Sigma$ whose only connective is implication, that is, ${\to}\in\Sigma^{(2)}$ and consider the $\Sigma$-Nmatrix $\MPt=\tuple{\{0,1\},\{1\},\cdot_\MPt}$ where $\to_\MPt=\rightsquigarrow_{\TWO'}$ as defined in Example~\ref{platypuses}. That is to say that $\tuple{\Sigma,\vdash_\MPt}$ is the logic axiomatized by the rule of \emph{modus ponens}.

Obviously, $\vdash^{\Int}_\MPt{=}\vdash_\mathsf{int}$, and Theorem~\ref{adding} tells us that it is characterized by the PNmatrix $\MPt\ast\Mt_{\Int}$. 
Since $\MPt$ is finite and $\Mt_{\Int}$ is denumerable, we can conclude that $\MPt\ast\Mt_{\Int}$ is denumerable, which means that intuitionistic propositional logic
can be characterized by a single denumerable PNmatrix. This is a remarkable property of PNmatrices, witnessing their compression abilities, and contrasts with the known fact that intuitionistic logic cannot be characterized by a countable matrix~\cite{godel,wronski74,Wojcicki88}.  %
\hfill $\triangle$
\end{example}

Nicer, finite-valued, semantics can be obtained in particularly well-behaved scenarios, which (unsurprisingly) do not include the examples above. We refer to reader to~\cite{addax} for further details.

\subsection{Axiomatizability by splitting}

Obtaining symbolic calculi, or axiomatizations, for logics of interest, particularly if originally presented by semantic means, is well-known to be a non-trivial task, even harder when one seeks particularly well-behaved calculi. Herein, we show that to some extent, our results about combined semantics can have an impact on the possibility of splitting this task into the problem of obtaining calculi for suitably defined syntactic fragments of the logic, which can then put together, in a modular way, to produce a calculus for the original logic. Indeed, in abstract, given a logic $\cL$ we will be looking for ways to obtain logics $\cL_1,\cL_2$ such that $\cL_1\bullet\cL_2=\cL$.

Concretely, let $\tuple{\Sigma,\Mt}$ with $\Mt=\tuple{V,D,\cdot_\Mt}$ be a PNmatrix, and split the associated logical language into signatures $\Sigma_1,\Sigma_2$ such that $\Sigma=\Sigma_1\cup\Sigma_2$.
For $i\in\{1,2\}$, consider the component PNmatrix $\tuple{\Sigma_i,\Mt_i}$ where $\Mt_i=\tuple{V,D,\cdot_{\Mt_i}}$ is the reduct of $\Mt$ to the subsignature $\Sigma_i$, that is $c_{\Mt_i}=c_\Mt$ for every connective $c\in\Sigma_i$. Under which conditions can we obtain an axiomatization of $\tuple{\Sigma,\propto_\Mt}$ by just putting together axiomatizations of $\tuple{\Sigma_1,\propto_{\Mt_1}}$ and $\tuple{\Sigma_2,\propto_{\Mt_2}}$, in either the single and multiple-conclusion settings?\smallskip

It is useful to introduce the following notation. Given $\tuple{\Sigma,\Mt}$, a formula $A$ with $\var(A)\subseteq\{p_1,\dots,p_n\}$, and values $x_1,\dots,x_n\in V$, we define $A_\Mt(x_1,\dots,x_n)=\{v(A):v\in\Val(\Mt),v(p_1)=x_1,\dots,v(p_n)=x_n\}$. Given formulas $B_1,\dots,B_n$ we will also use $A(B_1,\dots,B_n)$ to denote the formula $A^\sigma$ where $\sigma$ is a substitution such that $\sigma(p_1)=B_1,\dots,\sigma(p_n)=B_n$.\smallskip

Recall from~\cite{SS,SYNTH,wollic19} that $\tuple{\Sigma,\Mt}$ is said to be \emph{monadic} provided that 
for every two values $x,y\in V$ with $x\neq y$ there exists a one-variable formula $S\in L_\Sigma(\{p\})$, to which we call a \emph{separator of $x$ and $y$}, such that $S_\Mt(x),S_\Mt(y)\neq\emptyset$ and $S_\Mt(x)\subseteq D$ and $S_\Mt(y)\cap D=\emptyset$, or vice-versa. When all the separators $S$ can be found in $L_{\Sigma_0}(\{p\})$ for a subsignature $\Sigma_0\subseteq\Sigma$, we say that the PNmatrix is \emph{$\Sigma_0$-monadic}.

\begin{lemma}\label{splitlemma}
If $\tuple{\Sigma,\Mt}$ is $(\Sigma_1\cap\Sigma_2)$-monadic then $\BVal(\Mt)=\BVal(\Mt_1\ast\Mt_2)$.
\end{lemma}
\proof{Taking advantage of rexpansions, and the universal property of strict products from Proposition~\ref{productstrict}, it is clear that 
the identity function on $V$ establishes strict homomorphisms $i_1:\tuple{\Sigma,{\Mt}}\to\tuple{\Sigma_1,{\Mt_1}}$ and $i_2:\tuple{\Sigma,{\Mt}}\to\tuple{\Sigma_2,{\Mt_2}}$, and thus there exists a strict homomorphism $j:\tuple{\Sigma,{\Mt}}\to\tuple{\Sigma,{\Mt_1\ast\Mt_2}}$ (letting $i_1(x)=i_2(x)=x$ and $j(x)=(x,x)$ for every $x\in V$). Therefore, we know $\BVal(\Mt)\subseteq\BVal(\Mt_1\ast\Mt_2)$.

{ 
To prove the converse inclusion, we need to take into account that $\Mt_1\ast\Mt_2$ has many spurious values.
Indeed, all pairs $(x,y)$ with $x\neq y$ are spurious. To see this, note that  if $x\neq y$ then there exists a separator $S\in L_{\Sigma_1\cap\Sigma_2}(\{p\})$ of $x$ and $y$. Then, $v(S(A))\in S_{\Mt_1\ast\Mt_2}(v(A))=S_{\Mt_1\ast\Mt_2}(x,y)\subseteq S_{\Mt_1}(x)\times S_{\Mt_2}(y)=S_{\Mt}(x)\times S_{\Mt}(y)$ which is impossible because the strict product does not have pairs of values where one is designated but not the other.
Hence, it is easy to see that $v':L_\Sigma(P)\to V$ such that $v'(A)=x$ if $v(A)=(x,x)$, for each formula $A$,  defines a valuation $v'\in\Val(\Mt)$ compatible to any possible valuation $v\in\Val(\Mt)$. We conclude that $\BVal(\Mt_1\ast\Mt_2)\subseteq\BVal(\Mt)$.\qed}}\\

We can now state our split axiomatization result for multiple-conclusion logics.

\begin{theorem}\label{splitmult}
If $R_i$ is an axiomatization of $\tuple{\Sigma_i,\der_{\Mt_i}}$, for each $i\in\{1,2\}$, then $R_1\cup R_2$ is an axiomatization of $\tuple{\Sigma,\der_\Mt}$, 
provided that $\Mt$ is $(\Sigma_1\cap\Sigma_2)$-monadic.

\end{theorem}
\proof{Clearly, using Theorem~\ref{multiplefibring} along with Proposition~\ref{fibringcalculi}, we have 
$\tuple{\Sigma,\der_{\Mt_1\ast\Mt_2}}=\tuple{\Sigma_1,\der_{\Mt_1}}\bullet\tuple{\Sigma_2,\der_{\Mt_2}}=\tuple{\Sigma_1,\der_{R_1}}\bullet\tuple{\Sigma_2,\der_{R_2}}=\tuple{\Sigma,\der_{R_1\cup R_2}}$, that is, $R_1\cup R_2$ is an axiomatization of the multiple-conclusion logic characterized by the strict product $\Mt_1\ast\Mt_2$. Using the monadicity proviso, Lemma~\ref{splitlemma} ensures that $\BVal(\Mt)=\BVal(\Mt_1\ast\Mt_2)$ and therefore $\der_{\Mt_1\ast\Mt_2}{=}\der_\Mt$, which concludes the proof.\qed
}\\

The previous result allows us to obtain a multiple-conclusion calculus for the logic characterized by a given PNmatrix by simply putting together calculi for simpler logics, under the appropriate provisos. For instance, taking advantage of techniques such as those developed in~\cite{SYNTH}, one can even obtain incrementally analytic axiomatizations for the logic of a given monadic PNmatrix.\smallskip

It is worth noting that the result of Theorem~\ref{splitmult} has some simple immediate consequences. Indeed, it directly applies to any PNmatrix having no more than one designated value and one undesignated value, such as the Boolean-like PNmatrices in Examples~\ref{classical1} and~\ref{platypuses}, by simply using the separator $p$ (this actually shows that one can provide axiomatizations for each connective separately, not just for fragments of classical logic, but also covering less common two-valued non-deterministic connectives). Let us look at another interesting example.

\begin{example}\label{strongk}(Kleene's strong three-valued logic, revisited).\\
Recall the multiple-conclusion version of the implication-free fragment of Kleene's strong three-valued logic $\tuple{\Sigma,\der_\KSMt}$ from Example~\ref{kleene1}, defined on the signature $\Sigma$ containing the connectives $\wedge,\vee,\neg$. The four-valued $\Sigma$-Pmatrix $\KSMt=\tuple{\{0,a,b,1\},\{b,1\},\cdot_\KSMt}$ is monadic, using $p,\neg p$ as its separators. Indeed, $p$ separates the designated values $b,1$ from the undesignated values $0,a$, and $\neg p$ separates $b$ from $1$, and also $0$ from $a$. 

Consider the splitting corresponding to the signatures $\Sigma=\Sigma_1\cup\Sigma_2$ with $\Sigma_1$ containing $\wedge,\neg$ and $\Sigma_2$ containing $\vee,\neg$, and let $\KSMt_1$ and $\KSMt_2$ be the corresponding reducts of $\KSMt$. 
Since $\neg\in\Sigma_1\cap\Sigma_2$, Theorem~\ref{splitmult} tells us that we can obtain a calculus for $\tuple{\Sigma,\der_\KSMt}$ by simply joining calculi for the fragments $\tuple{\Sigma_1,\der_{\KSMt_1}}$ and $\tuple{\Sigma_2,\der_{\KSMt_2}}$. 

It is worth noting that in the multiple-conclusion calculus for $\der_\KSMt$ put forth in Example~\ref{kleene1}, the rules where $\vee$ does not appear constitute a calculus for $\der_{\KSMt_1}$, and the rules where $\wedge$ does not appear constitute a calculus for $\der_{\KSMt_2}$.
\hfill$\triangle$
\end{example}

In the following example we will see that monadicity in the shared language is only a sufficient condition for splitting axiomatizations, but that it really plays an important role in the result.

\begin{example}\label{luka}(Three-valued  \L ukasiewicz logic).\\
Expanding from Example~\ref{lk}, let us now consider the signature $\Sigma$ with $\Sigma^{(1)}=\{\neg,\nabla\}$, $\Sigma^{(2)}=\{\to\}$ and $\Sigma^{{n}}=\emptyset$ for $n\notin \{1,2\}$. The three-valued matrix of  \L ukasiewicz for these connectives is $\mathbb{L}=\tuple{\{0,\frac{1}{2},1\},\{1\},\cdot_\mathbb{L}}$ as defined below.

 \begin{center}
 \begin{tabular}{c | c |c c}
& $\neg_{\mathbb{L}}$& $\nabla_{\mathbb{L}}$ \\
\hline
$0$& $1$ & $0$\\
$\frac{1}{2}$ & $\frac{1}{2}$ & $1$\\
$1$ & $0$ & $1$
\end{tabular}
\qquad
 \begin{tabular}{c | c c c}
$\to_{\mathbb{L}}$& $0$ & $\frac{1}{2}$ & $1$\\
\hline
$0$& $1$ & $1$ & $1$\\
$\frac{1}{2}$ & $\frac{1}{2}$ & $1$ & $1$ \\
$1$ & $0$ & $\frac{1}{2}$ & $1$
\end{tabular}  
\end{center}

Besides the familiar connectives of negation and implication, $\nabla$ 
{  
 is a possibility operator that can be traced back to \L ukasiewicz and Tarski (see \cite{PolishLogic}).
 }
Note that $\nabla A$ is definable as $\neg A\to A$, and thus $\nabla p\der_{\mathbb{L}}\neg p \to p$.

The matrix $\mathbb{L}$ is monadic, for instance using $p,\neg p$ as its separators, but also alternatively using 
$p,\nabla p$, as both $\neg p$ and $\nabla p$ separate $0$ from $\frac{1}{2}$. However, one cannot separate these values using only implication. This said, we could of course apply Theorem~\ref{splitmult} to splittings of signatures sharing $\neg$, or $\nabla$, or both, as in the previous example. Instead, let us look at some more informative cases.\smallskip

{\bf(1)} Let us first consider the splitting corresponding to the signatures $\Sigma=\Sigma_1\cup\Sigma_2$ with $\Sigma_1$ containing $\neg,\to$ and $\Sigma_2$ containing $\nabla$, and let $\mathbb{L}_1$ and $\mathbb{L}_2$ be the corresponding reducts of $\mathbb{L}$. 
Since $\neg,\nabla\notin\Sigma_1\cap\Sigma_2=\emptyset$, Theorem~\ref{splitmult} cannot guarantee that we can obtain a calculus for $\tuple{\Sigma,\der_\mathbb{L}}$ by simply joining calculi for the fragments $\tuple{\Sigma_1,\der_{\mathbb{L}_1}}$ and $\tuple{\Sigma_2,\der_{\mathbb{L}_2}}$. Indeed, this can never be the case as $\tuple{\Sigma,\der_{\mathbb{L}_1\ast\mathbb{L}_2}}$ is strictly weaker that the  \L ukasiewicz logic $\tuple{\Sigma,\der_{\mathbb{L}}}$. To see this,  note that the strict product 
$\mathbb{L}_1\ast\mathbb{L}_2=\tuple{\{00,0\frac{1}{2},\frac{1}{2}0,\frac{1}{2}\frac{1}{2},11\},\{11\},\cdot_\ast}$ is defined by the tables below

 \begin{center}
  \begin{tabular}{c | c |c ccc}
&$\neg_{\ast}$& $\nabla_{\ast}$\\
\hline
$00$& $11$ & $00,\frac{1}{2}0$\\
$0\frac{1}{2}$ & $11$& $11$ \\ 
$\frac{1}{2}0$ & $\frac{1}{2}0,\frac{1}{2}\frac{1}{2}$ & $00,\frac{1}{2}0$\\ 
$\frac{1}{2}\frac{1}{2}$ & $\frac{1}{2}0,\frac{1}{2}\frac{1}{2}$ & $11$\\ 
$11$ & $00,0\frac{1}{2}$ & $11$ 
\end{tabular}
\qquad
 \begin{tabular}{c | c c ccc}
$\to_{\ast}$& $00$ & $0\frac{1}{2}$& $\frac{1}{2}0$& $\frac{1}{2}\frac{1}{2}$ & $11$\\
\hline
$00$& $11$ & $11$ & $11$& $11$ & $11$\\
$0\frac{1}{2}$ & $11$ & $11$ & $11$& $11$ & $11$\\ 
$\frac{1}{2}0$ & $\frac{1}{2}0,\frac{1}{2}\frac{1}{2}$ & $\frac{1}{2}0,\frac{1}{2}\frac{1}{2}$ & $11$& $11$ & $11$\\ 
$\frac{1}{2}\frac{1}{2}$ & $\frac{1}{2}0,\frac{1}{2}\frac{1}{2}$ & $\frac{1}{2}0,\frac{1}{2}\frac{1}{2}$ & $11$& $11$ & $11$\\ 
$11$ & $00,0\frac{1}{2}$ & $00,0\frac{1}{2}$ & $\frac{1}{2}0,\frac{1}{2}\frac{1}{2}$& $\frac{1}{2}0,\frac{1}{2}\frac{1}{2}$ & $11$

\end{tabular}
 \end{center}
and that $\nabla p\not\der_{\mathbb{L}_1\ast\mathbb{L}_2}\neg p \to p$, as witnessed by a valuation $v\in\Val(\mathbb{L}_1\ast\mathbb{L}_2)$ with $v(p)=0\frac{1}{2}$, $v(\neg p)=v(\nabla p)=11$, which necessarily has $v(\neg p \to p)\in\{00,0\frac{1}{2}\}$.

This example shows that although the given PNmatrix, in this case $\mathbb{L}$, is both $\Sigma_1$-monadic and $\Sigma_2$-monadic, the splitting may behave badly precisely because $\mathbb{L}$ is not $\Sigma_1\cap\Sigma_2$-monadic.\smallskip

{\bf(2)} Let us now consider the splitting corresponding to the signatures $\Sigma=\Sigma_1\cup\Sigma_2$ with $\Sigma_1$ containing $\neg,\to$ and $\Sigma_2$ containing $\nabla,\to$, and let $\mathbb{L}_1$ and $\mathbb{L}_2$ be the corresponding reducts of $\mathbb{L}$. 
Again, since $\neg,\nabla\notin\Sigma_1\cap\Sigma_2=\{\to\}$, Theorem~\ref{splitmult} cannot guarantee that we can obtain a calculus for $\tuple{\Sigma,\der_\mathbb{L}}$ by joining calculi for the fragments $\tuple{\Sigma_1,\der_{\mathbb{L}_1}}$ and $\tuple{\Sigma_2,\der_{\mathbb{L}_2}}$. However, in this particular case, it turns out that one can safely join split axiomatizations, simply because 
$\tuple{\Sigma,\der_{\mathbb{L}_1\ast\mathbb{L}_2}}=\tuple{\Sigma,\der_{\mathbb{L}}}$. To see this,  note that the strict product 
$\mathbb{L}_1\ast\mathbb{L}_2=\tuple{\{00,0\frac{1}{2},\frac{1}{2}0,\frac{1}{2}\frac{1}{2},11\},\{11\},\cdot_\ast}$ is 
similar to the one above, but with implication interpreted now as in the table below.

 \begin{center}
  \begin{tabular}{c | c c ccc}
$\to_{\ast}$& $00$ & $0\frac{1}{2}$& $\frac{1}{2}0$& $\frac{1}{2}\frac{1}{2}$ & $11$\\
\hline
$00$& $11$ & $11$ & $11$& $11$ & $11$\\
$0\frac{1}{2}$ & $\emptyset$ & $11$ & $\emptyset$& $11$ & $11$\\ 
$\frac{1}{2}0$ & $\emptyset$ & $\emptyset$ & $11$& $11$ & $11$\\ 
$\frac{1}{2}\frac{1}{2}$ & $\frac{1}{2}\frac{1}{2}$ & $\emptyset$ & $\emptyset$& $11$ & $11$\\ 
$11$ & $00$ & $0\frac{1}{2}$ & $\frac{1}{2}0$& $\frac{1}{2}\frac{1}{2}$ & $11$

\end{tabular}
 \end{center}

Clearly, all bivaluations in $\BVal(\mathbb{L})$ are also in $\BVal(\mathbb{L}_1\ast\mathbb{L}_2)$, as the tables of the total component $\{00,\frac{1}{2}\frac{1}{2},11\}$, depicted below, are simple renamings of the truth-tables of $\mathbb{L}$.
 
  \begin{center}
 \begin{tabular}{c | c |c c}
& $\neg_{\mathbb{L}}$& $\nabla_{\mathbb{L}}$ \\
\hline
$00$& $11$ & $00$\\
$\frac{1}{2}\frac{1}{2}$ & $\frac{1}{2}\frac{1}{2}$ & $11$\\
$11$ & $00$ & $11$
\end{tabular}
\qquad
 \begin{tabular}{c | c c c}
$\to_{\mathbb{L}}$& $00$ & $\frac{1}{2}\frac{1}{2}$ & $11$\\
\hline
$00$& $11$ & $11$ & $11$\\
$\frac{1}{2}\frac{1}{2}$ & $\frac{1}{2}\frac{1}{2}$ & $11$ & $11$ \\
$11$ & $00$ & $\frac{1}{2}\frac{1}{2}$ & $11$
\end{tabular}  
\end{center}

This example shows that the requirement that the given PNmatrix is $\Sigma_1\cap\Sigma_2$-monadic is sufficient, but not necessary, for the splitting axiomatization result to follow.
\hfill $\triangle$
\end{example}

Although obtaining axiomatizations for single-conclusion logics is known to be substantially harder, our line of reasoning may still apply, as long as we further demand saturation.

\begin{theorem}\label{splitsing}
If $R_i$ is an axiomatization of $\tuple{\Sigma_i,\vdash_{\Mt_i}}$, for each $i\in\{1,2\}$, then $R_1\cup R_2$ is an axiomatization of $\tuple{\Sigma,\vdash_\Mt}$, 
provided that $\Mt$ is $(\Sigma_1\cap\Sigma_2)$-monadic and both $\tuple{\Sigma_i,{\Mt_i}}$ are saturated.

\end{theorem}
\proof{
Using Theorem~\ref{saturatedfibring} and Proposition~\ref{fibringcalculi}, we can conclude that  
$\tuple{\Sigma,\vdash_{\Mt_1\ast\Mt_2}}=\tuple{\Sigma_1,\vdash_{\Mt_1}}\bullet\tuple{\Sigma_2,\vdash_{\Mt_2}}=\tuple{\Sigma_1,\vdash_{R_1}}\bullet\tuple{\Sigma_2,\vdash_{R_2}}=\tuple{\Sigma,\vdash_{R_1\cup R_2}}$, that is, $R_1\cup R_2$ is an axiomatization of the single-conclusion logic characterized by the strict product $\Mt_1\ast\Mt_2$. Using the monadicity proviso, now, Lemma~\ref{splitlemma} ensures that $\BVal(\Mt)=\BVal(\Mt_1\ast\Mt_2)$ and therefore $\vdash_{\Mt_1\ast\Mt_2}{=}\vdash_\Mt$, which concludes the proof.\qed
}\\

In Theorem~\ref{splitsing}, we should emphasize that if $\Mt$ is total then it is enough to require that $\tuple{\Sigma,\Mt}$ itself is saturated, as that will guarantee the saturation of both 
$\tuple{\Sigma_i,{\Mt_i}}$.
Naturally, the saturation proviso makes the result of the theorem harder to use, but these difficulties are in line with the results of~\cite{softcomp} regarding the (im)possibility of splitting single-conclusion axiomatizations of classical logic.

\begin{example}\label{classicalfin}(Classical logic can hardly be split).\\
Recall Example~\ref{classical1}, and the substantial differences between combined fragments of classical logic in the single-conclusion (Example~\ref{classical2}) and multiple-conclusion (Example~\ref{multi2}) settings. Expectedly, these differences impact the way  axiomatizations of classical logic may be split. Indeed, while such splitting is universally possible in the multiple-conclusion scenario, as we saw above, it becomes a rarity in the single-conclusion case.

We have seen that classical matrices are always monadic, but now Theorem~\ref{splitsing} further demands saturation. It turns out that the two-valued interpretation of most interesting classical connectives is not saturated, notably for fragments containing $\neg$, $\to$, or $\vee$, as we have seen previously, which agrees with the failure of split axiomatizations as already shown in Example~\ref{classical2}. 

Connectives whose interpretation is saturated include, however, $\wedge$, $\top$, and $\bot$. Hence, Theorem~\ref{splitsing} tells us that fragments of classical logic corresponding to signatures $\Sigma\subseteq\Sigma_\mathsf{cls}^{\wedge,\top,\bot}$ can be axiomatized by putting together axiomatizations for each of the connectives. As seen in Example~\ref{classical2}, joining single-conclusion calculi for $\vdash^\wedge_\mathsf{cls}$ and $\vdash^\top_\mathsf{cls}$ yields a calculus for $\vdash^{\wedge,\top}_\mathsf{cls}$.
\hfill $\triangle$
\end{example}

We next analyze another interesting example.

\begin{example}\label{processor2}(Axiomatizing information sources).\\
Recall from Example~\ref{processor}  the logic of information sources defined, over the signature  $\Sigma_S=\Sigma_\mathsf{cls}^{\wedge,\vee,\neg}$ containing the connectives $\wedge,\vee,\neg$, by the four-valued $\Sigma_S$-Nmatrix 
 $\mathbb{S} = \tuple{\{f,\bot,\top,t\},\{\top,t\},\cdot_{\mathbb{S}}}$. It is easy to see that $\mathbb{S}$ is $\{\neg\}$-monadic, with separators $p,\neg p$. 
 
According to Theorem~\ref{splitmult}, the observations above  imply that a multiple-conclusion calculus for  $\tuple{\Sigma_S,\der_{\mathbb{S}}}$ can be obtained by joining multiple-conclusion calculi for fragments based on splitting signatures $\Sigma_1,\Sigma_2\subseteq\Sigma_S$ such that $\Sigma_1\cup\Sigma_2=\Sigma_S$ as long as $\neg\in\Sigma_1\cap\Sigma_2$.
 
Further, the reduct Nmatrices $\mathbb{S}_i$ are saturated. Given a consistent theory $\Gamma$ of $\vdash_{\mathbb{S}_i}$, it is easy to see that $v^{-1}(\{\top,t\})=\Gamma$ where $v\in\Val(\mathbb{S}_i)$ is defined as follows.

$$v(A)=
\begin{cases}
 \top & \mbox{ if }A,\neg A\in\Gamma \\
 t & \mbox{ if }A\in\Gamma, \neg A\notin\Gamma \\ 
 f & \mbox{ if }A\notin\Gamma, \neg A\in\Gamma \\ 
 \bot & \mbox{ if }A,\neg A\notin\Gamma
 \end{cases}.
$$

Actually, the saturation of $\mathbb{S}$ (and of the $\mathbb{S}_i$ reducts) can also be seen as a consequence of the fact that $\der_\mathbb{S}$ is axiomatized by a single-conclusion calculus.

Hence, according to Theorem~\ref{splitsing}, we know also that a single-conclusion calculus for  $\tuple{\Sigma_S,\vdash_{\mathbb{S}}}$ can be obtained by joining single-conclusion calculi for fragments based on splitting signatures $\Sigma_1,\Sigma_2\subseteq\Sigma_S$ such that $\Sigma_1\cup\Sigma_2=\Sigma_S$ as long as $\neg\in\Sigma_1\cap\Sigma_2$. This is apparent in the calculus put forth in Example~\ref{processor}, if we consider $\neg,\wedge\in\Sigma_1$ and $\neg,\vee\in\Sigma_2$.
\hfill$\triangle$
\end{example}

\subsection{Decidability preservation}

Transfer theorems have always been a main drive of the research in combining logics. Decidability is certainly one of the most desirable properties a logic should have, opening the way for the development of tool support for logical reasoning. There are some known results~\cite{mconiglio:acs:css:10,igpl} regarding the particular case of disjoint combinations, but it is worth looking carefully at the semantic characterizations developed in the previous section, and analyzing their contribution to decidability preservation in general. That is, when given two decidable logics, under which conditions can we guarantee that their combination is still decidable? \smallskip

\subsubsection{Deciding multiple-conclusion combined logics}

Let us first look at the decision problem for multiple-conclusion logics. We will say that a multiple-conclusion logic $\tuple{\Sigma,\der}$ is \emph{decidable} if there exists an algorithm ${\tt D}$, which terminates when given any finite sets $\Gamma,\Delta\subseteq L_\Sigma(P)$ as input, and outputs  ${\tt D}(\Gamma,\Delta)=\textrm{yes}$ if $\Gamma\der\Delta$, and ${\tt D}(\Gamma,\Delta)=\textrm{no}$ if $\Gamma\not\der\Delta$. According to this definition it is clear that one is actually deciding the compact version $\tuple{\Sigma,\der_{\textrm{fin}}}$ of the logic, and hence we will henceforth assume, without loss of generality, that the logic at hand is compact.

Of course, any logic characterized by a finite PNmatrix is decidable~\cite{Baaz2013}. This case covers the combination of any two logics when they are each characterized by a finite PNmatrix, given that the strict product operation preserves finiteness. But we can go beyond the finite-valued case.\smallskip

Corollary~\ref{multiplechar} is quite appealing, and mathematically clean, but a decision procedure based on it would require (potentially) running through all partitions of the set of all formulas. As the similarity with cut for sets is striking,  one may try to obtain a more usable version inspired by cut for a suitable finite set, somehow related to the input, which we will dub \emph{context}. 

In general, we demand a \emph{context function} $\ctx:\wp(L_{\Sigma_1\cup\Sigma_2}(P))\to\wp(L_{\Sigma_1\cup\Sigma_2}(P))$ such that $\Omega\subseteq\ctx(\Omega)$. 
Aiming at decidability preservation, of course, we will further require that $\ctx(\Omega)$ is finite for finite $\Omega\subseteq L_{\Sigma_1\cup\Sigma_2}(P)$.\smallskip

Let $\calB_1\subseteq\BVal(\Sigma_1)$ and $\calB_2\subseteq\BVal(\Sigma_2)$ be sets of bivaluations closed under substitutions. 
We say that $\calB_1,\calB_2$ are \emph{$\ctx$-extensible} when, for any finite set 
$\Omega\subseteq L_{\Sigma_1\cup\Sigma_2}(P)$,
if there exist $b_1\in\calB_1^{\Sigma_1\cup\Sigma_2}$ and  $b_2\in\calB_2^{\Sigma_1\cup\Sigma_2}$ such that $b_1(A)=b_2(A)$ for every $A\in\ctx(\Omega)$, then there must exist $b\in\calB_{12}^{\mult}=\calB_1^{\Sigma_1\cup\Sigma_2}\cap\calB_1^{\Sigma_1\cup\Sigma_2}$ such that $b(A)=b_1(A)=b_2(A)$ for every $A\in\Omega$.

We say that multiple-conclusion logics $\tuple{\Sigma_1,\der_1}$, $\tuple{\Sigma_2,\der_2}$ are \emph{$\ctx$-extensible} when $\calB_1$ and $\calB_2$ are the sets of bivaluations characterizing them, that is, ${\der_1}={\der_{\calB_1}}$ and ${\der_2}={\der_{\calB_2}}$, and $\calB_1,\calB_2$ are themselves $\ctx$-extensible. 

This definition has a more abstract alternative characterization.

\begin{lemma}\label{ctxformultiple}
Let $\tuple{\Sigma_1,\der_1}$, $\tuple{\Sigma_2,\der_2}$ be multiple-conclusion logics, $\tuple{\Sigma_1\cup\Sigma_2,\der_{12}}=\tuple{\Sigma_1,\der_1}\bullet\tuple{\Sigma_2,\der_2}$ their combination, and $\ctx$ a context function. The following are equivalent:
\begin{itemize}
\item $\tuple{\Sigma_1,\der_1}$, $\tuple{\Sigma_2,\der_2}$ are \emph{$\ctx$-extensible};
\item given any partition $\tuple{\underline{\Omega},\overline{\Omega}}$ of $\ctx(\Omega)$ for finite $\Omega$,  if 
$\underline{\Omega}\not\der_1^{\Sigma_1\cup\Sigma_2}\overline{\Omega}$ and $\underline{\Omega}\not\der_2^{\Sigma_1\cup\Sigma_2}\overline{\Omega}$ then  
$\Omega\cap\underline{\Omega}\not\der_{12}\Omega\cap\overline{\Omega}$.
\end{itemize}
\end{lemma}
\proof{The result is immediate, taking into account Proposition~\ref{bivmultiple} and the fact that, for each $k$, the only set of bivaluations characterizing 
$\tuple{\Sigma_k,\der_k}$ is $\calB_k=\{b\in\BVal(\Sigma_k):b^{-1}(1)\not\der_k b^{-1}(0)\}$.\qed}\\

We can now obtain a more decidability-friendly version of Corollary~\ref{multiplechar}.

\begin{lemma}\label{multipledecisionlemma}
Let $\tuple{\Sigma_1,\der_1}$, $\tuple{\Sigma_2,\der_2}$ be $\ctx$-extensible multiple-conclusion logics, and consider their combination $\tuple{\Sigma_1\cup\Sigma_2,\der_{12}}=\tuple{\Sigma_1,\der_1}\bullet\tuple{\Sigma_2,\der_2}$. For every finite $\Gamma,\Delta\subseteq L_{\Sigma_1\cup\Sigma_2}(P)$, we have:
$$\Gamma\der_{12}\Delta$$
$$\textrm{if and only if}$$ 
$$\textrm{for each partition }\tuple{\underline{\Omega},\overline{\Omega}}\textrm{ of }\ctx(\Gamma\cup\Delta)\textrm{, there is }k\in\{1,2\}\textrm{ such that }$$
$$
\Gamma\cup\underline{\Omega}\der_k^{\Sigma_1\cup\Sigma_2}\overline{\Omega}\cup\Delta.$$
\end{lemma}
\proof{Using Corollary~\ref{multiplechar}, if $\Gamma\not\der_{12}\Delta$ then there exists a partition $\tuple{\underline{\Omega},\overline{\Omega}}$ of 
$L_{\Sigma_1\cup\Sigma_2}(P)$ such that $\Gamma\cup\underline{\Omega}\not\der_1^{\Sigma_1\cup\Sigma_2}\overline{\Omega}\cup\Delta$ and 
$\Gamma\cup\underline{\Omega}\not\der_2^{\Sigma_1\cup\Sigma_2}\overline{\Omega}\cup\Delta$.
It is easy to see that $\tuple{\underline{\Omega}\cap\ctx(\Gamma\cup\Delta),\overline{\Omega}\cap\ctx(\Gamma\cup\Delta)}$ is a partition of $\ctx(\Gamma\cup\Delta)$, and by dilution, 
$\Gamma\cup(\underline{\Omega}\cap\ctx(\Gamma\cup\Delta))\not\der_1^{\Sigma_1\cup\Sigma_2}(\overline{\Omega}\cap\ctx(\Gamma\cup\Delta))\cup\Delta$ and 
$\Gamma\cup(\underline{\Omega}\cap\ctx(\Gamma\cup\Delta))\not\der_2^{\Sigma_1\cup\Sigma_2}(\overline{\Omega}\cap\ctx(\Gamma\cup\Delta))\cup\Delta$.

Conversely, if there is a partition $\tuple{\underline{\Omega},\overline{\Omega}}$ of $\ctx(\Gamma\cup\Delta)$ such that $\Gamma\cup\underline{\Omega}\not\der_1^{\Sigma_1\cup\Sigma_2}\overline{\Omega}\cup\Delta$ and $\Gamma\cup\underline{\Omega}\not\der_2^{\Sigma_1\cup\Sigma_2}\overline{\Omega}\cup\Delta$ then, directly from $\ctx$-extensibility and Lemma~\ref{ctxformultiple}, we can conclude that $(\Gamma\cup\Delta)\cap(\Gamma\cup\underline{\Omega})\not\der_{12}^{\Sigma_1\cup\Sigma_2}(\Gamma\cup\Delta)\cap(\overline{\Omega}\cup\Delta)$, or equivalently, since monotonicity implies that $\Gamma\subseteq\underline{\Omega}$ and $\Delta\subseteq\overline{\Omega}$, that $\Gamma\not\der_{12}\Delta$.\qed}\\

We now apply these ideas toward decidability preservation we assume that the context function $\ctx$ is computable.

\begin{theorem}\label{multipledecision}
Let $\tuple{\Sigma_1,\der_1},\tuple{\Sigma_2,\der_2}$ be $\ctx$-extensible logics.
If $\ctx$ is computable and $\tuple{\Sigma_1,\der_1},\tuple{\Sigma_2,\der_2}$ are both decidable then their combination 
$\tuple{\Sigma_1\cup\Sigma_2,\der_{12}}$ is also decidable.
\end{theorem}
\proof{Let ${\tt D_1,D_2}$ be algorithms deciding $\tuple{\Sigma_1,\der_1},\tuple{\Sigma_2,\der_2}$, respectively.
To decide $\tuple{\Sigma_1\cup\Sigma_2,\der_{12}}$ consider the following non-deterministic algorithm ${\tt D}$.
$$\begin{array}{lll}
{\tt D}&:&{\tt input}\;\Gamma,\Delta\\
&& {\tt compute\;\;}\Omega=\ctx(\Gamma\cup\Delta)\\
&& {\tt non-deterministically\; choose\; partition\;}\tuple{\underline{\Omega},\overline{\Omega}}\;{\tt of\;}\Omega\\
&& {\tt if\;}{\tt D_1}(\Gamma\cup\underline{\Omega},\Delta\cup\overline{\Omega})={\tt yes\;or\;}{\tt D_2}(\Gamma\cup\underline{\Omega},\Delta\cup\overline{\Omega})={\tt yes}\;{\tt then}\\
&& \quad {\tt output\; yes}\\
&& {\tt else}\\
&&\quad {\tt output\; no}\\
\end{array}
$$

The correctness of ${\tt D}$ is an immediate consequence of Lemma~\ref{multipledecisionlemma}, noting that the {\tt no} output is always correct, independently of the non-deterministic choice of the partition, whereas the {\tt yes} output is only correct when it holds for all choices. 
\qed}

\begin{corollary}\label{multiplecorol}
Let $\tuple{\Sigma_1,\Mt_1},\tuple{\Sigma_2,\Mt_2}$ be PNmatrices such that their strict product $\tuple{\Sigma_1\cup\Sigma_2,\Mt_1\ast\Mt_2}$ is total.
If $\tuple{\Sigma_1,\der_{\Mt_1}},\tuple{\Sigma_2,\der_{\Mt_2}}$ are both decidable then their combination $\tuple{\Sigma_1\cup\Sigma_2,\der_{\Mt_1\ast\Mt_2}}$ is also decidable.
\end{corollary}
\proof{Given Theorem~\ref{multipledecision}, it is enough to observe that $\BVal(\Mt_1),\BVal(\Mt_2)$ are $\sub$-extensible whenever $\tuple{\Sigma_1\cup\Sigma_2,\Mt_1\ast\Mt_2}$ is total. Note also that $\sub(\Omega)$ is computable. 

Just note that given a set $\Omega\subseteq L_{\Sigma_1\cup\Sigma_2}(P)$ and valuations $v_1\in\Val(\Mt_1^{\Sigma_1\cup\Sigma_2}),v_2\in\Val(\Mt_2^{\Sigma_1\cup\Sigma_2})$ such that $v_1(A)$ is compatible with $v_2(A)$ for every $A\in\sub(\Omega)$, the function $f:\sub(\Omega)\to V_*$ defined by $f(A)=(v_1(A),v_2(A))$ is a prevaluation of $\Mt_1\ast\Mt_2$. Therefore, as $\Mt_1\ast\Mt_2$ is total, $f$ extends to a valuation $v\in\Val(\Mt_1\ast\Mt_2)$. Lemma~\ref{strictvals} thus guarantees the envisaged $\sub$-extensibility property.\qed}\\

It should be noted that $\sub$-extensibility, or totality of the strict product in the case of PNmatrices, is a sufficient condition for preserving decidability, but further research needs to be done in order to find tighter conditions. 

\begin{example}\label{multipledisjointcompact}(Decidability preservation for disjoint combinations).\\
A major result of~\cite{igpl} was the preservation of decidability for disjoint combinations of single-conclusion logics. We will get back to this particular result in the next subsection, but for now we will show that decidability is preserved also by disjoint combination of multiple-conclusion logics. Assume that both $\tuple{\Sigma_1,\der_1},\tuple{\Sigma_2,\der_2}$ are decidable, and $\Sigma_1\cap\Sigma_2=\emptyset$.\smallskip

Some particular cases of this problem have quite easy solutions. For instance, if there exist (total) Nmatrices $\tuple{\Sigma_1,\Mt_1},\tuple{\Sigma_2,\Mt_2}$ such that $\der_{\Mt_1}{=}\der_1$ and $\der_{\Mt_2}{=}\der_2$, and both Nmatrices have designated and undesignated values, it follows from Definition~\ref{strictproduct} that $\tuple{\Sigma_1\cup\Sigma_2,M_1\ast M_2}$ is also a (total) Nmatrix, and the combined logic $\tuple{\Sigma_1\cup\Sigma_2,\der_1\bullet \der_2}$ is decidable as a consequence of Corollary~\ref{multiplecorol}. 
Another particularly straightforward case happens when the truth-values in either of the two Nmatrices are all designated (or all undesignated), as the corresponding component logic will be decidable, in a trivial way, and the same applies to the resulting combined logic. 
In order to tackle the general case, we can use Theorem~\ref{multipledecision}.\smallskip

In order to prove that $\tuple{\Sigma_1\cup\Sigma_2,\der_1\bullet \der_2}$ is decidable
using Theorem~\ref{multipledecision}, it suffices to show that $\tuple{\Sigma_1,\der_1},\tuple{\Sigma_2,\der_2}$ are $\ctx$-extensible for a suitable 
computable context function. Let $\underline{X}$ be a finite theorem-set of either $\der_1$ or $\der_2$, that is, $\emptyset\der_{i}\underline{X}$ for some $i\in\{1,2\}$, if such a set exists. In that case, if $\cB_i$ is the set of bivaluations characterizing $\tuple{\Sigma_i,\der_i}$ and $b\in\cB_i$ then it is clear that $1\in b(\underline{X})$. When none of the component logics has a finite theorem-set then $\underline{X}=\emptyset$, and compactness implies that ${\tt 0}\in\cB_i$. Symmetrically,  let $\overline{X}$ be a finite anti-theorem-set of either $\der_1$ or $\der_2$, that is, $\overline{X}\der_{i}\emptyset$ for some $i\in\{1,2\}$, if such a set exists. In that case, if $b\in\cB_i$ then it is clear that $0\in b(\underline{X})$. As before, $\overline{X}=\emptyset$ if none of the component logics has a finite anti-theorem-set, in which case compactness implies that ${\tt 1}\in\cB_i$. 

We consider $\mathsf{ctx}(\Omega)=\sub(\Omega\cup\underline{X}\cup\overline{X})$, which can clearly be computed.
Suppose now that $\cB_1,\cB_2$ are the sets of bivaluations (closed under substitutions) characterizing the component logics $\tuple{\Sigma_1,\der_1},\tuple{\Sigma_2,\der_2}$, 
and that $b_1\in\cB_1^{\Sigma_1\cup\Sigma_2},b_2\in\cB_2^{\Sigma_1\cup\Sigma_2}$ agree in $\ctx(\Omega)$ for some finite set $\Omega\subseteq L_{\Sigma_1\cup\Sigma_2}(P)$, that is $b_1(A)=b_2(A)$ for every $A\in\ctx(\Omega)$.

If $b_1(A)=b_2(A)=0$ for every $A\in\ctx(\Omega)$ then none of the component logics has a finite theorem-set, $\underline{X}=\emptyset$,  and hence ${\tt 0}\in\calB_{12}^{\mult}=\calB_1^{\Sigma_1\cup\Sigma_2}\cap\calB_1^{\Sigma_1\cup\Sigma_2}$. On the other hand, if  $b_1(A)=b_2(A)=1$ for every $A\in\ctx(\Omega)$ then none of the component logics has a finite anti-theorem-set,  $\overline{X}=\emptyset$, and hence ${\tt 1}\in\calB_{12}^{\mult}$. Thus, we proceed knowing that we can fix formulas $F_0,F_1\in\ctx(\Omega)$ such that $b_1(F_0)=b_2(F_0)=0$ and $b_1(F_1)=b_2(F_1)=1$. 

First we modify $b_1,b_2$ so that they also agree on $P\setminus\var(\ctx(\Omega))$ (for simplicity, we chose to evaluate them all to $0$). Consider the substitution $\sigma:P\to L_{\Sigma_1\cup\Sigma_2}(P)$ such that
$$\sigma(p)=
\begin{cases}
p &\mbox{if } p\in \var(\ctx(\Omega))\\
F_0 & \mbox{otherwise}
\end{cases}
$$
and let $b_1^0=b_1\circ\sigma$ and $b_2^0=b_2\circ\sigma$. Easily, $b_1^0\in\cB_1^{\Sigma_1\cup\Sigma_2},b_2^0\in\cB_2^{\Sigma_1\cup\Sigma_2}$ extend $b_1,b_2$ on $\ctx(\Omega)$, and further agree on all formulas in $\Omega^0=P\cup\ctx(\Omega)$. 

Recursively, let $\Omega^{k+1}=\Omega^k\cup\{c(A_1,\dots,A_n)\in L_{\Sigma_1\cup\Sigma_2}(P):A_1,\dots,A_n\in\Omega_k\}$ and obtain bivaluations $b_1^{k+1}\in\cB_1^{\Sigma_1\cup\Sigma_2},b_2^{k+1}\in\cB_2^{\Sigma_1\cup\Sigma_2}$ extending the previous on $\Omega^k$, and further agreeing on $\Omega^{k+1}$. For each formula $A\in\Omega^{k+1}\setminus\Omega^k$ with $\head(A)\notin\Sigma_i$ we evaluate $b_{3-i}^k(A)$ and modify the skeleton variable $\skel_i(A)$ accordingly when building $b_i^{k+1}$. Hence, consider for each $i\in\{1,2\}$ the substitution $\sigma_i^{k}:P\to L_{\Sigma_i}(P)$ such that
$$\sigma_i^{k}(p)=
\begin{cases}
\skel_i(F_{x}) &\mbox{if } p=\skel_i(A), A\in\Omega^{k+1}\setminus\Omega^k, \head(A)\notin\Sigma_i,b_{3-i}^k(A)=x\\
p & \mbox{otherwise}
\end{cases}
$$
and set $b_i^{k+1}=b_i^k{\circ}\unskel_i{\circ}\sigma_i^{k}{\circ}\skel_i$. If $A\in\Omega^k$, or $A\in\Omega^{k+1}\setminus\Omega^k$ and $\head(A)\in \Sigma_i$, 
then $\skel_i(A)^{\sigma_i^{k}}=\skel_i(A)$, and thus $b_i^{k+1}(A)=b_i^{k}(A)$.
Furthermore, if $A\in \Omega^{k+1}\setminus\Omega^k$ and $\head(A)\notin\Sigma_i$ then 
$b_i^{k+1}(A)=b_i^k(\unskel_i(\sigma_i^{k}(\skel_i(A))))=b_i^k(\unskel_i(\skel_i(F_{b_{3-i}^k(A)})))=b_i^k(F_{b_{3-i}^k(A)})=b_{3-i}^k(A)=b_{3-i}^{k+1}(A)$.

Partitioning each $\Omega^k$ in two disjoint parts $\Gamma^k=\{A\in\Omega^k:b_1^k(A)=1\}=\{A\in\Omega^k:b_2^k(A)=1\}$, and $\Delta^k=\{A\in\Omega^k:b_1^k(A)=0\}=\{A\in\Omega^k:b_2^k(A)=0\}$, for $k\in\nats_0$, it is clear that $\Gamma^k\subseteq\Gamma^{k+1}$, $\Delta^k\subseteq\Delta^{k+1}$, and also that $\Gamma=\bigcup_{k\in\nats_0}\Omega^k$ and $\Delta=\bigcup_{k\in\nats_0}\Delta^k$ are a partition of $L_{\Sigma_1\cup\Sigma_2}(P)=\bigcup_{k\in\nats_0}\Omega^k$. Further, each $b_i^k$ shows precisely that $\Gamma^k\not\der_i^{\Sigma_1\cup\Sigma_2}\Delta^k$. The compactness of both $\der_1$ and $\der_2$ then implies that 
$\Gamma\not\der_1^{\Sigma_1\cup\Sigma_2}\Delta$ and $\Gamma\not\der_2^{\Sigma_1\cup\Sigma_2}\Delta$. Therefore, the bivaluation $b:L_{\Sigma_1\cup\Sigma_2}(P)\to\{0,1\}$ such that $b(A)=1$ if $A\in\Gamma$, and $b(A)=0$ if $A\in\Delta$ (resulting as a \emph{limit} of the above sequences of bivaluations) is such that $b\in\cB_1^{\Sigma_1\cup\Sigma_2}\cap\cB_2^{\Sigma_1\cup\Sigma_2}=\calB_{12}^{\mult}$. Of course, $b$ agrees with $b_1,b_2$ on $\ctx(\Omega)\supseteq\Omega$, which concludes the argument.\hfill$\triangle$
\end{example}

We will further illustrate the use of such general results to concrete combined logics in the more familiar case, below, of single-conclusion logics. 
Of course, our general decidability preservation technique using a suitable context function resembles many of the decidability preservation techniques used in concrete examples, such as fusions of modal logics~\cite{gkwz03}, or even beyond in combining equational and first-order theories~\cite{pig74,nelson79,baader,tinelli,rasga:sernadas}.  
A deeper account of the scope of the abstraction we propose is beyond the reach of this paper.\smallskip

\subsubsection{Deciding single-conclusion combined logics}

As before, the single-conclusion case can now be addressed as an application of the same ideas. A single-conclusion logic $\tuple{\Sigma,\vdash}$ is \emph{decidable} if there exists an algorithm ${\tt D}$, which terminates when given any finite set $\Gamma\subseteq L_\Sigma(P)$ and formula $A\in L_\Sigma(P)$ as input, and outputs  ${\tt D}(\Gamma,A)={\tt yes}$ if $\Gamma\vdash A$, and ${\tt D}(\Gamma,A)={\tt no}$ if $\Gamma\not\vdash A$. As before, we will henceforth assume with lost of generality that the logic at hand is compact, as this definition is equivalent to deciding the compact version $\tuple{\Sigma,\vdash_{\textrm{fin}}}$ of the logic.

\smallskip

The following very simple result will help us to apply the ideas used in the multiple-conclusion scenario, to the single-conclusion case.

\begin{proposition}\label{decisions}
The following implications hold.
\begin{enumerate}
\item[(a)]
If a multiple-conclusion logic $\tuple{\Sigma,\der}$ is decidable, then so is its single-conclusion companion $\sing(\tuple{\Sigma,\der})$.
\item[(b)]
If a single-conclusion logic $\tuple{\Sigma,\vdash}$ is decidable, then so is its minimal multiple-conclusion counterpart $\posit(\tuple{\Sigma,\vdash})$.
\end{enumerate}
\end{proposition}
\proof{We consider each property.
\begin{enumerate}
\item[(a)]
Let $\sing(\tuple{\Sigma,\der})=\tuple{\Sigma,\vdash}$, and ${\tt D}$ be an algorithm deciding $\tuple{\Sigma,\der}$.
To decide $\tuple{\Sigma,\vdash}$ consider the following algorithm ${\tt D'}$.
$$\begin{array}{lll}
{\tt D'}&:&{\tt input}\;\Gamma,A\\
&& {\tt output\;}{\tt D}(\Gamma,\{A\})
\end{array}
$$
Clearly, $\Gamma\vdash A$ if and only if $\Gamma\der\{A\}$, and ${\tt D'}$ decides $\tuple{\Sigma,\vdash}$. 
\item[(b)] 
Let $\posit(\tuple{\Sigma,\vdash})=\tuple{\Sigma,\der}$, and ${\tt D}$ be an algorithm deciding $\tuple{\Sigma,\vdash}$.
To decide $\tuple{\Sigma,\der}$ consider the following algorithm ${\tt D'}$.
$$\begin{array}{lll}
{\tt D'}&:&{\tt input}\;\Gamma,\Delta\\
&& {\tt res}={\tt no}\\
&& {\tt for\;each\;}B\in\Delta:\\
&&\quad {\tt if\;}{\tt D}(\Gamma,B)={\tt yes}\;{\tt then}\\
&&\quad\quad {\tt\;res=yes}\\
&& {\tt output\; res}
\end{array}
$$
Clearly, $\Gamma\der\Delta$ if and only if $\Gamma\der B$ for some $B\in\Delta$. Hence, ${\tt D'}$ decides $\tuple{\Sigma,\der}$.
\qed\smallskip
\end{enumerate}
}

Now, to go directly to the results, we say that single-conclusion logics $\tuple{\Sigma_1,\vdash_1}$, $\tuple{\Sigma_2,\vdash_2}$ are \emph{$\ctx$-extensible} when $\calB_1$ and $\calB_2$ are sets of bivaluations characterizing them, that is, ${\vdash_1}={\vdash_{\calB_1}}$ and ${\vdash_2}={\vdash_{\calB_2}}$, and their (uniquely determined) meet-closures $\calB_1^\cap,\calB_2^\cap$ are $\ctx$-extensible.

Again, this definition has a more abstract alternative characterization, based on theories.

\begin{lemma}\label{ctxforsingle}
Let $\tuple{\Sigma_1,\vdash_1}$, $\tuple{\Sigma_2,\vdash_2}$ be single-conclusion logics, $\tuple{\Sigma_1\cup\Sigma_2,\vdash_{12}}=\tuple{\Sigma_1,\vdash_1}\bullet\tuple{\Sigma_2,\vdash_2}$ their combination, and $\ctx$ a context function. The following are equivalent:
\begin{itemize}
\item $\tuple{\Sigma_1,\vdash_1}$, $\tuple{\Sigma_2,\vdash_2}$ are \emph{$\ctx$-extensible};
\item given $\Omega$ and theories $\Delta_k$ of $\tuple{\Sigma_1\cup\Sigma_2,\vdash_k^{\Sigma_1\cup\Sigma_2}}$ for $k\in\{1,2\}$,  if 
$\Delta_1\cap\ctx(\Omega)=\Delta_2\cap\ctx(\Omega)$ then there exists a theory $\Delta$ of $\tuple{\Sigma_1\cup\Sigma_2,\vdash_{12}}$ such that $\Delta\cap\Omega=\Delta_1\cap\Omega=\Delta_2\cap\Omega$.
\end{itemize}
\end{lemma}
\proof{The result is immediate, taking into account Proposition~\ref{bivsingle} and the fact that, for each $k$, the only meet-closed set of bivaluations characterizing 
$\tuple{\Sigma_k,\vdash_k}$ is $\calB^\cap_k=\{b\in\BVal(\Sigma_k):b^{-1}(1) \textrm{ is a theory of } \tuple{\Sigma_k,\vdash_k}\}$.\qed}\\

As before, we obtain a general decidability preservation result.

\begin{theorem}\label{singledecision}
Let $\tuple{\Sigma_1,\vdash_1},\tuple{\Sigma_2,\vdash_2}$ be $\ctx$-extensible logics. If $\ctx$ is computable and $\tuple{\Sigma_1,\vdash_1},\tuple{\Sigma_2,\vdash_2}$ are both decidable then their combination $\tuple{\Sigma_1\cup\Sigma_2,\vdash_{12}}$ is also decidable.
\end{theorem}
\proof{Let $\calB_1\subseteq\BVal(\Sigma_1)$ and $\calB_2\subseteq\BVal(\Sigma_2)$, both closed under substitutions, be such that ${\vdash_1}={\vdash_{\calB_1}}={\vdash_{\calB_1^\cap}}$ and ${\vdash_2}={\vdash_{\calB_2}}={\vdash_{\calB_2^\cap}}$. From Proposition~\ref{decisions}~(b), we know that $\tuple{\Sigma_1,\der_{\calB_1^\cap}}$ and $\tuple{\Sigma_1,\der_{\calB_2^\cap}}$ are both decidable. From the hypothesis of $\ctx$-extensibility of $\tuple{\Sigma_1,\vdash_1},\tuple{\Sigma_2,\vdash_2}$, we know that $\calB_1^\cap,\calB_2^\cap$ are $\ctx$-extensible, and thus also  
$\tuple{\Sigma_1,\der_{\calB_1^\cap}},\tuple{\Sigma_1,\der_{\calB_2^\cap}}$. Then, Theorem~\ref{multipledecision} guarantees that  $\tuple{\Sigma_1,\der_{\calB_1^\cap}}\bullet\tuple{\Sigma_1,\der_{\calB_2^\cap}}$ is decidable.
But we know that $\tuple{\Sigma_1,\vdash_1}\bullet\tuple{\Sigma_2,\vdash_2}=\tuple{\Sigma_1\cup\Sigma_2,\vdash_{12}}=\sing(\tuple{\Sigma_1,\der_{\calB_1^\cap}}\bullet\tuple{\Sigma_1,\der_{\calB_2^\cap}})$, and Proposition~\ref{decisions}~(a) guarantees that the combined logic is decidable.

Putting together the algorithms obtained in the proofs of Theorem~\ref{multipledecision} and Proposition~\ref{decisions}, if ${\tt D_1,D_2}$ are algorithms deciding $\tuple{\Sigma_1,\vdash_1},\tuple{\Sigma_2,\vdash_2}$, respectively, we obtain the following
non-deterministic algorithm ${\tt D}$  for deciding $\tuple{\Sigma_1\cup\Sigma_2,\vdash_{12}}$.
$$\begin{array}{lll}
{\tt D}&:&{\tt input}\;\Gamma,A\\
&& {\tt compute\;\;}\Omega=\ctx(\Gamma\cup\{A\})\\
&& {\tt non-deterministically\; choose\; partition\;}\tuple{\underline{\Omega},\overline{\Omega}}\;{\tt of\;}\Omega\\
&& {\tt res}={\tt no}\\
&& {\tt for\; each\;}B\in(\overline{\Omega}\cup\{A\})\\
&& \quad{\tt if\;}{\tt D_1}(\Gamma\cup\underline{\Omega},B)={\tt yes\; or\;}{\tt D_2}(\Gamma\cup\underline{\Omega},B)={\tt yes}\;{\tt then}\\
&& \quad\quad {\tt res=yes}\\
&&{\tt output\; res}
\end{array}
$$

\vspace*{-6mm}\qed
}\smallskip

Let us illustrate some particular applications of Theorem~\ref{singledecision}.

\begin{corollary}\label{singlecorol}
Let $\tuple{\Sigma_1,\Mt_1},\tuple{\Sigma_2,\Mt_2}$ be saturated PNmatrices such that their strict product $\tuple{\Sigma_1\cup\Sigma_2,\Mt_1\ast\Mt_2}$ is total.
If $\tuple{\Sigma_1,\vdash_{\Mt_1}}$, $\tuple{\Sigma_2,\vdash_{\Mt_2}}$ are both decidable then also their combination $\tuple{\Sigma_1\cup\Sigma_2,\vdash_{\Mt_1\ast\Mt_2}}$ is decidable.
\end{corollary}
\proof{Given Theorem~\ref{singledecision}, it is enough to observe that $\BVal(\Mt_1)^\cap,\BVal(\Mt_2)^\cap$ are $\sub$-extensible whenever $\tuple{\Sigma_1\cup\Sigma_2,\Mt_1\ast\Mt_2}$ is total.

By saturation, Lemma~\ref{theone} guarantees that $\BVal(\Mt_k)^\cap=\BVal(\Mt_k)\cup\{{\tt 1}\}$ for each $k\in\{1,2\}$. Their $\sub$-extensibility follows  from the fact that also $\BVal(\Mt_1),\BVal(\Mt_2)$ are $\sub$-extensible, implied by the hypothesis that the strict product is total, as in the proof of Corollary~\ref{multiplecorol}.
\qed}\\

The case of disjoint combinations is also worth revisiting.

%%%%%%% EXAMPLES

\begin{example}\label{disjointcompact}(Decidability preservation for disjoint combinations, again).\\
We are now able to obtain a simple proof of the the major result of~\cite{igpl}: the preservation of decidability for disjoint combinations of single-conclusion logics. Let $\tuple{\Sigma_1,\vdash_1},\tuple{\Sigma_2,\vdash_2}$ be decidable and 
$\Sigma_1\cap\Sigma_2=\emptyset$. Using Proposition~\ref{decisions} we know that $\posit(\tuple{\Sigma_1,\vdash_1}),\posit(\tuple{\Sigma_2,\vdash_2})$ are decidable as well, and the result in Example~\ref{multipledisjointcompact} tells us that $\posit(\tuple{\Sigma_1,\vdash_1})\bullet\posit(\tuple{\Sigma_2,\vdash_2})$ is also decidable. Finally, note that using again Proposition~\ref{decisions}, we have that 
$\sing(\posit(\tuple{\Sigma_1,\vdash_1})\bullet\posit(\tuple{\Sigma_2,\vdash_2}))$ is decidable, and we know that 
$\sing(\posit(\tuple{\Sigma_1,\vdash_1})\bullet\posit(\tuple{\Sigma_2,\vdash_2}))=\sing(\posit(\tuple{\Sigma_1,\vdash_1}))\bullet\sing(\posit(\tuple{\Sigma_2,\vdash_2}))=\tuple{\Sigma_1,\vdash_1}\bullet\tuple{\Sigma_2,\vdash_2}$.\smallskip

A more direct proof of the same result could instead be obtained from Theorem~\ref{singledecision}, along an argument similar to the one used in Example~\ref{multipledisjointcompact} for the multiple-conclusion setting.\hfill$\triangle$
\end{example}

Ultimately, our results followed simply as applications of the corresponding results in the multiple-conclusion case. In any case, a result analogous to Lemma~\ref{multipledecisionlemma} can still be obtained, which would imply Theorem~\ref{singledecision} and Corollary~\ref{singlecorol}, as well, and which can be understood as a decidability-friendly version of Corollary~\ref{singlechar}.

\begin{lemma}\label{singledecisionlemma}
Let $\tuple{\Sigma_1,\vdash_1}$, $\tuple{\Sigma_2,\vdash_2}$ be $\ctx$-extensible single-conclusion logics, and consider their combination $\tuple{\Sigma_1\cup\Sigma_2,\vdash_{12}}=\tuple{\Sigma_1,\vdash_1}\bullet\tuple{\Sigma_2,\vdash_2}$. For every finite $\Gamma\cup\{A\}\subseteq L_{\Sigma_1\cup\Sigma_2}(P)$, we have:
$$\Gamma\vdash_{12}A$$
$$\textrm{if and only if}$$
$$\textrm{ for each }\Gamma\subseteq\Omega\subseteq \ctx(\Gamma\cup\{A\}),$$
$$\textrm{if }
((\Omega^{\vdash_1^{\Sigma_1\cup\Sigma_2}}\cup\Omega^{\vdash_2^{\Sigma_1\cup\Sigma_2}})\cap\ctx(\Gamma\cup\{A\}))\subseteq\Omega\textrm{ then }A\in\Omega.$$
\end{lemma}
\proof{Using Corollary~\ref{singlechar}, if $\Gamma\not\vdash_{12}A$ then there exists $\Gamma\subseteq\Omega\not\owns A$ such that $\Omega$ is a theory of both 
$\tuple{\Sigma_1\cup\Sigma_2,\vdash_1^{\Sigma_1\cup\Sigma_2}}$ and 
$\tuple{\Sigma_1\cup\Sigma_2,\vdash_2^{\Sigma_1\cup\Sigma_2}}$. Then, one has ${((\Omega\cap\ctx(\Gamma\cup\{A\}))^{\vdash_k^{\Sigma_1\cup\Sigma_2}}}{\cap}\ctx(\Gamma\cup\{A\}))\subseteq
%{(\Omega^{\vdash_k^{\Sigma_1\cup\Sigma_2}})\cap\ctx(\Gamma\cup\{A\})}=
{\Omega\cap\ctx(\Gamma\cup\{A\})}$ for each $k\in\{1,2\}$, 
and $\Gamma\subseteq{\Omega\cap\ctx(\Gamma\cup\{A\})}\not\owns A$.

Conversely, if there is  $\Gamma\subseteq\Omega\subseteq \ctx(\Gamma\cup\{A\})$ such that $A\notin\Omega$, but with $\Omega^{\vdash_k^{\Sigma_1\cup\Sigma_2}}\cap{\ctx(\Gamma\cup\{A\}})\subseteq\Omega$ for each $k\in\{1,2\}$ then 
it follows that $\Delta_1=\Omega^{\vdash_1^{\Sigma_1\cup\Sigma_2}}$ 
and $\Delta_2=\Omega^{\vdash_2^{\Sigma_1\cup\Sigma_2}}$ are theories such that $\Delta_1\cap{\ctx(\Gamma\cup\{A\}})=\Delta_2\cap{\ctx(\Gamma\cup\{A\}})$. 
Thus, directly from $\ctx$-extensibility and Lemma~\ref{ctxforsingle}, we can conclude that 
there exists a theory $\Delta$ of $\tuple{\Sigma_1\cup\Sigma_2,\vdash_{12}}$ such that 
$\Delta\cap (\Gamma\cup\{A\})=\Delta_1\cap (\Gamma\cup\{A\})=\Delta_2\cap (\Gamma\cup\{A\})$.
It follows that $\Gamma\subseteq\Delta\not\owns A$ and we conclude that $\Gamma\not\vdash_{12}A$.\qed}

\section{Conclusion}\label{sec:conc}

The main contribution of this paper is the definition of a first modular semantics for combined logics, in both the single and multiple-conclusion scenarios. 
It is worth emphasizing again that the analysis of the latter scenario was crucial for the development, due to its tight connection with bivaluations. Of course, bivaluations could be simply seen as partitions of the language, or as theories, but looking at them as semantic functions helps to smoothen the path to many-valued interpretations. Naturally, the adoption of PNmatrices as models was also fundamental, as partiality enables us to deal with possibly conflicting interpretations of shared language, whereas non-determinism makes it straightforward to accomodate language extensions.

Discovering the finiteness-preserving strict-product operation on PNmatrices and its universal property were certainly central to the results obtained. 
Further, it should be said that it provides a solution to the problem at hand which is quite close to the abstract idea of Gabbay's fibring function~\cite{gabbay1996,GAbFib99,ccal:car:jfr:css:04}, and a very natural modular version of Schechter's proposal~\cite{schechter}. Unfortunately, when combining single-conclusion logics, this is simply not enough. The solution we found uses the $\omega$-power operation for saturation purposes, but at the expense of losing finite-valuedness. This path hinders the easy application of these techniques in practice for single-conclusion logics, in general, also because checking saturation does not seem to be a trivial task. Nevertheless, we have seen that $\omega$-powers may ultimately be too radical a solution, namely as finite powers are sometimes sufficient. Clearly, a better understanding of saturation is necessary, including its connection with admissibility of rules, as studied for instance in~\cite{IEMHOFF,CINTULA2010162,bezhanishvili_bezhanishvili_iemhoff_2016,BezhGabGhiJib}.

At this point we should remark that the fact that we only consider logics defined by a single PNmatrix, instead of a collection of PNmatrices, is a simplifying assumption with almost no loss of generality, as we have shown that partiality allows for summing collections of PNmatrices whenever one has at least one connective with more than one place. \smallskip

The three applications developed taking into account our semantic characterizations are also worth mentioning as valuable contributions. 

First, the construction of a semantics for strengthening a given many-valued logic with additional axioms is interesting in its implications, but just a reinterpretation of a result in~\cite{addax}, where less general but finite semantics for certain well-behaved particular cases are also obtained. Nevertheless, obtaining a denumerable PNmatrix for intuitionistic propositional logic is a worth example of the compressing power of partiality and non-determinism, and deserves to be further explored, particularly in connection with ideas for aproximating logics (see, for instance,~\cite{baazach,dagostino,ghilardi,masacci,finger}).

Secondly, studying conditions under which the problem of obtaining a calculus for a logic may be split into the problem of obtaining axiomatizations for suitable syntactically defined fragments of the logic is quite crucial for a modular understanding of combining logics. A related approach was considered in~\cite{Coniglio07}, but aiming at a disjoint split with additional axioms. On the contrary, the results we obtain here are inline with ideas explored in~\cite{softcomp} concerning axiomatizations of classical logic, and fit well with our running track of research on extracting calculi for PNmatrices in an automated way, by taking advantage of monadicity requirements~\cite{SYNTH,wollic19}.

Third, and last, having a clear semantics for combined logics is a crucial tool for studying their decidability. The very general criteria we obtained already cover previous results regarding disjoint combinations, as we have shown. We believe they are also powerful enough to encompass other results in the literature, like the decidability of fusions of modal logics~\cite{Wolter98,gkwz03}, and even adaptations of Nelson-Oppen-like techniques for deciding certain equational and first-order theories~\cite{pig74,nelson79,baader,tinelli,rasga:sernadas}. A thourough analysis of the complexity of the obtained decision algorithms was intentionally not addressed in this paper, but certainly deserves future attention.

\bibliographystyle{plain}
\bibliography{biblio}

\end{document}